\documentclass[final, reqno]{amsart}
\usepackage{amsmath}
\usepackage{amsfonts}
\usepackage{amssymb}
\usepackage{amsthm}
\title{A Centre-Stable Manifold\\for the Focussing Cubic NLS in $\set R^{1+3}$}
\author{Marius Beceanu}
\thanks{Department of Mathematics, University of Chicago, 5734 S. University Ave., Chicago, IL 60637 ({\tt mbeceanu@uchicago.edu}). This work is part of the author's Ph.\ D.\ thesis at the University of Chicago.}

\newtheorem{theorem}{Theorem}[section]
\newtheorem{lemma}[theorem]{Lemma}
\newtheorem{proposition}[theorem]{Proposition}
\newtheorem{corollary}[theorem]{Corollary}

\newtheorem{assumption}{Assumption}
\newtheorem{definition}{Definition}
\newcommand{\set}{\mathbb}
\newcommand{\dl}{\nabla}
\renewcommand{\frak}{\mathfrak}

\newcommand{\mc}{\mathcal}
\newcommand{\be}{\begin{equation}}
\newcommand{\ee}{\end{equation}}

\newcommand{\ba}{\begin{array}}
\newcommand{\ds}{\displaystyle}
\newcommand{\ea}{\end{array}}
\newcommand{\bpm}{\begin{pmatrix}}
\newcommand{\epm}{\end{pmatrix}}
\newcommand{\lb}{\label}
\DeclareMathOperator{\sgn}{sgn}

\newcommand{\ov}{\overline}
\newcommand{\dd}{{\,}{d}}
\newcommand{\dg}{\dagger}
\usepackage [frenchlinks=true]{hyperref}

\begin{document}
\numberwithin{equation}{section}
\begin{abstract}
Consider the focussing cubic nonlinear Schr\"odinger equation in $\set R^3$:
\be
i\psi_t+\Delta\psi = -|\psi|^2 \psi.
\lb{NLS}
\ee
It admits special solutions of the form $e^{it\alpha}\phi$, where $\phi \in \mc S(\set R^3)$ is a positive ($\phi>0$) solution of
\be
-\Delta \phi + \alpha\phi = \phi^3.
\lb{phi}
\ee
The space of all such solutions, together with those obtained from them by rescaling and applying phase and Galilean coordinate changes, called standing waves, is the $8$-dimensional manifold that consists of functions of the form $e^{i(v \cdot + \Gamma)} \phi(\cdot - y, \alpha)$.

We prove that any solution starting sufficiently close to a standing wave in the $\Sigma = W^{1, 2}(\set R^3) \cap |x|^{-1}L^2(\set R^3)$ norm and situated on a certain codimension-one local Lipschitz manifold exists globally in time and converges to a point on the manifold of standing waves. 

Furthermore, we show that $\mc N$ is invariant under the Hamiltonian flow, locally in time, and is a centre-stable manifold in the sense of Bates, Jones \cite{bates}.

The proof is based on the modulation method introduced by Soffer and Weinstein for the $L^2$-subcritical case and adapted by Schlag to the $L^2$-supercritical case. An important part of the proof is the Keel-Tao endpoint Strichartz estimate in $\set R^3$ for the nonselfadjoint Schr\"odinger operator obtained by linearizing (\ref{NLS}) around a standing wave solution.
\end{abstract}
\maketitle
\pagestyle{myheadings}
\thispagestyle{plain}
\markboth{M. BECEANU}{Focussing Cubic NLS in 3D}

\section{Introduction}
\subsection{Main result}
For a parameter path $\pi = (v_k, D_k, \alpha, \Gamma)$ such that $\|\dot \pi\|_{\infty} + \|\langle t \rangle \dot\pi(t)\|_1 < \infty$, define the nonuniformly moving soliton $W(\pi(t))$ by
\be\begin{aligned}
W(\pi(t))(x) &= e^{i \theta(t, x)} \phi(x-y(t), \alpha(t)) \\
\theta(t, x) &\ds= v(t) x - \int_0^t(|v(s)|^2 - \alpha(s)) \dd s + \gamma(t)\\
y(t) &\ds= 2\int_0^t v(s) \dd s + D(t)\\
\gamma(t) &= \Gamma(t) - (v(t)-v(\infty))D(\infty) + \int_t^{\infty} 2s \dot v(s) v(\infty) \dd s.
\end{aligned}
\ee
\begin{theorem}[Main result]
There exists a local codimension-one Lipschitz manifold $\mc N$ in $\Sigma = H^1 \cap |x|^{-1} L^2$, containing the 8-dimensional manifold of standing waves, such that equation (\ref{NLS}) has a global dispersive solution $\Psi$ if we start with initial data $\Psi(0)$ on the manifold $\mc N$.

Furthermore, the solution depends Lipschitz continuously on the initial data and decomposes into a moving soliton and a dispersive term: $\Psi = W(\pi(t)) + R(t)$, with
\be
\|\dot\pi\|_{\infty} + \|\langle t \rangle \dot\pi(t)\|_1 \leq C \|\Psi(0) - W(\pi(0))\|_{\Sigma}
\ee
and
\be
\|R\|_{L^{\infty}_t L^2_x \cap L^2_t L^6_x \cap \langle t \rangle^{-1/2} L^2_t L^{6+\infty}_x} \leq C \|\Psi(0) - W(\pi(0))\|_{\Sigma}.
\ee
The dispersive term scatters: $R(t) = e^{it\Delta} f_0 + o_{L^2}(1)$, for some $f_0 \in L^2$.

Moreover, for a solution $\Psi$ of initial data $\Psi(0) \in \mc N$, one has that $\Psi(t) \in \Sigma$ for all $t$ and $\psi(t) \in \mc N$ for sufficiently small~$t$.

Finally, $\mc N$ is a centre-stable manifold for this equation in the sense of Bates, Jones \cite{bates}.
\lb{tichie}
\end{theorem}

\subsection{Background}
Consider the focussing nonlinear cubic Schr\"odinger equation (\ref{NLS}). It admits a particular class of solutions of the form $e^{it\alpha}\phi$, where $\phi = \phi(\cdot, \alpha) \in \mc S$, $\phi>0$, are solutions of (\ref{phi}).

These solutions exist for all time and are periodic. Positive, smooth solutions $\phi$ to (\ref{phi}) are called ground states and solutions to (\ref{NLS}) obtained from $e^{it\alpha} \phi$ by Galilean coordinate changes, phase changes, or scaling are called standing waves. All these transformations are symmetries of the equation, subsumed by the following formula:
\be
\mc G(t)(f(x, t)) = e^{i(\Gamma+vx-t|v|^2)} f(\alpha^{1/2}x-2tv-D, \alpha t).
\ee

A natural question is whether standing waves are stable under small perturbations.

From a physical point of view, the NLS equation in $\set R^3$ with cubic nonlinearity and the focussing sign (\ref{phi}) describes, to a first approximation, the self-focussing of optical beams due to the nonlinear increase of the refraction index. As such, the equation appeared for the first time in the physical literature in 1965, in \cite{kelley}.

\subsection{Known stability results in other cases} Concerning the general NLS problem, more results have been obtained in the defocussing case or for $L^2$-subcritical and $L^2$-critical power nonlinearities in the focussing case. A few negative results have been established as well.

Cazenave and Lions \cite{cazenave} and Weinstein \cite{wein1}, \cite{wein2} used the method of modulation to prove the orbital stability of standing waves in the focussing $L^2$-subcritical case. Asymptotic stability results have been first obtained by Soffer, Weinstein \cite{soffer1}, \cite{soffer2}, then by Pillet, Wayne, \cite{pillet}, Buslaev, Perelman \cite{buslaev1}, \cite{buslaev2}, \cite{buslaev3}, Cuccagna \cite{cuc}, Rodnianski, Schlag, Soffer,  \cite{rod2}, \cite{rod3}, Schlag \cite{schlag}, and Krieger, Schlag \cite{krieger}. Grillakis, Shatah, and Strauss \cite{gril1}, \cite{gril2} developed a general theory of stability of solitary waves for Hamiltonian evolution equations, which, when applied to the Schr\"odinger equation, shows the dichotomy between the $L^2$-subcritical and critical or supercritical cases.

If the nonlinearity is $L^2$-critical or supercritical and focussing, negative energy $\langle x\rangle^{-1} H^1$ initial data leads to solutions that blow up in finite time, due to the virial identity (see Glassey \cite{glassey}). For weakening the condition on initial data and for a survey of this topic see \cite{sulem} and \cite{caz2}. Berestycki, Cazenave \cite{bercaz} showed that blow-up can occur for arbitrarily small perturbations of ground states. Recent results concerning the blowup of the critical and supercritical equation include Merle, Raphael \cite{mer2} and Krieger, Schlag \cite{kri3}.

In 1993, Merle \cite{merle} showed in the $L^2$-critical case the existence of a minimal blow-up mass for $H^1$ solutions, equal to that of the standing wave solution, such that any solution with smaller mass has global existence and dispersive behavior. A comparable result was achieved in 2006 by Kenig, Merle \cite{kenig} for the $\dot H^1$-critical equation.

A similar statement is possible concerning the cubic nonlinearity studied here (which is $\dot H^{1/2}$-critical). The present paper does not address this question, but is a first step in that direction.

\subsection{The theory of Bates and Jones} In 1989, Bates, Jones \cite{bates} proved that the space of solutions decomposes into an unstable and a centre-stable manifold, for a large class of semilinear equations. As far as it concerns this paper, their result is the following: consider a Banach space $X$ and the semilinear equation
\be
u_t = A u + f(u),
\ee
under the assumptions
\begin{enumerate}
\item[H1] $A:X \to X$ is a closed, densely defined linear operator that generates a $C_0$ group.
\item[H2] The spectrum of $A$ decomposes into $\sigma(A) = \sigma_s(A) \cup \sigma_c(A) \cup \sigma_u(A)$ situated in the left half-plane, on the imaginary axis, and in the right half-plane respectively and $\sigma_s(A)$ and $\sigma_u(A)$ are bounded.
\item[H3] The nonlinearity $f$ is locally Lipschitz, $f(0) = 0$, and $\forall \epsilon>0$ there exists a neighborhood of $0$ on which $f$ has Lipschitz constant $\epsilon$.
\end{enumerate}
Furthermore, let $X^u$, $X^c$, and $X^s$ be the $A$-invariant subspaces corresponding to $\sigma_u$, $\sigma_c$, and respectively $\sigma_s$ and let $S^c(t)$ be the evolution generated by $A$ on $X^c$. Bates and Jones further assume that
\begin{enumerate}
\item[C1-2] $dim X^u$, $dim X^s < \infty$.
\item[C3] $\forall \rho>0$ $\exists M>0$ such that $\|S^c(t)\| \leq M e^{\rho|t|}$.
\end{enumerate}

Let $\Phi$ be the flow associated to the nonlinear equation. We call $\mc N \subset U$ $t$-invariant if, whenever $\Phi(s)v \in U$ for $s \in [0, t]$, $\Phi(s)v \in \mc N$ for $s \in [0, t]$.

Let $W^u$ be the set of $u$ for which $\Phi(t) u \in U$for all $t<0$ and decays exponentially as $t \to -\infty$. Also, consider the natural direct sum projection $\pi^{cs}$ on $X^c \oplus X^s$.
\begin{definition}
A \emph{centre-stable} manifold $\mc N \subset U$ is a Lipschitz manifold with the property that $\mc N$ is $t$-invariant relative to $U$, $\pi^{cs}(\mc N)$ contains a neighborhood of $0$ in $X^c \oplus X^s$, and $\mc N \cap W^u = \{0\}$.
\lb{centr}
\end{definition}

The result of \cite{bates} is then
\begin{theorem}
Under assumptions H1-H3 and C1-C3, there exists an open neighborhood $U$ of $0$ such that $W^u$ is a Lipschitz manifold which is tangent to $X^u$ at $0$ and there exists a centre-stable manifold $W^{cs} \subset U$ which is tangent to $X^{cs}$.
\lb{t1}
\end{theorem}

Gesztesy, Jones, Latushkin, Stanislavova \cite{ges} proved in 2000 that the abstract Theorem \ref{t1} applies to the semilinear Schr\"odinger equation. More precisely, their main result was that
\begin{theorem}
Given the equation
\be
iu_t - \Delta u - f(x, |u|^2)u - \beta u = 0
\ee
and assuming that
\begin{enumerate}
\item{H1} $f$ is $C^3$ and all derivatives are bounded on $U \times \set R^3$, where $U$ is a neighborhood of $0$;
\item{H2} $f(x, 0) \to 0$ exponentially as $x\to\infty$;
\item{H3} $\beta<0$;
\item{H4} $u_0$ is an exponentially decaying stationary solution to the equation (standing wave);
\end{enumerate}
then there exists a neighborhood of $u_0$ that decomposes into a centre-stable and an unstable manifold.
\end{theorem}

While providing an interesting answer to the problem, the main drawback of this approach is that one cannot infer the global in time behavior of the solutions on the centre-stable manifold. Indeed, once a solution leaves the specified neighborhood of $0$, one cannot say anything more about it, not even concerning its existence.

\subsection{The result of Schlag}
In \cite{schlag}, Schlag extended the method of modulation to the $L^2$-supercritical case and proved that in the neighborhood of each ground state of equation (\ref{NLS}) there exists a codimension-one Lipschitz submanifold of $H^{1}(\set R^3) \cap W^{1,1}(\set R^3)$ such that initial data on the submanifold lead to global solutions.

The method used in \cite{schlag} and applied in the current paper with some enhancements is the following: write the solution to equation (\ref{NLS}) as $\Psi = W + R$, where $W = e^{i\theta}\phi(x-y, \alpha)$ is a nonlinearly moving standing wave, determined by the parameter path $\pi = (\Gamma, D, \alpha, v)$ as in (\ref{def_W}), while $R$ is an error term that needs to be controlled. One obtains the nonlinear Schr\"odinger equation (\ref{Z_neliniar}) in $Z = \bpm R \\ \overline R \epm$, with the nonselfadjoint Hamiltonian
\be
\mc H_Z^{\pi} = \bpm \Delta+2|W(\pi)|^2 & W(\pi)^2 \\ -\overline W(\pi)^2 & -\Delta-2|W(\pi)|^2 \epm
\ee
and localized quadratic and nonlocalized cubic nonlinear terms on the right-hand side.

The spectrum of the Hamiltonian determines the properties of the equation. Following an appropriate transformation, it becomes real-valued and takes the form
\be
\mc H = \bpm \Delta + 2\phi(\cdot, \alpha)^2 - \alpha & \phi(\cdot, \alpha)^2 \\ -\phi(\cdot, \alpha)^2 & -\Delta - 2 \phi(\cdot, \alpha)^2 + \alpha\epm.
\lb{ham}
\ee
For the rest of this paper, we make the following standard \emph{spectral assumption}:
\begin{assumption}
$\mc H$ has no embedded eigenvalues in the interior of its essential spectrum for any $\alpha>0$.
\lb{assum}
\end{assumption}
Such assumptions are routinely made in the proof of asymptotic stability results, as for example in \cite{buslaev1}, \cite{cuc}, \cite{rod3}.

Even though Assumption \ref{assum} is expected to be true, it has not been proved to hold. Nevertheless, the assumption is most likely true generically, in the sense that embedded eigenvalues should, as a rule, vanish under perturbations by turning into resonances in the upper-half plane (by Fermi's rule), see \cite{cuc2}. Thus, even if Assumption~\ref{assum} fails in some particular case, one should be able to reinstate it by means of perturbations.

Under this assumption, we completely describe the spectrum of $\mc H$ following \cite{schlag}, with the proof delayed until the next section. It consists of an absolutely continuous part $(-\infty, -\alpha] \cup [\alpha, \infty)$ supported on the real axis, a generalized eigenspace at $0$ with $4$ eigenvectors and $4$ generalized eigenvectors. To each disconnected component of the spectrum there corresponds a Riesz projection (namely $P_c$, $P_{root}$, and $P_{im} = P_+ + P_-$ respectively) given by a Cauchy integral.

In the course of the proof, Schlag used the method of modulation. The necessity for it arises because the projection of the solution onto the generalized eigenspace of the Hamiltonian at zero does not disperse or satisfy Strichartz estimates. Physically, this corresponds to the fact that a nonzero displacement of the solution $\Psi$ relative to the soliton $W$ does not go away in time and that even a small  relative velocity can lead to a large displacement in finite time. Since the right-hand side terms of the equation keep introducing small perturbations, one constantly needs to adjust the soliton path in order to eliminate them from the generalized zero eigenspace.

One of the main contributions of Schlag \cite{schlag} was adapting the modulation method to the $L^2$-supercritical case. In this case, the main difficulty lies in dealing with the unstable mode of the equation, which corresponds to the imaginary eigenvalue $i\sigma$ of $\mc H$. To address this, \cite{schlag} showed that the solution of the linearized equation does not grow exponentially in time if and only if  the initial data $Z(0)$ is on a certain codimension-one manifold, tangent to $Ker (P_+(0))$. This choice eliminates the effect of the unstable eigenvalue.

In this manner, Schlag \cite{schlag} proved global existence and decay properties for the linearized equation with $H^{1} \cap W^{1, 1}$ initial data on a codimension-one manifold. A fixed point argument allowed him to go back to the nonlinear equation.

The main result of \cite{schlag} states the following:
\begin{theorem}
Impose the spectral Assumption \ref{assum} and fix $\alpha_0>0$. Then there exist a small $\delta > 0$ and a Lipschitz manifold $\mc N$ of size $\delta$ inside $W^{1, 2} \cap W^{1, 1}$, of codimension one, so that $\phi(\cdot, \alpha_0) \in \mc N$, with the following property: for any choice of initial data $\psi(0) \in \mc N$, the NLS equation (\ref{NLS}) has a
global $H^{1}$ solution $\psi(t)$ for $t \geq 0$. Moreover,
\be
\psi(t) = W(t, \cdot) + R(t)
\ee
where $W$ as in (\ref{def_W}) is governed by a path $\pi(t)$ of parameters so that $|\pi(t) - (0, 0, 0, \alpha_0)| \leq \delta$ and which converges to some terminal vector $\pi(\infty)$ such that $\sup_{t\geq0} |\pi(t) - \pi(\infty)| \leq C\delta$. Finally,
\be
\|R(t)\|_{H^1} \leq C\delta,\ \|R(t)\|_{\infty} \leq C\delta t^{-3/2}
\ee
for all $t > 0$, and there is scattering:
\be
R(t) = e^{it\Delta} f_0 + o_{L^2}(1) \text{ as } t \to \infty
\ee
for some $f_0 \in L^2(\set R^3)$.
\end{theorem}

The main problem here is that the $H^1 \cap W^{1, 1}$ space is not preserved by the flow. Starting with a function $\psi(0) \in \mc N$ of finite $H^{1} \cap W^{1, 1}$ norm at $t=0$ as initial data, there is no guarantee that $\psi(t)$ will still have finite $W^{1, 1}$ norm for any $t \ne 0$. Therefore, the question whether the manifold $\mc N$ is invariant under the Hamiltonian flow does not make sense in this context. One can replace the $W^{1, 2}(\set R^3) \cap W^{1,1}(\set R^3)$ norm with the stronger invariant $\Sigma^{5/2+\epsilon} = H^{5/2+\epsilon} \cap |x|^{-5/2-\epsilon} L^2$ norm, but this weakens the result considerably.

Another example of the same phenomenon, in the case of the wave equation, is given by Krieger, Schlag \cite{kri4}. For a more general survey of this topic, see \cite{schlag2}.

\subsection{Current paper} The result of this paper represents an improvement over that of Schlag \cite{schlag}, in that it holds in the $H^{1/2} \cap L^{4/3-\epsilon}$ norm, which is strictly weaker than the invariant $\Sigma = H^{1} \cap |x|^{-1} L^2$ space, a somewhat natural choice for equation (\ref{NLS}). In this space, the question concerning the manifold's invariance under the flow becomes meaningful and it turns out that the answer is affirmative.

This paper follows the method of proof of \cite{schlag} (namely the method of modulation, adapted to the $L^2$-supercritical case), but some important details differ.

The choice of $H^{1/2}$ for initial data is sharp and corresponds to the fact that the equation is $\dot H^{1/2}$-critical. It is possible only due to Keel-Tao endpoint Strichartz estimates for the linearized Hamiltonian. The endpoint corresponds exactly to using half a derivative to bound the nonlocalized cubic right-hand side term of the linearized equation.

The $L^{4/3-\epsilon}$ condition on the initial data leads to a $t^{-1}$ decay in $L^2$ in time of the solution that compensates for the possibility of linear growth in the modulation equations. This problem arises because of the generalized eigenspace of the Hamiltonian at $0$.

This $L^2$ in time decay bound is not sharp. We expect that, due to the oscillatory nature of the integrand, further improvements are achievable by using conditionally convergent integrals, instead of absolutely converging ones as in the current paper.

\subsection{Linear estimates} The first dispersive estimates concerning NLS with nonselfadjoint Hamiltonians are present in \cite{buslaev1}.

More recently, Erdogan, Schlag \cite{erdogan2} considered Hamiltonians of the form $\mc H = \mc H_0 + V$, where
\be \mc H_0 = \bpm -\Delta + \mu & 0\\0&\Delta-\mu\epm,\ V = \bpm -U & -W \\ W & U\epm.
\lb{H_general}
\ee
They made the following assumptions: that $-\sigma_3 V$ is a positive matrix, that $L_- = -\Delta + \mu + U + W \geq 0$, that $U$ and $W$ have polynomial decay, and the spectral Assumption \ref{assum}. Here $\sigma_3$ denotes the Pauli matrix $\bpm 1 & 0 \\ 0 & -1\epm$.

Under these conditions, Erdogan, Schlag \cite{erdogan2} proved the $L^2$ boundedness of the evolution $e^{it\mc H}$ for $|V| \leq C \langle x \rangle^{-1-\epsilon}$. In \cite{schlag}, Schlag proved the $L^1 \to L^{\infty}$ dispersive estimate and Strichartz nonendpoint estimates, for $\langle x \rangle^{-3-\epsilon}$ potential decay and under the further assumption that the edges of the spectrum are neither eigenvalues nor resonances. Erdogan, Schlag \cite{erdogan2} obtained corresponding results for nonselfadjoint Hamiltonians in the presence of a resonance or eigenvalues at the edges of the essential spectrum, if the potential decays like $\langle x \rangle^{-10-\epsilon}$. Yajima \cite{yajima} proved independenty the same result, assuming less decay on $V$.

This paper establishes the following Keel-Tao endpoint Strichartz estimates for a nonselfadjoint Hamiltonian of the form (\ref{H_general}):
\begin{corollary}
Suppose that $\mc H = \mc H_0 + V$, where
\be \mc H_0 = \bpm -\Delta + \mu & 0\\0&\Delta-\mu\epm,\ V = \bpm -U & -W \\ W & U\epm,
\ee
that $-\sigma_3 V$ is a positive matrix, that $L_- = -\Delta + \mu + U + W \geq 0$, that $|V|  \leq C \langle x \rangle^{-7/2-}$, that the spectral Assumption \ref{assum} holds, and that the edges of the spectrum $\pm \mu$ are neither eigenvalues nor resonances.

Then the evolution $e^{it\mc H}P_c$ satisfies the following Strichartz-type estimates:
\be\ba{rl}
\|e^{it\mc H} P_c f\|_{L^q_t L^r_x} &\leq C \|f\|_2,\\
\ds\Big\|\int e^{-is\mc H} P_c F(s) \dd s\Big\|_2 &\leq C \|F\|_{L^{q'}_t L^{r'}_x},\\
\ds\Big\|\int_{s<t} e^{it\mc H}e^{-is\mc H^*} P_c F \dd s\Big\|_{L^q_t L^r_x} &\leq C \|F\|_{L^{\tilde q'}_t L^{\tilde r'}_x},\\
\ds\Big\|\int_{s<t} e^{i(t-s) \mc H} P_c F \dd s\Big\|_{L^q_t L^r_x} &\leq C \|F\|_{L^{\tilde q'}_t L^{\tilde r'}_x}.
\ea\ee
for any sharply admissible $(q, r)$ (that is, such that $2\leq q$, $r \leq \infty$, $\ds\frac 1 q + \frac 3 {2r} = \frac 3 4$) and $(\tilde q, \tilde r)$. The same estimates hold after swapping $\mc H$ and $\mc H^*$.
\lb{Strichartz}
\end{corollary}

Note that Corollary \ref{Strichartz} is not an immediate consequence of \cite{tao}, because $e^{it\mc H} e^{-is \mc H^*} \neq e^{i(t-s) \mc H}$.

The exact rate of decay of $V$ does not matter for the purpose of this paper, since we deal only with exponentially decaying potentials. However, it is important to have the endpoint Strichartz estimate for two reasons. Firstly, by linearizing the equation one obtains small localized linear terms on the right-hand side and it is useful to be able to bound their contribution using the $L^2_t \to L^2_t$ endpoint Strichartz estimate. Secondly, as mentioned before, the sharp estimate allows one to use exactly half a derivative in handling the nonlocalized cubic terms on the right-hand side.

The difficulty in the proof lies in the fact that $\mc H$ is not selfadjoint, so the usual $L^1 \to L^{\infty}$ dispersive estimate does not imply the endpoint estimate Corollary \ref{Strichartz}. Therefore, we use the following strengthened version of it:
\begin{proposition}
Under the assumptions of Corollary \ref{Strichartz},
\be
\|e^{it\mc H} P_c e^{-is\mc H^*} P_c^*\|_{1 \to \infty} \leq C|t-s|^{-3/2}.
\lb{II2}
\ee
\lb{dispers}
\end{proposition}
The proof of this statement is a generalization of the one given in \cite{schlag} for  the usual dispersive estimate. The argument uses the spectral representation of the evolution from that paper and the finite Born sum expansion of the resolvent for both the operator and its adjoint.

Once established, the estimate (\ref{II2}), together with the $L^2$ theory of \cite{erdogan2}, makes possible to apply the methods of \cite{tao}, leading to Corollary \ref{Strichartz}. 

Now we return to the nonlinear problem. Without loss of generality, take any standing wave $W(0)$ and transform it, by means of a symmetry transformation $\mc G_0$, into a positive ground state $\phi(\cdot, \alpha_0)$ of equation (\ref{NLS}). Then let $P_+(0)$ and $P_{root}(0)$ be the Riesz projections onto the eigenspace corresponding to the eigenvalue $i\sigma$ of positive imaginary part and respectively onto the generalized zero eigenspace of the linearized Hamiltonian $\mc H$ (\ref{ham}) at time $0$. Furthermore, let $f^+(0)$, $\tilde f^+(0)$, and $\eta_F(0)$, $\xi_F(0)$ be the normalized eigenvectors of $\mc H$ and $\mc H^*$ at $i\sigma$ and the generalized eigenvectors of $\mc H$ and $\mc H^*$ at $0$, respectively. All are exponentially decaying Schwartz functions and
\be
P_+(0) = \langle \cdot, \tilde f^+(0) \rangle f^+(0),\ P_{root} = \sum_F \langle \cdot, \xi_F(0) \rangle \eta_F (0).
\ee

In the sequel we use the notation $L^{p \cap q} = L^p \cap L^q$ and $L^{p+q} = L^p + L^q$.

Following these preparations, we state a more technical result from which the main theorem follows almost immediately. For simplicity, we first state it in the case when the initial data is in the neighborhood of a positive ground state and its projection on the generalized zero eigenspace vanishes.

\begin{theorem} Assume that $W(0) = \phi(\cdot, \alpha_0)$ is a positive ground state of equation (\ref{phi}). For $1\leq q <4/3$, let $S_{\delta}$ be given by
\be
S_{\delta} = \Big\{R_0 \in H^{1/2}(\set R^3) \cap L^{q}(\set R^3) \mid \|R_0\|_{H^{1/2} \cap L^{q}} < \delta,\ (P_+(0) + P_{root}(0)) \bpm R_0 \\ \overline R_0\epm = 0\Big\}.
\ee
Then, for some small $\delta$, there exists a map $F:S_{\delta} \rightarrow H^{1/2}(\set R^3) \cap L^{q}(\set R^3)$, whose range is spanned by a Schwartz function, given by
\be
\mc F(R_0) = h(R_0) {f^+(0)} = \bpm \mc F_1 \\ \mc F_2 \epm
\ee
such that
\begin{enumerate}
\item
$\|\mc F(R_0)\| \leq C\|R_0\|^2_{H^{1/2} \cap L^{q}}$
\item
$\|\mc F(R_0^1) - \mc F(R_0^2)\| \leq C \delta \|R_0^1-R_0^2\|_{H^{1/2} \cap L^{q}}$
\end{enumerate}
and, for every $R_0 \in S(\delta)$, the equation (\ref{NLS}) having $\Psi(R_0)(0) = W(0) + R_0 + \mc F_1(R_0)$ as initial data admits a global solution $\Psi(R_0)$. Moreover, the solution $\Psi(R_0)$ has the following properties:
\begin{enumerate}
\item
$\Psi(R_0)$ depends Lipschitz continuously on $R_0$,
\be
\|\Psi(R_0^1) - \Psi(R_0^2)\|_{\langle t \rangle^{1/2-\epsilon} L^{2}_t \langle x \rangle L^{6+\infty}_x} \leq C \|R_0^1 - R_0^2\|_{2}.
\ee
\item
There exists a parameter path $\pi$ with $\pi(0) = (0, 0, 0, \alpha_0)$ and $\|t\dot\pi(t)\|_1 + \|\dot \pi\|_{1\cap \infty} < \infty$ such that $\Psi(R_0)$ stays close to $W(\pi)$ for all time $t \geq 0$: $\Psi(R_0) = W(\pi) + R$, where
\be
\|R\|_{L^{\infty}_t H^{1/2}_x \cap L^2_t W^{1/2, 6}_x \cap \langle t \rangle^{-1/2-\epsilon} L^{2}_t L^{6+\infty}_x} \leq \delta
\ee
and one has scattering: for some $f_0$ in $L^2$,
\be
R(t) = e^{it\Delta} f_0 + o_{L^2}(1).
\ee
\end{enumerate}
\lb{main}
\end{theorem}

The map $R_0 \mapsto \Psi(R_0)(0) = W(0) + R_0 + \mc F_1(R_0)$ takes $S_{\delta}$ to a codimension-nine submanifold $\mc N_9$ of $H^{1/2} \cap L^q$. Indeed, the map is Lipschitz bicontinuous for sufficiently small $\delta$ and $S_{\delta}$ is an open set in a codimension-nine linear space.

Since we want to extend this result to more general standing waves instead of just ground states, we conjugate everything by means of symmetry transformations. Also note that the codimension-nine manifold $\mc N_9$ provided by Theorem \ref{main} becomes, after applying symmetry transformations, a codimension-one submanifold. These observations are summarized in the following corollary:

\begin{corollary}
Consider any standing wave $W(0)$. Under the same assumptions as in Theorem \ref{main}, there exists a codimension-one Lipschitz manifold $\mc N_{L^q \cap H^{1/2}}$ in $L^q \cap H^{1/2}(\set R^3)$, $1\leq q < 4/3$, given locally by
\be
\mc N_{L^q \cap H^{1/2}} = \bigcup_{\frak g} \frak g (\mc N_9),
\ee
whose tangent space at $W(0)$ is $Ker(P_+(0))$, such that for initial data $\Psi(0)$ on the manifold $\mc N_{L^q \cap H^{1/2}}$ the equation has a global dispersive solution $\Psi$, with the same properties as in Theorem \ref{main}, but with respect to some more general standing wave $\frak g(W(0))$, such that $|\frak g| \leq C \|R_0\|_{L^q \cap H^{1/2}}$, instead of simply $W(0)$.
\label{cor17}
\end{corollary}

A straightforward consequence is that the same result holds in the strictly stronger norm of $\Sigma^1 = H^{1} \cap |x|^{-1} L^2$, which has the advantage of being locally invariant under the flow. Furthermore, in this topology one can identify $\mc N$ as the centre-stable manifold of \cite{bates}, from the previous discussion. This leads to the main Theorem \ref{tichie}, stated on the first page.

We remark that $\Sigma^1$ can be replaced in this statement with any invariant $\Sigma^s = H^s \cap |x|^{-s} L^2$ space, for $s>3/4$, this being the minimal requirement so that $\Sigma^s \subset L^q \cap H^{1/2}$ for some $q<4/3$.

Acknowledgment: I would like to thank Professor Wilhelm Schlag for his suggestions and for his very careful reading of this paper.

\section{Proof of the Nonlinear Results}
\subsection{Formulation of the problem}
We aim to prove that there exists a codimension-one submanifold of $H^{1/2} \cap L^{4/3-\epsilon}$ on which the focussing cubic nonlinear Schr\"odinger equation (\ref{NLS}) has global solutions. Throughout this section we employ the Keel-Tao endpoint Strichartz estimates of Section \ref{lin_est}.

In the $L^2$-subcritical case, Cazenave, Lions \cite{cazenave} and Weinstein \cite{wein2} proved that stability occurs for any solution that starts in a sufficiently small neighborhood of the standing wave manifold. However, the presence of an unstable eigenstate of the linearization precludes one from achieving such a result in the $L^2$-supercritical case and Berestycki, Lions \cite{bere} prove that arbitrarily small perturbations of the ground state may lead to blowup in finite time. The best that one can hope for is the existence of a codimension-one manifold on which the evolution does not lead to blowup. This is indeed the result proved by Schlag \cite{schlag} and improved here.

Let $\phi = \phi(\cdot, \alpha)$ be the radially symmetric ground state (meaning $\phi>0$) of the semilinear Schr\"odinger operator corresponding to energy $\alpha>0$, that is a solution of (\ref{phi}). The existence of such solutions to equation (\ref{phi}) was proved by Berestycki and Lions in \cite{bere}, who further showed that they are infinitely differentiable and exponentially decaying. Uniqueness was established by Coffman \cite{coffman} for (\ref{phi}) and Kwong \cite{kwong} and McLeod, Serrin \cite{mcl} for more general nonlinearities.

In the particular case of the cubic nonlinearity, the equation (\ref{phi}) and its solutions have the scaling invariance $\phi(x, \alpha) = \alpha^{1/2} \phi(\alpha^{1/2} x, 1)$.

Note that $e^{it\alpha} \phi(x, \alpha)$ is a $1$-parameter family of periodic solutions for equation (\ref{NLS}). Starting from it, one can obtain more solutions by taking advantage of the symmetries of equation (\ref{NLS}). Applying the following family of transformations
\be
\mc G(t)(f(x, t)) = e^{i(\Gamma+vx-t|v|^2)} \alpha^{1/2} f(\alpha^{1/2}x-2tv-D, \alpha t)
\lb{coord}
\ee
to $e^{it} \phi(\cdot, 1)$, the result is a wider 8-parameter family of solutions to (\ref{NLS})
\be
\mc G(t)(e^{it} \phi(x, 1)) = e^{i(\Gamma + vx - t|v|^2 + \alpha t)} \alpha^{1/2} \phi(\alpha^{1/2}x-2tv-D, 1)
\ee
or, after reparametrizing,
\be
\mc G(t)(e^{it} \phi(x, 1)) = e^{i(\Gamma + vx - t|v|^2 + \alpha t)} \phi(x-2tv-D, \alpha),
\ee
which we call standing waves.

Here $\mc G$ is composed of a Galilean coordinate change, with six degrees of freedom corresponding to $v$ and $D$, a phase change represented by $\Gamma$, and a rescaling embodied by $\alpha$. Henceforth we call such $\mc G$ as in (\ref{coord}) symmetry transformations, since they correspond to the symmetries of equation (\ref{NLS}).

In the sequel we consider the pairs made of a function and its conjugate instead of just the function alone. For example, by a standing wave we will also mean the pair $\bpm \mc G(t)(e^{it} \phi(x, 1)) \\ \ov{\mc G(t)(e^{it} \phi(x, 1))} \epm$. There is an obvious correspondence between the pair and its first component, as long as the components are conjugate to one another. All the column two-vectors that apear in this paper will have this property, related to the fact that the vector form of equation (\ref{NLS}) is $\ds\bpm 0 & 1 \\ 1 & 0\epm$-invariant. 

The question arises whether standing waves are stable under small perturbations. We seek perturbed solutions of the form $\psi = W(\pi)+R$ with small $R$, where
\be\begin{aligned}
W(\pi(t))(x) &= e^{i \theta(t, x)} \phi(x-y(t), \alpha(t)) \\
\theta(t, x) &\ds= v(t) x - \int_0^t(|v(s)|^2 - \alpha(s)) \dd s + \gamma(t)\\
y(t) &\ds= 2\int_0^t v(s) \dd s + D(t).
\end{aligned}
\lb{def_W}
\ee
$W(\pi)$ represents a moving soliton governed by the parameter path $\pi = (\Gamma, \alpha, D_i, v_i)$. We look for solutions $\psi$ that remain close to the 8-dimensional manifold of solitons for all positive times $t>0$, hence to a moving soliton like $W(\pi)$.

\subsection{Setting up the contraction scheme}
Assume that all the parameters describing $W(\pi)$, namely $\gamma$, $\alpha$, $D_i$, and $v_i$, have limits as $t \to \infty$, denoted $\gamma(\infty)$ etc.. It is more convenient in the sequel to consider an alternative to $\gamma$, namely a new parameter $\Gamma$ such that
\be
\dot\Gamma = \dot\gamma + \dot v (2tv(\infty) + D(\infty)),
\ee
more precisely $\Gamma(t) = \gamma(t) + (v(t)-v(\infty)) D(\infty) - \int_t^{\infty} 2s \dot v(s) v(\infty) \dd s$.

Henceforth, we assume that
\be
\|\dot\pi\|_{\infty} + \|\langle t \rangle \dot \pi(t)\|_1 < \infty,
\ee
where $\pi(t) = (v_k(t), D_k(t), \alpha(t), \Gamma(t))$. Note that $\gamma$ can be recovered from $\pi$ and that $\|\dot\gamma\|_1<\infty$ under our assumption; also $\gamma(\infty) = \Gamma(\infty)$.

For $F \in \{v_k, D_k, \alpha, \Gamma\}$, denote by $\xi_F^Z$ the following family of vectors:
\be\begin{split}
\xi_{D_k}^Z(t) &= \bpm e^{i\theta(x, t)} x_k \phi(x-y(t), \alpha(t)) \\ e^{-i\theta(x, t)} x_k \phi(x-y(t), \alpha(t))\epm\\
\xi_{v_k}^Z(t) &= \bpm ie^{i\theta(x, t)} \partial_k \phi(x-y(t), \alpha(t)) \\ -ie^{-i\theta(x, t)} \partial_k \phi(x-y(t), \alpha(t))\epm\\
\xi_{\Gamma}^Z(t) &= \bpm ie^{i\theta(x, t)}\partial_{\alpha} \phi(x-y(t), \alpha(t)) \\ -ie^{-i\theta(x, t)}\partial_{\alpha} \phi(x-y(t), \alpha(t))\epm\\
\xi_{\alpha}^Z(t) &= \bpm e^{i\theta(x, t)} \phi(x-y(t), \alpha(t)) \\ e^{-i\theta(x, t)} \phi(x-y(t), \alpha(t))\epm.
\end{split}\ee
Also let $\eta^Z_F = -i\sigma_3 \xi^Z_F = \Big(\ba{cc} -i & 0 \\ 0 & i \ea\Big) \xi^Z_F$.

Their immediate importance is that $\eta^Z_F(t)$ span the tangent space of the $8$-dimensional standing wave manifold at $\bpm W(\pi(t))\\\ov{W(\pi(t))}\epm$, for each individual $t \geq 0$, and $\xi_F^Z(t)$ form a dual basis with respect to the usual dot product.

From another perspective, note that if $v=D=\gamma=0$ and $W$ is a positive ground state of the equation, then $\eta_F^Z$ span the generalized eigenspace of the linearized Hamiltonian $\mc H$ (\ref{ham}) at zero and $\xi_F^Z$ fulfill the same function for its adjoint $\mc H^*$. However, the property of being an eigenvector is not preserved under symmetry transformations, so this characterization is no longer true when $W$ is a  more general standing wave instead of a positive ground state.

The following lemma exhibits the equation satisfied by the error term $R$:
\begin{lemma}
$\Psi = W(\pi) + R$ is a solution of equation (\ref{NLS}) if and only if $Z = \bpm z_1=R\\z_2=\overline R\epm$ is a solution to
\be
i \partial_t Z + \mc H_{\pi}^Z Z = -i\dot\pi \partial_{\pi} \Big(\ba{c}W(\pi) \\ \ov{W(\pi)}\ea\Big) + N^Z(Z, \pi),
\lb{Z_neliniar}
\ee
where
\be
\mc H_{\pi}^Z = \Big(\ba{cc} \Delta+2|W(\pi)|^2 & W(\pi)^2 \\ -\overline W(\pi)^2 & -\Delta-2|W(\pi)|^2 \ea\Big),
\ee
\be\begin{split}
\dot\omega \partial_{\pi} \bpm W(\pi) \\ \ov{W(\pi)}\epm = & \dot v^{\omega} \bpm ie^{i\theta}x \phi(x-y, \alpha) \\ -ie^{-i\theta}x \phi(x-y, \alpha)\epm + \dot \gamma^{\omega} \bpm ie^{i\theta}\phi(x-y, \alpha) \\ -ie^{-i\theta}\phi(x-y, \alpha)\epm + \\
&+ \dot\alpha^{\omega} \bpm e^{i\theta}\partial_{\alpha} \phi(x-y, \alpha) \\ e^{-i\theta}\partial_{\alpha} \phi(x-y, \alpha)\epm + \dot D^{\omega} \bpm -e^{i\theta}\dl\phi(x-y, \alpha) \\ -e^{-i\theta}\dl\phi(x-y, \alpha)\epm\\
= & \ds\dot\alpha^{\omega} \eta^Z_{\Gamma} - \dot\gamma^{\omega} \eta^Z_{\alpha} - \sum_k (\dot D_k^{\omega} \eta^Z_{v_k} + \dot v_k^{\omega} \eta^Z_{D_k})
\end{split}\lb{Wpi}\ee
(with $\theta = \theta(x, t)$, $\alpha=\alpha(t)$, $y = y(t)$) and
\be
N^Z(Z, \pi) = \bpm -2|z_1|^2 W(\pi) - \overline W(\pi) z_1^2 - |z_1|^2 z_1 \\ 2|z_2|^2 \overline W(\pi) + W(\pi) z_2^2 + |z_2|^2 z_2 \epm.
\ee
\end{lemma}

We wrote $\dot \pi \partial_{\pi} W(\pi)$ in a more general form that becomes convenient when linearizing the equation. In the linearized equation, the family of vectors $\eta_F$ depends on one path and the coefficients $\dot v$ etc. depend on another.

\begin{proof} By direct computation. $\Psi = W(\pi)+R$ satisfies equation (\ref{NLS}),
\begin{multline}
i\partial_t (W+R) + \Delta(W+R) = -|W+R|^2 (W+R)\\
= -|W|^2 W - 2|W|^2 R -W^2 \ov R - 2W |R|^2  - 2 \ov W R^2 - |R|^2 R,
\end{multline}
while $W(\pi)$ fulfills the identity
\be
i\partial_t W(\pi) + \Delta W(\pi) = -|W(\pi)|^2 W(\pi) + i \dot \pi \partial_{\pi} W(\pi).
\ee
Subtracting the two relations, one has that
\be
i\partial_t R + \Delta R = - 2|W|^2 R -W^2 \ov R - 2W |R|^2  - 2 \ov W R^2 - |R|^2 R - i \dot \pi \partial_{\pi} W(\pi).
\ee
Joining this equation with its conjugate, one obtains exactly (\ref{Z_neliniar}).
\end{proof}

Consider the following linearized version of equation (\ref{Z_neliniar}), in which we partly replace $\pi$ and $Z$ with the auxiliary functions $\pi^0$ and $Z^0$:
\be
i Z_t + \mc H^Z_{\pi^0} Z = -i\dot\pi \partial_{\pi} \bpm W(\pi^0) \\ \ov{W(\pi^0)}\epm + N^Z(Z^0, \pi^0).
\lb{Z}
\ee

We choose $W(\pi)$ such that at each time $t$ it satisfies the orthogonality condition
\be
\langle Z(t), \xi^Z_F(t)\rangle = 0,
\lb{ort_Z}
\ee
which leads to a system of modulation equations for the path $\pi$. If the standing wave $W(\pi)$ is a positive ground state of the equation, this simply means that $P_{root} Z(t) = 0$, that is the projection of $Z(t)$ onto the generalized eigenspace at $0$ of the Hamiltonian $\mc H$ (\ref{ham}) is $0$. Otherwise, the condition takes a more complicated meaning.

If we try to approximate $\mc H^Z_{\pi^0}$ by a constant Hamiltonian in order to solve equation (\ref{Z_neliniar}), the problem is that the potential moves with nonzero velocity along the path described by $\pi$. Thus, we need to change the reference frame to one that moves with the same speed as $\mc H_{\pi^0}^Z$. However, since $\mc H^Z_{\pi^0}$ does not move with constant speed, we have (in order to avoid gradient terms in the equation) to choose a reference frame moving at the constant speed that best approximates the speed of $\mc H^Z_{\pi^0}$. The same considerations apply to the phase of $\mc H^Z_{\pi^0}$. Therefore, we need to determine the uniform movement path that best approximates $\pi^0$.

Define the following limit values associated to any path $\pi$:
\be\begin{split}
\Gamma_{\infty} &\ds= \gamma(\infty) - \int_0^{\infty} (v^2(\infty) - v^2(s) - \alpha(\infty) + \alpha(s)) \dd s \\
v_{\infty} &= v(\infty)\\
D_{\infty} &\ds= D(\infty) - 2\int_0^{\infty} (v(\infty)-v(s)) \dd s\\
\alpha_{\infty} &= \alpha(\infty)\\
\phi_{\infty} &= \phi(\cdot, \alpha_{\infty}).
\end{split}
\ee

Given some parameter path $\pi$ such that $\|\dot \pi\|_{1\cap\infty} + \|t \dot \pi\|_1 < \infty$, one can distinguish a symmetry transformation (as in (\ref{coord})) $\mc G_{\pi}$ associated to $\pi$,
\be
\mc G_{\pi}(t)(f(x)) = e^{-i\theta_{\infty}(t, x+y_{\infty}(t))} f(x+y_{\infty}(t)),
\ee
where
\be\begin{split}
y_{\infty}(t) &= 2tv_{\infty}+D_{\infty}\\
\theta_{\infty}(t, x) &= \Gamma_{\infty} + v_{\infty} x - t(|v_{\infty}|^2 - \alpha_{\infty}).
\end{split}
\ee
Also consider the corresponding transformation for column two-vectors,
\be
\frak g_{\pi}(t) \bpm u_1 \\ u_2 \epm = \bpm \mc G_{\pi}(t) u_1 \\ \overline {\mc G_{\pi}(t) \overline u_2} \epm.
\ee
Note that $\mc G_{\pi}$ and $\frak g_{\pi}$ only depend on the terminal values $D_{\infty}$, $v_{\infty}$, $\Gamma_{\infty}$, and $\alpha_{\infty}$.

For future reference, let
\be
\rho_{\infty}(t, x) = \theta(t, x+y_{\infty}) - \theta_{\infty}(t, x+y_{\infty}).
\ee

$\mc G_{\pi} W(\pi)$ is close to a constant ground state $\phi(\cdot, \alpha)$ and it turns out that the best uniformly moving approximation to $\pi$ is provided by the constant path $(\Gamma_{\infty}, \alpha_{\infty}, D_{\infty}, v_{\infty})$.

Therefore, we apply the transformation $\frak g_{\pi^0}$ to the linearized equation (\ref{Z}). In this context it is natural to introduce the families of functions
\be\ba{c}
\eta_{\alpha}(t) = \bpm -i e^{i\rho_{\infty}} \phi(x+y_{\infty}-y, \alpha) \\ ie^{-i\rho_{\infty}} \phi(x+y_{\infty}-y, \alpha)\epm = \frak g_{\pi^0}(t) \eta^Z_{\alpha}(t)\\
\eta_{\Gamma}(t) = \bpm e^{i\rho_{\infty}} \partial_{\alpha} \phi(x+y_{\infty}-y, \alpha) \\ e^{-i\rho_{\infty}} \partial_{\alpha} \phi(x+y_{\infty}-y, \alpha)\epm = \frak g_{\pi^0}(t) \eta^Z_{\Gamma}(t)\\
\eta_{v_k}(t) = \bpm e^{i\rho_{\infty}} \partial_k \phi(x+y_{\infty}-y, \alpha) \\ e^{-i\rho_{\infty}} \partial_k \phi(x+y_{\infty}-y, \alpha)\epm = \frak g_{\pi^0}(t) \eta^Z_{v_k}(t)\\
\eta_{D_k}(t) = \bpm -ix_k e^{i\rho_{\infty}} \phi(x+y_{\infty}-y, \alpha) \\ ix_k e^{-i\rho_{\infty}} \phi(x+y_{\infty}-y, \alpha)\epm = \frak g_{\pi^0}(t) \eta^Z_{D_k}(t) - (y_{\infty})_k \frak g_{\pi^0}(t) \eta^Z_{\alpha}(t)
\ea\ee
(where $\rho_{\infty} = \rho_{\infty}(x, t)$, $y_{\infty} = y_{\infty}(t)$, $y=y(t)$, $\alpha = \alpha(t)$) and
\be
\xi_F(t) = i\sigma_3 \eta_F(t) = \left\{\ba{ll}\frak g_{\pi^0}(t) \xi^Z_F(t), &\text{for } F \in \{\alpha, \Gamma, v_k\} \\ \frak g_{\pi^0}(t) \xi^Z_{D_k}(t) - (y_{\infty})_k \frak g_{\pi^0}(t) \xi^Z_{\alpha}(t), &\text{for } F=D_k. \ea\right.
\ee
We made a change in the definitions of $\eta_{D_k}$ and $\xi_{D_k}$ so that these functions would be uniformly bounded in time, instead of linearly increasing as they would have been if we had just applied the symmetry transformation $\frak g_{\pi^0}$. The two families $\eta_F$ and $\xi_F$ span the generalized $0$ eigenspaces of $\mc H$ and $\mc H^*$ respectively, if $W(\pi(^0t))$ is a standing wave. Furthermore, $\eta_F$ span the tangent space of the eight-dimensional standing wave manifold at $\frak g_{\pi^0}(t) \bpm W(\pi^0(t))\\\ov{W(\pi^0(t))}\epm$ at each individual $t\geq 0$.

\begin{lemma}
$(Z, \pi)$ is a solution of equation (\ref{Z}) if and only if $U = \frak g_{\pi^0} Z$ and $\pi$ satisfy
\be
i U_t + \mc H_{\pi^0} U = -i\dot \pi \partial_{\pi}\frak g_{\pi^0} \bpm W(\pi^0) \\ \ov{W(\pi^0)}\epm + N(U^0, \pi^0) + V_{\pi^0}U,
\lb{nlsu}
\ee
where $U^0 = \frak g_{\pi^0} Z^0$ and we used the notations
\be
\mc H_{\pi^0}(x) = \bpm \Delta + 2\phi_{\infty}^2(x) - \alpha_{\infty} & \phi_{\infty}^2(x) \\ -\phi_{\infty}^2(x) & -\Delta - 2 \phi_{\infty}^2(x) + \alpha_{\infty} \epm,
\ee
\be
V_{\pi^0} = \bpm 2(\phi_{\infty}^2(x) - \phi^2(x+y_{\infty}-y, \alpha)) & \phi_{\infty}^2(x) - e^{2i\rho_{\infty}} \phi^2(x+y_{\infty}-y, \alpha) \\ -\phi_{\infty}^2(x) + e^{-2i\rho_{\infty}}\phi^2(x+y_{\infty}-y, \alpha) & -2(\phi_{\infty}^2(x) - \phi^2(x+y_{\infty}-y, \alpha)) \epm,
\ee
\be\ba{rl}
\dot \omega \partial_{\pi}\frak g_{\pi} \bpm W(\pi) \\ \ov{W(\pi)} \epm = 
&\ds\dot\alpha^{\omega} \eta_{\Gamma} - \dot\gamma^{\omega} \eta_{\alpha} - \sum_k (\dot D_k^{\omega} \eta_{v_k} + \dot v_k^{\omega} \eta_{D_k}),\lb{pi}
\ea\ee
and
\be
N(U, \pi) = \bpm -2|u_1|^2 e^{i\rho_{\infty}} \phi(x+y_{\infty}-y) - u_1^2 e^{-i\rho_{\infty}} \phi(x+y_{\infty}-y) - |u_1|^2 u_1 \\ 2|u_2|^2 e^{-i\rho_{\infty}} \phi(x+y_{\infty}-y) + u_2^2 e^{-i\rho_{\infty}} \phi(x+y_{\infty}-y) + |u_2|^2 u_2 \epm.
\lb{defNUpi}
\ee
\lb{transformare}
\end{lemma}
Here $U = \bpm u_1 \\ u_2 \epm$ is a $\set C^2$-valued function and $\pi$, $\pi^0$, $\omega$ are paths. We wrote the term $\ds\dot \pi \partial_{\pi}\frak g_{\pi^0} \bpm W(\pi^0) \\ \ov{W(\pi^0)}\epm$ in a more general form, in order to exhibit its dependence on two paths, $\pi$ and $\pi^0$. We also recall the notations $\rho_{\infty} = \rho_{\infty}(x, t)$, $y_{\infty} = y_{\infty}(t)$, $y=y(t)$, $\alpha = \alpha(t)$, $\phi_{\infty} = \phi(\cdot, \alpha_{\infty}$.

Henceforth, $\rho_{\infty}$, $y$, $y_{\infty}$, $\phi$, and $\phi_{\infty}$ will refer to quantities derived from $\pi^0$.

Let $\mc H = \mc H_{\pi^0} + V_{\pi^0}$. $\mc H$ is the Hamiltonian of equation (\ref{nlsu}), but we split it into a constant part $\mc H_{\pi^0}$ and an error term $V_{\pi^0}$.

\begin{proof} Firstly, we compute the following commutators:
\be
[\partial_t, \frak g_{\pi}] = \big(-(|v_{\infty}|^2 + \alpha_{\infty}) i\sigma_3 + 2 v_{\infty} \dl\big) \frak g_{\pi}
\ee
and
\be
[\Delta \sigma_3, \frak g_{\pi}] = \big(-|v_{\infty}|^2 \sigma_3 - 2 v_{\infty} i\dl \big) \frak g_{\pi}.
\lb{comut2}
\ee
Plugging $Z = \frak g_{\pi^0}^{-1} U$ and $Z^0 = \frak g_{\pi^0}^{-1} U^0$ into (\ref{Z}), we then have
\begin{multline}
\frak g_{\pi^0}^{-1} (i U_t + \Delta \sigma_3 U) + [i\partial_t + \Delta \sigma_3, \frak g_{\pi^0}^{-1}] U = \\
= \bpm 2|W|^2 & W^2 \\ -\ov{W^2} & -2|W|^2 \epm \frak g_{\pi^0}^{-1} U -i \dot\pi \partial_{\pi} \bpm W(\pi^0) \\ \ov{W(\pi^0)} \epm + N^Z(\frak g_{\pi^0}^{-1} U^0, \pi^0).
\end{multline}
Therefore, applying $\frak g_{\pi^0}$ to this equation and taking into account the fact that
\begin{multline}
\frak g_{\pi^0} \bpm 2|W|^2 & W^2 \\ -\ov{W^2} & -2|W|^2 \epm \frak g_{\pi^0}^{-1} =\\
= \bpm 2\phi^2(x+y_{\infty}-y, \alpha) & e^{2i\rho_{\infty}} \phi^2(x+y_{\infty}-y, \alpha) \\ -e^{-2i\rho_{\infty}} \phi^2(x+y_{\infty}-y, \alpha) & -2\phi^2(x+y_{\infty}-y, \alpha) \epm
\end{multline}
and
\be
\frak g_{\pi^0} ([i\partial_t + \Delta \sigma_3, \frak g_{\pi^0}^{-1}]) = -[i\partial_t + \Delta\sigma_3, \frak g_{\pi^0}] \frak g_{\pi^0}^{-1} = -\alpha_{\infty} \sigma_3
\ee
we obtain
\begin{multline}
i U_t + \Delta \sigma_3 U - \alpha_{\infty} \sigma_3 U = \\
\begin{aligned}
=&\bpm2\phi^2(x+y_{\infty}-y, \alpha) & e^{2i\rho_{\infty}} \phi^2(x+y_{\infty}-y, \alpha) \\ -e^{-2i\rho_{\infty}} \phi^2(x+y_{\infty}-y, \alpha) & -2\phi^2(x+y_{\infty}-y, \alpha) \epm U - \\
&- i\dot\pi \partial_{\pi} \frak g_{\pi^0} W(\pi^0) + \frak g_{\pi^0} N^Z(\frak g_{-\pi^0} U^0, \pi^0).
\end{aligned}\end{multline}
However, by definition
\be
\Delta \sigma_3 U - \alpha_{\infty} \sigma_3 U - \bpm2\phi^2(x+y_{\infty}-y, \alpha) & e^{2i\rho_{\infty}} \phi^2(x+y_{\infty}-y, \alpha) \\ -e^{-2i\rho_{\infty}} \phi^2(x+y_{\infty}-y, \alpha) & -2\phi^2(x+y_{\infty}-y, \alpha) \epm U = \mc H U.
\ee
Note that
\be\ba{rl}
\dot\pi \partial_{\pi} \frak g_{\pi^0} \bpm W(\pi^0) \\ \ov{W(\pi^0)} \epm = & \ds\dot\alpha \frak g_{\pi^0} \eta^Z_{\Gamma} - \dot\gamma \frak g_{\pi^0} \eta^Z_{\alpha} - \sum_k (\dot D_k \frak g_{\pi^0} \eta^Z_{v_k} + \dot v_k \frak g_{\pi^0} \eta^Z_{D_k}).
\ea\ee
Here the important fact is that
\begin{multline}
\sum_k \dot v_k \frak g_{\pi^0} \eta^Z_{D_k} + \dot\gamma \frak g_{\pi^0} \eta^Z_{\alpha} = \\
\ba{l}\ds= \sum_k \dot v_k \bpm -i(x+y_{\infty}) e^{i\rho_{\infty}} \phi(x+y_{\infty}-y, \alpha) \\ i(x + y_{\infty}) e^{-i\rho_{\infty}} \phi(x+y_{\infty}-y, \alpha)\epm + \dot \gamma \bpm -i e^{i\rho_{\infty}} \phi(x+y_{\infty}-y, \alpha) \\  ie^{-i\rho_{\infty}} \phi(x+y_{\infty}-y, \alpha)\epm \\
\ds= \sum_k \dot v_k \bpm -ix e^{i\rho_{\infty}} \phi(x+y_{\infty}-y, \alpha) \\ ix e^{-i\rho_{\infty}} \phi(x+y_{\infty}-y, \alpha)\epm + \dot \Gamma \bpm -i e^{i\rho_{\infty}} \phi(x+y_{\infty}-y, \alpha) \\  ie^{-i\rho_{\infty}} \phi(x+y_{\infty}-y, \alpha)\epm.
\ea\end{multline}
Finally, a simple computation shows that
\be
\frak g_{\pi^0} N^Z(\frak g_{\pi^0}^{-1} U^0, \pi^0) = N(U^0, \pi^0)
\ee
and thus we have retrieved all the terms of equation (\ref{nlsu}).
\end{proof}

In the next three lemmas we examine in more detail the properties of $\xi_F$ and $\eta_F$.

\begin{lemma}
$\eta_F$ and $\xi_F$ are biorthogonal, in the sense that
\be\ba{rl}
\langle \eta_{\alpha}, \xi_{G} \rangle &= -\frac 1 2 \alpha^{-1} \|\phi\|_2^2 \text{ for } G = \Gamma \text{ and } 0 \text{ otherwise }\\
\langle \eta_{\Gamma}, \xi_{G} \rangle &= \frac 1 2 \alpha^{-1} \|\phi\|_2^2 \text{ for } G = \alpha \text{ and } 0 \text{ otherwise }\\
\langle \eta_{v_k}, \xi_{G} \rangle &= -2\|\phi\|_2^2 \text{ for } G = D_k \text{ and } 0 \text{ otherwise }\\
\langle \eta_{D_k}, \xi_{G} \rangle &= 2\|\phi\|_2^2 \text{ for } G = v_k \text{ and } 0 \text{ otherwise.}
\ea\ee
\lb{orto}
\end{lemma}
Note that all these functions are related to $\pi^0$, not to $\pi$.
\begin{proof}
By direct computation. Most of the integrals cancel simply as the product of even and odd functions. The only nontrivial computation is that
\be
\langle \eta^F_{\alpha}, \xi^F_{\Gamma} \rangle = \int_{\set R^3} 2 \phi(x, \alpha) \partial_{\alpha} \phi(x, \alpha) \dd x = \partial_{\alpha} \|\phi(\cdot, \alpha)\|_2^2 = -\frac 1 2 \alpha^{-1} \|\phi(\cdot, \alpha)\|^2.
\ee
\end{proof}

Let
\be
E = \bpm -|v(t) - v_{\infty}|^2 - 2i(v(t)-v_{\infty})\dl & 0 \\ 0 & |v(t) - v_{\infty}|^2 - 2i(v(t)-v_{\infty})\dl\epm
\ee
and $\mc H = \mc H_{\pi^0} + V_{\pi^0}$.

\begin{lemma}
\be\ba{l}
\mc H^* \xi_{\alpha} = E \xi_{\alpha}\\
\mc H^* \xi_{\Gamma} = -2i \xi_{\alpha} + E \xi_{\Gamma}\\
\mc H^* \xi_{v_k} = E \xi_{v_k}\\
\mc H^* \xi_{D_k} = -2i\xi_{v_k} + E \xi_{D_k}.
\ea\lb{Hstea}\ee
\end{lemma}
\begin{proof} By direct computation.

For $t=\infty$ we note that $\xi_F$ actually become generalized eigenvectors for $\mc H^*$, because the symmetry transformation $\frak g_{\pi^0}$ was chosen so that $\frak g_{\pi^0} \bpm W(\pi^0) \\ \ov{W(\pi^0)} \epm$ becomes a ground state in the limit. Otherwise, there is a small error term.
\end{proof}

\begin{lemma}
If $U$ is a solution of (\ref{nlsu}), $\pi$ satisfies the modulation equations
\be\begin{aligned}
\dot \alpha &= 2\alpha\|\phi\|_2^{-2} (-\partial_t \langle U, \xi_{\alpha} \rangle + \langle U, \dot \xi_{\alpha}\rangle - i\langle U, E\xi_{\alpha}\rangle - i\langle N(U^0, \pi^0), \xi_{\alpha}\rangle)\\
\dot \Gamma &= 2\alpha\|\phi\|_2^{-2} (-\partial_t \langle U, \xi_{\Gamma} \rangle + 2i\langle U, \xi_{\alpha}\rangle + \langle U, \dot \xi_{\Gamma}\rangle - i\langle U, E\xi_{\Gamma}\rangle - i\langle N(U^0, \pi^0), \xi_{\Gamma}\rangle)\\
\dot v_k &= \|\phi\|_2^{-2} (-\partial_t \langle U, \xi_{v_k} \rangle + \langle U, \dot \xi_{v_k}\rangle - i\langle U, E\xi_{v_k}\rangle - i\langle N(U^0, \pi^0), \xi_{v_k}\rangle)\\
\dot D_k &= \|\phi\|_2^{-2} (-\partial_t \langle U, \xi_{D_k} \rangle + 2i\langle U, \xi_{v_k}\rangle + \langle U, \dot \xi_{D_k}\rangle - i\langle U, E\xi_{D_k}\rangle - i\langle N(U^0, \pi^0), \xi_{D_k}\rangle).
\lb{mod}\end{aligned}\ee
\lb{lemma25}
\end{lemma}
\begin{proof}
Consider the Schr\"odinger equation (\ref{nlsu}) and take its dot product with $\xi_F$ for each $F$:
\be
i \langle U_t, \xi_F \rangle + \langle U, \mc H^* \xi_F \rangle = -i\langle \dot \pi \partial_{\pi}\frak g_{\pi^0} \bpm W(\pi^0) \\ \ov{W(\pi^0)} \epm, \xi_F \rangle + \langle N(U^0, \pi^0), \xi_F \rangle.
\ee
Write $\langle U_t, \xi_F \rangle = \partial_t \langle U, \xi_F \rangle - \langle U, \dot \xi_F \rangle$. Lemma (\ref{Hstea}) helps evaluate the second term and the identity (\ref{pi}) and Lemma \ref{orto} help evaluate the third term, leading to the equations (\ref{mod}) above.
\end{proof}

We restate the orthogonality condition (\ref{ort_Z}) for $U$ in the form
\be
\langle U(t), \xi_{F}(t) \rangle = 0 \text{ for } F \in \{\alpha, \Gamma, v_k, D_k\},
\lb{ort_U}
\ee
for every $t \geq 0$, which is equivalent to condition (\ref{ort_Z}) on $Z$. Indeed, they follow from one another by applying the unitary transformations $\frak g_{\pi^0}(t)$.

Let
\be
L_{\pi^0} U = \sum_{F \in \{\alpha, \Gamma, v_i, D_i\}} \|\phi\|_2^{-2} (\langle U, \dot \xi_F\rangle - i\langle U, E\xi_F\rangle) \eta_F
\ee
and
\be
N_{\pi^0}(U^0, \pi^0) = - \sum_{F \in \{\alpha, \Gamma, v_i, D_i\}} \|\phi\|_2^{-2} i\langle N(U^0, \pi^0), \xi_F\rangle \eta_F.
\ee
The modulation equations can then be rewritten as
\be
-i \dot \pi \partial_{\pi} \frak g_{\pi^0} W(\pi^0) = L_{\pi^0} U + N_{\pi^0}(U^0, \pi^0).
\lb{modula}
\ee
$L_{\pi^0} U$ represents the part that is linear in $U$ and $N_{\pi^0}(U^0, \pi^0)$ represents the nonlinear part $\langle N(U^0, \pi^0), \xi_{F}\rangle$.

Note that the orthogonality condition for all times $t$ is equivalent to the modulation equations (\ref{modula}) together with the orthogonality condition at time $0$.

Next, we estimate a few useful quantities that appear in the right-hand side terms of the equations. Let
\be
\nu(T) = \|\langle t \rangle\dot\pi(t)\|_{L^1(T, \infty)} + \|\dot\pi\|_{L^{\infty}(T, \infty)}
\ee
and likewise $\nu^0$ for $\pi^0$. We still assume that $\|\dot\pi\|_{1\cap\infty} + \|t\dot\pi\|_1\leq\nu(0)<\infty$ is bounded and we justify this assumption later.

Now we state very general estimates that are used in the proof:
\begin{lemma} For any path $\pi$,
\be\begin{aligned}
\ds\int_T^{\infty} |v(t)-v_{\infty}| \dd t &\leq \|t \dot v(t)\|_{L^1(T, \infty)} \leq \nu(T),\\
T|v(T)-v_{\infty}| &\leq \|t \dot v(t)\|_{L^1(T, \infty)} \leq \nu(T),\\
\ds\int_T^{\infty} |v(t)-v_{\infty}|^2 \dd t &\leq \|\dot v(t)\|_{L^1(T, \infty)} \|t \dot v(t)\|_{L^1(T, \infty)} \leq \nu^2(T),\\
\ds\int_T^{\infty} t |v(t)-v_{\infty}|^2 \dd t &\leq \|t \dot v(t)\|_{L^1(T, \infty)}^2 \leq \nu^2(T).
\end{aligned}\ee
\lb{lemma26}
\end{lemma}
\begin{proof} All of these estimates are straightforward.
\end{proof}

Then, there are some more specific estimates that we need:
\begin{lemma}
\be\begin{aligned}
|\dot\rho_{\infty}| \leq& |\dot\Gamma| + |x| |\dot v| + (|v|+|v_{\infty}|) |v-v_{\infty}| + |\alpha-\alpha_{\infty}|\\
|\dl \rho_{\infty}| \leq& |v-v_{\infty}| \\
|\dot y - \dot y_{\infty}| \leq& |\dot D| + 2|v-v_{\infty}| \\
|\dot \xi_F + i E\xi_F| \leq& |\dot\Gamma| + |\dot D| + (1+|x|) |\dot v| + |v-v_{\infty}| + |\dot\alpha| + \\
&+|\alpha-\alpha_{\infty}|\\
\|V_{\pi^0}\|_{1+\infty \to 1\cap\infty} \leq& (|y-y_{\infty}| + |\rho_{\infty}| + |\alpha-\alpha_{\infty}|) \leq \nu^0(t)\\
\|\dl V_{\pi^0}\|_{1+\infty \to 1\cap\infty} \leq& |v-v_{\infty}| \leq \nu^0(t)\\
\|L_{\pi^0}\|_{1+\infty \to 1\cap\infty} \leq& C\|\langle \cdot, \dot \xi_F\rangle - i\langle \cdot, E\xi_F\rangle\| \\
\leq& C(|\dot y- \dot y_{\infty}| + |\dot \rho_{\infty}| + |\dot \alpha|) \leq \nu^0(t)\\
\|\dl L_{\pi^0} U(t)\|_{1 \cap \infty} \leq& C|\langle U(t), \dot \xi_F\rangle - i\langle U, E\xi_F\rangle|(1 + |v-v_{\infty}|).
\end{aligned}\ee
\lb{lemma27}
\end{lemma}
\begin{proof}
\begin{multline}
\rho_{\infty}(t) = \theta(t, x+y_{\infty}) - \theta_{\infty}(t, x+y_{\infty}) \\
= (v(t) - v(\infty)) (x + 2tv(\infty) + D(\infty)) - \int_0^t (|v(s)|^2 - |v(\infty)|^2 - \alpha(s) + \alpha(\infty)) \dd s + \\
+ \gamma(t) - \Gamma_{\infty}\\
= (v(t)-v(\infty)) (x + 2tv(\infty) + D(\infty)) + \int_t^{\infty} (|v(s)|^2 - |v(\infty)|^2 - \alpha(s) + \alpha(\infty)) \dd s + \\
+ \gamma(t) - \Gamma(\infty).
\end{multline}
Therefore,
\be
\dot \rho_{\infty}(t) = \dot v(t) x - (|v(t)|^2-|v(\infty)|^2) + \alpha(t) - \alpha(\infty) + \dot \Gamma(t).
\ee
Finally, note that
\be\begin{aligned}
\dot \xi_F = &i\dot\rho \sigma_3 \xi_F + (\dot y_{\infty} - \dot y) \dl \xi_F + \dot \alpha \partial_{\alpha} \xi_F\\
= &i(\dot v(t) x - |v(t)|^2 + |v(\infty)|^2 + \alpha(t) - \alpha(\infty) + \dot \Gamma(t)) \sigma_3 \xi_F + \\
&+ 2(v(\infty) - v(t)) \dl\xi_F + \dot \alpha \partial_{\alpha} \xi_F,
\end{aligned}\ee
so
\be\begin{aligned}
\dot \xi_F + i E \xi_F = & i(\dot v(t) x - |v(t)|^2 + |v(\infty)|^2 + \alpha(t) - \alpha(\infty) + \dot \Gamma(t)) \sigma_3 \xi_F + \\
&+ i|v(t)-v_{\infty}|^2 \sigma_3 \xi_F + \dot \alpha \partial_{\alpha} \xi_F\\
= &i(\dot v(t) x + 2v_{\infty} (v-v_{\infty}) + \alpha(t) - \alpha(\infty) + \dot \Gamma(t)) \sigma_3 \xi_F + \dot \alpha \partial_{\alpha} \xi_F.
\end{aligned}\ee
The other formulae follow by straightforward computations.
\end{proof}

\subsection{Spectrum of the Hamiltonian} \lb{spectru}

Without loss of generality, we perform a symmetry transformation in the nonlinear equation (\ref{NLS}) and assume that the initial data is in the neighborhood of a positive ground state of the equation $W(\pi(0)) = \phi(\cdot, \alpha)$  instead of a more general standing wave. This is possible because standing waves are, by definition, an orbit of the action of symmetry transformations. Even though symmetry transformations change the spectrum, the information gained in the manner is still useful in the general case.

Consider the operators 
\be
\mc H = \bpm -\Delta + \alpha - 2\phi^2(\cdot, \alpha) & -\phi^2(\cdot, \alpha) \\ \phi^2(\cdot, \alpha) & \Delta - \alpha + 2\phi^2(\cdot, \alpha)\epm = \mc H_0 + V.
\ee
By rescaling, one sees that all these operators have the same spectrum up to dilation and similar spectral properties.

We restate the known facts about the spectrum of $\mc H$. As proved by Buslaev, Perelman \cite{buslaev1} and also Rodnianski, Schlag, Soffer in \cite{rod3}, under fairly general assumptions, $\sigma(\mc H) \subset \set R \cup i\set R$ and is symmetric with respect to the coordinate axes and all eigenvalues are simple with the possible exception of $0$. Furthermore, by Weyl's criterion $\sigma_{ess}(\mc H) = (-\infty, -\alpha] \cup [\alpha, +\infty)$.

Grillakis, Shatah, Strauss \cite{gril1} and Schlag \cite{schlag} showed that there is only one pair of conjugate imaginary eigenvalues  $\pm i\sigma$ and that the corresponding eigenvectors decay exponentially. For the decay see Hundertmark, Lee \cite{hund}. The pair of conjugate imaginary eigenvalues $\pm i\sigma$ reflects the $L^2$-supercritical nature of the problem.

The generalized eigenspace at $0$ arises due to the symmetries of the equation, which is invariant under Galilean coordinate changes, phase changes, and scaling. It is relatively easy to see that each of these symmetries gives rise to a generalized eigenvalue of the Hamiltonian $\mc H$ at $0$, but proving the converse is much harder and was done by Weinstein in \cite{wein1}, \cite{wein2}.

Schlag \cite{schlag} showed, using ideas of Perelman \cite{perel}, that if the operators
\be
L_{\pm} = -\Delta + \alpha - \phi^2(\cdot, \alpha) \mp 2 \phi^2(\cdot, \alpha)
\ee
that arise by conjugating $\mc H$ with $\bpm 1 & i \\ 1 & -i \epm$ have no eigenvalue in $(0, \alpha]$ and no resonance at $\alpha$, it implies that the real discrete spectrum of $\mc H$ is $\{0\}$ and that the edges $\pm \alpha$ are neither eigenvalues nor resonances. A paper of Demanet, Schlag \cite{demanet} proved numerically that the scalar operators meet these conditions. Therefore, there are no eigenvalues in $[-\alpha, \alpha]$ and $\pm \alpha$ are neither eigenvalues nor resonances for $\mc H$.

Furthermore, the method of Agmon \cite{agmon}, adapted to the matrix case, enabled Erdogan, Schlag \cite{erdogan2} and independently \cite{cuc2} to prove that any resonances embedded in the interior of the essential spectrum (that is, in $(-\infty, -\alpha) \cup (\alpha, \infty)$) have to be eigenvalues, under very general assumptions.

Under the spectral Assumption \ref{assum} we now have a complete description of the spectrum of $\mc H$. It consists of a pair of conjugate purely imaginary eigenvalues, a generalized eigenspace at $0$, and the essential spectrum $(-\infty, -\alpha] \cup [\alpha, \infty)$.

It helps in the proof to exhibit the discrete eigenspaces of $\mc H$. Denote by $f^{\pm}$ and $\tilde f^{\pm}$ the normalized eigenfunctions of $\mc H$ and respectively $\mc H^*$ corresponding to the $\pm i \sigma$ eigenvalues. Also observe that $\eta_F$ are the generalized eigenfunctions at zero of $\mc H$ and $\xi_F$ fulfill the same role for $\mc H^*$.

Furthermore, now we can express the Riesz projections, following Schlag \cite{schlag}, as
\be
P_{im} = P_+ + P_-,\ P_{\pm} = \langle \cdot, \tilde {f^{\pm}} \rangle f^{\pm},
\ee
\be
P_{root} = \langle \cdot, \xi_{\alpha} \rangle \eta_{\Gamma} + \langle \cdot, \xi_{\Gamma} \rangle \eta_{\alpha} + \sum_k (\langle \cdot, \xi_{v_k} \rangle \eta_{D_k} + \langle \cdot, \xi_{D_k} \rangle \eta_{v_k}),
\ee
and
\be
P_c = 1-P_{im}-P_{root}.
\ee
Even though we do not have an explicit form of the imaginary eigenvectors, Schlag \cite{schlag} proved that $f^{\pm}$, in the $L^2$ norm, and $\sigma$ are locally Lipschitz continuous as a function of $\alpha$ and that $f^{\pm}$ are exponentially decaying.

Concerning the continuous spectrum, the absence of embedded eigenvalues, following the spectral Assumption \ref{assum}, permitted Erdogan, Schlag \cite{erdogan2} to state the limiting absorbtion principle in the following form:
\begin{lemma}
Assume that the thresholds of the spectrum of $\mc H = \mc H_0 + V$ (\ref{ham}) are regular, meaning that the operators $1+(\mc H_0-(\pm \alpha \pm i0))^{-1} V$ are invertible from the weighted Sobolev space $\langle x \rangle^{1+\epsilon} L^2$ to itself for any $\epsilon>0$. Then there exists $0<\alpha'<\alpha$ such that
\be
\sup_{|\lambda|\geq\alpha', \epsilon>0} |\lambda|^{1/2} \|(\mc H-(\lambda \pm i\epsilon))^{-1}\|_{\langle x \rangle^{-1-\epsilon} L^2 \to \langle x \rangle^{1+\epsilon} L^2} < \infty
\ee
and
\be
\sup_{\substack{|\lambda|\geq\alpha', \epsilon>0\\ \ell=1, 2}} \|\partial_{\lambda}^{\ell} (\mc H-(\lambda \pm i\epsilon))^{-1}\|_{\langle x \rangle^{-1-\ell-\epsilon} L^2 \to \langle x \rangle^{1+\ell+\epsilon} L^2} < \infty
\ee
if $|V| \leq C \langle x \rangle^{-7/2-\epsilon}$.
\lb{limabs}
\end{lemma}
The fact that the thresholds are neither eigenvalues nor resonances implies their regularity.

\subsection{Proof of the main result}
\begin{proof}
To recapitulate, we are interested in finding solutions to equation (\ref{NLS}), starting from initial data in a neighborhood of the soliton $W(0)$, which remain close to stationary solutions for all times. Furthermore, we take $W(0)$ to be a positive ground state of the equation.

We prove that, to a first approximation, $\Big\{R_0 \mid (P_+ + P_{root}) \bpm R_0 \\ \overline R_0\epm = 0\Big\}$ is the stable submanifold. Quadratic corrections are needed, as the statement of Theorem~\ref{main} makes clear.

Let, for some $\delta\leq 1$, $1 \leq q<4/3$, small $\epsilon>0$,
\be\begin{split}
X_{\delta} = \big\{(U, \pi) \mid &\pi(0) = (0, 0, 0, \alpha),\ \|\langle t \rangle\dot\pi(t)\|_{1} + \|\dot\pi\|_{\infty} \leq \delta,\\
&\|U\|_{L^2_t W^{1/2,6}_x \cap L^{\infty}_t H^{1/2}_x} +  \|U\|_{\langle t \rangle^{1-2/q+\epsilon} L^{2}_t L^{6+\infty}_x} \leq \delta\big\}.
\lb{defx}
\end{split}\ee

Let $\beta = 2/q-1-\epsilon$. Note that, for $U(x, t) = \frak g_{\pi^0}(t) Z(x, t)$, one has $(U, \pi) \in X_{\delta}$ if and only if $(Z, \pi) \in X_{\delta}$, since $\frak g_{\pi^0}(t)$ is an isometry. In the sequel we deal with both $U$ and $Z$, as necessary.

We define a mapping $\Phi$ that takes the pair $(Z^0, \pi^0)$ to the solution $(Z, \pi)$ of the linearized equation (\ref{Z}) corresponding to the initial data
\be\begin{aligned}
Z(0) = & \bpm R_0\\ \overline{R_0}\epm + h f^{+}(0),\\
\pi(0)= & (0, 0, 0, \alpha),
\lb{Z_initial}
\end{aligned}\ee
where $h$ and $a_F(0)$ will be chosen later depending on $R_0$, so that $Z$ fulfills the orthogonality condition (\ref{ort_U}) and is globally bounded in time. The first condition can be equivalently formulated in terms of $U$,
\be
U(0) = \frak g_{\pi^0}(0) \Big(\bpm R_0\\ \overline{R_0}\epm + h f^{+}(0) \Big).
\lb{U_initial}
\ee
Further note that, by interpolation between $L^2$ and $H^1$, for $\delta<1$
\be
\|U(0)\|_{L^q \cap H^{1/2}} \leq C(\|R_0\|_{L^q \cap H^{1/2}}+|h|).
\ee

Now we run a fixed point argument, showing that $\Phi$ is a contraction in $X_{\delta}$ for small $\delta$. This is achieved in two steps, by proving firstly that $\Phi$ takes $X_{\delta}$ to itself and secondly that it is distance-decreasing in a weaker metric.

\subsection{Stability}

Here we prove that if $\delta \leq 1$ is sufficiently small and $(Z^0, \pi^0) \in X_{\delta}$, then $\Phi(Z^0, \pi^0) = (Z, \pi) \in X_{\delta}$.

Since $\frak g$ is an isometry, this is equivalent to proving that if $(U^0, \pi^0) \in X_{\delta}$, then the solution $(U, \pi)$ of (\ref{nlsu}) is in $X_{\delta}$, for small $\delta$. It is more convenient to prove the statement for $U$ than for $Z$.

Note that, after making $\delta$ as small as needed, $\alpha(t)$ always belongs to a fixed compact set centered at its initial value and therefore all the Sobolev norms of $\phi(\cdot, \alpha(t))$ and $f^{\pm}$, $\tilde f^{\pm}$ are uniformly bounded.

Replace $\dot\pi$ on the right-hand side of (\ref{nlsu}) by its expression given by (\ref{pi}) and assume the orthogonality condition $\langle U(t), \xi_{F}(t)\rangle = 0$ in order to obtain the system of equations
\be\begin{aligned}
iU_t + \mc H_{\pi^0} U &= L_{\pi^0} U + V_{\pi^0} U + N(U^0, \pi^0) + N_{\pi^0}(U^0, \pi^0)\\
\dot F &= 2\alpha_0 \|\phi_0\|_2^{-2} (\langle U, \dot \xi_{F_0}\rangle - i\langle U, E_0\xi_{F_0}\rangle - i\langle N(U^0, \pi^0), \xi_{F_0}\rangle),\  F \in \{\alpha, \Gamma\}\\
\dot F &= \|\phi_0\|_2^{-2} (\langle U, \dot \xi_{F_0}\rangle - i\langle U, E_0\xi_{F_0}\rangle - i\langle N(U^0, \pi^0), \xi_{F_0}\rangle),\ F\in \{D_k, v_k\}.
\end{aligned}\lb{ec_liniara}\ee
Here we made the replacement $-i\dot \pi \frak g_{\pi^0} \partial_{\pi^0} \bpm W(\pi^0) \\ \ov{W(\pi^0)} \epm = L_{\pi^0} U + N_{\pi^0}(U^0, \pi^0)$ by virtue of the modulation equations (\ref{mod}), (\ref{modula}). $L_{\pi^0} U$ represents the part that is linear in $U$ and $N_{\pi^0}(U^0, \pi^0)$ represents the nonlinear part $\langle N(U^0, \pi^0), \xi_{F_0}\rangle$.

The initial data is given by condition (\ref{Z_initial}).

The orthogonality condition at time $0$ is true by definition, regardless of the value of $R_0 \in S_{\delta}$ and $h$, because $\bpm R_0 \\ \ov{R_0} \epm \perp \xi_F(0)$ and, due to the spectral decomposition, $f^+(0) \perp \xi_F(0)$ too, for $F \in \{v_k, D_k, \alpha, \Gamma\}$.

Global existence of the solution $U$ to the linearized equation (\ref{ec_liniara}) follows by a standard fixed point argument. Introduce a second auxiliary function $U_1$ and write the equation as
\be
iU_t + \mc H_{\pi^0} U = V_{\pi^0} U_1 + L_{\pi^0} U_1 + N(U^0, \pi^0) + N_{\pi^0}(U^0, \pi^0).
\ee
Note that $\|e^{it \mc H_{\pi^0}}\|_{2 \to 2} \leq C e^{\sigma|t|}$.

For any $T_1$ consider a small time interval $[T_1, T_2]$ of length at most $1$. Assume that $\|U_1\|_{L^{\infty}(T_1, T_2; L^2_x)} \leq r$. One has that
\begin{multline}
\|U\|_{L^{\infty}(T_1, T_2; L^2_x)} \leq \\
\begin{aligned}
\leq &\|U(T_1)\|_2 + \int_{T_1}^{T_2} e^{\sigma(t-T_1)} \|RHS(t)\|_2 \dd t\\
\leq &\|U(T_1)\|_2 + C(T_2-T_1) (\|L_{\pi^0}+V_{\pi^0}\|_{L^2_x \to L^2_x} \|U_1\|_{L^{\infty}(T_1, T_2; L^2_x)} + \\
&+ \|N(U^0, \pi^0) + N_{\pi^0}(U^0, \pi^0)\|_{L^{\infty}_t L^2_x}) \leq r,
\end{aligned}
\end{multline}
if $r$ is chosen such that $r \geq C (\|U(T_1)\|_2 + \delta)$. Likewise, subtracting two copies of the equation, with two solutions $U$ and $\tilde U$ corresponding to auxiliary functions $U_1$ and $\tilde U_1$, one obtains
\be
\|U-\tilde U\|_{L^{\infty}(T_1, T_2; L^2_x)} \leq C(T_2-T_1) \|L_{\pi^0}+V_{\pi^0}\|_{L^2_x \to L^2_x} \|U_1-\tilde U_1\|_{L^{\infty}(T_1, T_2; L^2_x)}.
\ee
Thus, if $T_2-T_1$ is sufficiently small, the mapping that associates $U$ to $U_1$ is a contraction in the set $\{U_1 \mid \|U_1\|_{L^{\infty}(T_1, T_2; L^2_x)} \leq r\}$. If $r$ is sufficiently large the set is stable under the mapping.

One obtains a fixed point that is a solution to (\ref{ec_liniara}) on $(T_1, T_2)$, but the length $T_2-T_1$ for which this happens does not depend on the initial data.  Therefore, one can iterate and obtain a global in time solution $U$ of (\ref{ec_liniara}).

Next, we prove that the global solution $U$ is in $X_{\delta}$ and thus globally bounded for some unique value of the parameter $h$.

The operator $\mc H_{\pi^0}$ induces the time-independent decomposition $1=P_c + P_{root} + P_{im}$ on $L^2(\set R^3)$ corresponding to the decomposition of its spectrum into the absolutely continuous part, the generalized eigenspace at zero, and the imaginary eigenvalues, respectively. Since the range and cokernel of $P_{root}$ and $P_{im}$ are spanned by finitely many Schwartz functions, they are bounded from $L^p$ to $L^q$, for any $1 \leq p, q \leq \infty$. Therefore $P_c=1-P_{root}-P_{im}$ is bounded on $L^p$, $1\leq p \leq \infty$, and one can write
\be
P_{root} U(t) = \sum_F a_F(t) \eta_F(\infty),\ P_{im} U(t) = b^+(t) f^+ + b^-(t) f^-.
\lb{partea_finita}
\ee
We will bound each projection $P_{root} U$, $P_{im} U$, and $P_c U$, as well as $\dl P_{root} U$, $\dl P_{im} U$, and $\dl P_c U$ (six estimates in total).

One can bound the zero generalized eigenspace component in a straightforward manner. Expanding the orthogonality condition $\langle U(t), \xi_{F_0}(t)\rangle = 0$, one has for every $G$ that
\be
0 = \sum_F a_F(t) \langle \eta_{F}(\infty), \xi_{G}(t)\rangle + \langle (P_{im} + P_c) U(t), \xi_G(t)\rangle.
\lb{sist}
\ee
Since $|\xi_{G_0}(t)-\xi_{G_0}(\infty)| \leq C (|\rho_{\infty}| + |y-y_{\infty}|) \leq C \delta$ and the matrix with entries $\langle \eta_F(\infty), \xi_G(\infty)\rangle$ is invertible, the matrix with entries $\langle\eta_{F_0}(\infty), \xi_{G_0}(t)\rangle$ is invertible with bounded norm for small $\delta$. Therefore, by solving the system (\ref{sist}) one obtains that
\be
\|P_{root} U(t)\|_{1\cap\infty} \leq C \|(P_c+P_{im}) U(t)\|_{1+\infty}.
\label{root1}
\ee
Since the range of $P_{root}$ is spanned by Schwartz functions, the same holds with derivatives or weights:
\be
\|\dl P_{root} U(t)\|_{1\cap\infty} \leq C \|(P_c+P_{im}) U(t)\|_{1+\infty}.
\label{root2}
\ee

As for the other two components of $U$, one has that
\be
i\partial_t{P_c U} + \mc H_{\pi^0} P_c U = P_c (L_{\pi^0}(U) + V_{\pi^0} U + N(U^0, \pi^0) + N_{\pi^0}(U^0, \pi^0))
\lb{disp}
\ee
and
\be
i\partial_t{P_{im} U} + \mc H_{\pi^0} P_{im} U = P_{im} (L_{\pi^0}(U) + V_{\pi^0} U + N(U^0, \pi^0) + N_{\pi^0}(U^0, \pi^0)).
\lb{hyperb}
\ee

Using the explicit form (\ref{partea_finita}) of $P_{im} U(t) = b_-(t) f^- + b_+(t) f^+$, the corresponding equation (\ref{hyperb}) becomes
\be
\partial_t  \bpm b_-\\ b_+\epm + \bpm \sigma & 0 \\ 0 & -\sigma\epm \bpm b_-\\ b_+\epm = \bpm N_-\\ N_+\epm,
\lb{hiper}
\ee
where $|N_{\pm}(t)| \leq \|P_{im}(L_{\pi^0}(U) + V_{\pi^0} U + N(U^0, \pi^0)(t) + N_{\pi^0}(U^0, \pi^0))\|_{1 + \infty}$.
Here $\pm i\sigma$ are the imaginary eigenvalues of $\mc H_{\pi^0}$, as in our discussion of its spectrum in Section \ref{spectru}.

Now we state a standard elementary lemma, see \cite{schlag}. It characterizes the unique bounded solution of the two-dimensional ODE (\ref{hiper}).

\begin{lemma}
Consider the equation
\be
\dot x - \bpm \sigma & 0 \\ 0 & -\sigma\epm x = f(t),
\ee
where $f \in L^{1 \cap \infty}$. Then $x$ is bounded on $[0, \infty)$ if and only if
\be
0 = x_1(0) + \int_0^{\infty} e^{-t \sigma} f_1(t) \dd t.
\lb{condi}
\ee
In this case,
\be
x_1(t) = -\int_t^{\infty} e^{(t-s)\sigma} f_1(s) \dd s,\ x_2(t) = e^{-t \sigma} x_2(0) + \int_0^t e^{-(t-s)\sigma} f_2(s) \dd s
\ee
for all $t \geq 0$.
\lb{hyp}
\end{lemma}
\begin{proof}
Any solution will be a linear combination of the exponentially increasing and the exponentially decaying ones and we want to make sure that the exponentially increasing one is absent. It is always true that
\be
x_1(t) = e^{t \sigma} \Big(x_1(0) + \int_0^t e^{-s\sigma} f_1(s) \dd s\Big),\ x_2(t) = e^{-t \sigma} x_2(0) + \int_0^t e^{-(t-s)\sigma} f_2(s) \dd s.
\ee
Thus, if $x_1$ is to remain bounded, the expression between parantheses must converge to $0$, hence (\ref{condi}). Conversely, if (\ref{condi}) holds, then
\be
x_1(t) = -\int_t^{\infty} e^{(t-s)\sigma} f_1(s) \dd s
\ee
tends to $0$.
\end{proof}

Consequently, equation (\ref{hiper}) has a bounded solution if and only if
\be
0 = b_+(0) + \int_0^{\infty} e^{-t \sigma} N_+(t) \dd t.
\lb{conditie}
\ee

Now we establish the relation between $b_+(0)$ and $h$. The initial assignment (\ref{U_initial}) implies that
\be\begin{split}
b_+(0) &= \langle U(0), \tilde f^+ \rangle\\
&= \Big\langle \frak g_{\pi^0}(0) \bpm R_0 \\ \ov{R_0} \epm + h \frak g_{\pi^0}(0) f^+(0), \tilde f^+ \Big\rangle \\
&= \Big\langle \bpm R_0 \\ \ov{R_0} \epm, (\frak g_{\pi^0}(0)^{-1} \tilde f^+) - \tilde f^+(0) \Big\rangle + h \big(1+ \langle (\frak g_{\pi^0} f^+(0)) - f^+, \tilde f^+ \rangle\big).
\lb{bplus}
\end{split}
\ee
Taking into account the fact that
\be
\|(\frak g_{\pi^0} f^+) - f^+(0)\|_{\langle x \rangle L^2_x} \leq \|(\frak g_{\pi^0} f^+) - f^+\|_{\langle x \rangle L^2_x} + \|f^+ - f^+(0)\|_{\langle x \rangle L^2_x} \leq C\delta,
\ee
it follows that if (\ref{conditie}) holds then, for sufficiently small $\delta$, one can solve equation (\ref{bplus}) for $h$ and
\be
|h| \leq C(|b_+(0)|+ \delta \|R_0\|_{1+\infty}).
\ee
Clearly, condition (\ref{conditie}) is then fulfilled by a suitable choice of $h$. $U$ is globally bounded by the definition (\ref{defx}) of $X_{\delta}$, which implies the boundedness of each component, in particular $U_{im}$. Proceeding henceforth under this assumption, we get
\be
|b_+(0)| \leq C\|N_+\|_{(L^1_t + L^{\infty}_t)(L^1_x+L^{\infty}_x)} \leq C\delta\|U\|_{L^{\infty}_t L^2_x \cap L^2_t L^6_x}+C\delta^2,
\ee
\be
|h| \leq C\delta(\|R_0\|_{1+\infty} + \|U\|_{L^{\infty}_t L^2_x \cap L^2_t L^6_x})+C\delta^2.
\lb{h}
\ee

Note that $\sigma$ depends Lipschitz continuously on $\alpha$. Then $\sigma$ belongs to a compact subset of $(0, \infty)$, because $\alpha$ belongs to a compact subset of $(0, \infty)$. In this case, since $f^{\pm}$ are Schwartz functions, one has that
\be\begin{aligned}
\|P_+ U(t)\| _{1 \cap \infty \cap H^1} \leq & C \Big(\int_t^{2t} e^{(t-s)\sigma}\|RHS(s)\|_{1+\infty} \dd s + \int_{2t}^{\infty} e^{(t-s)\sigma}\|RHS(s)\|_{1+\infty} \dd s \Big)\\
\leq & C \Big(\int_t^{2t} e^{(t-s)}\|RHS(s)\|_{1+\infty} \dd s + e^{-t\sigma}(\delta \|U\|_{L^2_t L^6_x} + \delta^2)\Big).
\end{aligned}\ee
We deal with the two expressions separately. The latter poses no problem, due to the exponential decay. As for the former, we bound it in $\langle t \rangle^{-\beta} L^2_t L^{1\cap\infty}_x \cap H^1_x$ (the precise norm in $x$ does not matter, due to these being Schwartz functions) by means of the bilinear estimate
\be
\Big|\iint_{t<s<2t} \langle t \rangle^{\beta} \langle s \rangle^{-\beta} e^{t-s} g(t) f(s) \dd s \dd t\Big| \leq \|f\|_2 \|g\|_2.
\ee
Indeed, note that $\langle t \rangle^{\beta} \langle s \rangle^{-\beta}$ is bounded from above and the conclusion follows after a dyadic decomposition. The estimate implies that
\be
\Big\|\int_t^{2t} e^{(t-s)}\|RHS(s)\|_{1+\infty} \dd s\Big\|_{\langle t \rangle^{-\beta} L^2_t L^{1\cap\infty}_x} \leq C \|RHS\|_{\langle t \rangle^{-\beta} L^2_t L^{1\cap\infty}_x}.
\ee
The pointwise in time norm is easier to bound, since
\be
\int_t^{2t} e^{(t-s)}\|RHS(s)\|_{1+\infty} \dd s \leq C \max_{s \in [t, 2t]} \|RHS(s)\|_{1+\infty} \leq C (\delta \|U\|_{L^{\infty}_t L^2_x} + \delta^2).
\ee
The same works for $P_- U$, where we also have to take into account the contribution of the initial data. This leads to
\begin{multline}
\|P_{im} U\|_{L^{\infty}_t L^2_x \cap L^2_t L^6_x} + \|\dl P_{im} U\|_{L^{\infty}_t L^2_x \cap L^2_t L^6_x} \leq\\
\leq C(\|U(0)\|_{1+\infty} + \delta\|U\|_{L^{\infty}_t L^2_x \cap L^2_t L^6_x}+\delta^2).
\label{im}
\end{multline}
and
\be
\|\langle t \rangle^{\beta}P_{im} U\|_{L^2_t L^{1\cap\infty}_x} \leq C (\|U(0)\|_{1+\infty} + \delta\|U\|_{\langle t \rangle^{-\beta} L^2_t L^{6+\infty}_x}+\delta^2).
\lb{eq_imw}
\ee

Now we turn to $P_c U$, the projection on the continuous spectrum. One has, by Lemmas \ref{lemma26} and \ref{lemma27}, that
\be\begin{split}
&\|N(\pi^0, U)\|_{L^2_t L^{6/5}_x} \leq C(\|U\|_{L^{\infty}_t L^6_x} \|U\|_{L^{\infty}_t L^2_x} \|U\|_{L^2_t L^6_x} + \|U\|_{L^{\infty}_t L^2_x} \|U\|_{L^2_t L^6_x})\\
&\|V_{\pi^0} U\|_{L^2_t L^{6/5}_x} \leq C\nu^0(1+\nu^0) \|U(t)\|_{L^2_t L^6_x}\\
&\|L_{\pi^0} U\|_{L^2_t L^{6/5}_x} + \|\dl L_{\pi^0} U\|_{L^2_t L^{6/5}_x} \leq C\nu^0(1+\nu^0) \|U\|_{L^2_t L^6_x}\\
&\|N_{\pi^0}(U, \pi^0)\|_{L^2_t L^{6/5}_x} + \|\dl N_{\pi^0}(U, \pi^0)\|_{L^2_t L^{6/5}_x} \leq \\
&\leq C(\|U\|_{L^{\infty}_t L^6_x} \|U\|_{L^{\infty}_t L^2_x} \|U\|_{L^2_t L^6_x} + \|U\|_{L^{\infty}_t L^2_x} \|U\|_{L^2_t L^6_x}).
\end{split}\ee
We recall that, by the definition (\ref{defx}) of $X_{\delta}$, $\nu^0 \leq \delta$.

Applying endpoint Strichartz estimates to equation (\ref{disp}) yields that
\be\begin{aligned}
\|P_c U\|_{L^2_t L^6_x \cap L^{\infty}_tL^2_x} &\leq C(\|U(0)\|_2 + \|P_c (L_{\pi^0}(U) + V_{\pi^0} U + N(U^0, \pi^0))\|_{L^2_t L^{6/5}_x + L^1_tL^2_x}) \\
&\leq C(\|U(0)\|_2 + \delta\|U\|_{L^{\infty}_t L^2_x \cap L^2_t L^6_x}) +C\delta^2.
\label{cont1}
\end{aligned}\ee

We now establish the $H^{1/2}$ bounds for $P_c U$. Note that the components of $[\dl, \mc H]$ are Schwartz functions and therefore $\|[\dl, \mc H_{\pi^0}] U\|_{L^2_t L^{6/5}_x} \leq C \|U\|_{L^2_t L^6_x}$. By trilinear interpolation between
\be\begin{aligned}
\|P_c U\|_{L^{\infty}_t L^2_x \cap L^2_t L^6_x} &\leq C(\|U(0)\|_2 + \|i (P_c U)_t + \mc H P_c U\|_{L^1_t L^2_x + L^2_t L^{6/5}_x}) \\
&\leq C(\|U(0)\|_2 + \|RHS\|_{L^2_t L^{6/5}_x + L^1_tL^2_x})
\end{aligned}\ee
and
\begin{multline}
\|P_c U\|_{L^{\infty}_t H^1_x \cap L^2_t W^{1, 6}_x} \leq C(\|U(0)\|_{H^1_x} + \|i (P_c U)_t + \mc H P_c U\|_{L^1_t H^1_x + L^2_t W^{1, 6/5}_x}) \\
\leq C(\|U(0)\|_{H^1_x} + \|RHS\|_{L^2_t W^{1, 6/5}_x + L^1_tH^1_x} + \|U\|_{L^2_t L^6_x})
\end{multline}
it follows that
\begin{multline}
\|P_c U\|_{L^{\infty}_t H^{1/2}_x \cap L^2_t W^{1/2, 6}_x}\leq C(\|U(0)\|_{H^{1/2}_x} + \|RHS\|_{L^2_t W^{1/2, 6/5}_x + L^1_tH^{1/2}_x} + \|U\|_{L^2_t L^6_x}).
\end{multline}

The fractional Sobolev spaces $H^s$ and $W^{s, p}$ arise naturally by interpolation and are given by $H^s = W^{s, 2}$, $W^{s, p} = (-\Delta+1)^{-s/2} L^p$. In the sequel we use the Sobolev embeddings $W^{1/2, 3} \hookrightarrow L^6$ and $H^{1/2} \hookrightarrow L^3$.

Now we examine each term on the right-hand side of (\ref{nlsu}). We use the fractional Leibniz rule, as stated, for example, in \cite[p. 105]{taylor}:
\be\begin{aligned}
\|f_1 f_2 f_3\|_{W^{1/2, 6/5}} &\leq C (\sum_{i, j, k} \|f_i\|_{W^{1/2, 3}} \|f_j\|_{3} \|f_k\|_{6} + \sum_{i, j, k} \|f_i\|_{3} \|f_j\|_3 \|f_k\|_{W^{1/2, 6}})\\
&\leq C (\sum_{i, j, k} \|f_i\|_{W^{1/2, 3}} \|f_j\|_{3} \|f_k\|_{W^{1/2, 3}} + \sum_{i, j, k} \|f_i\|_{3} \|f_j\|_3 \|f_k\|_{W^{1/2, 6}}).
\lb{frac_L}
\end{aligned}\ee
Making all the $f$'s equal, one gets
\be\begin{aligned}
\|U^3\|_{L^2_t W^{1/2, 6/5}_x} &\leq C(\|U\|_{L^4_t W^{1/2, 3}_x}^2 \|U\|_{L^{\infty}_t L^3_x} + \|U\|_{L^2_t W^{1/2, 6}_x} \|U\|_{L^{\infty}_t L^3_x}^2)\\
&\leq C \|U\|_{L^{\infty}_t H^{1/2}_x \cap L^2_t W^{1/2, 6}_x}^3.
\end{aligned}\ee
Thus, in estimating this cubic term we had to use both the Keel-Tao endpoint Strichartz estimate and the critical half-derivative, which prevents us from doing any better (that is, from lowering the number of derivatives).

The localized quadratic terms can be handled similarly with the help of (\ref{frac_L}), the conclusion being
\be\begin{aligned}
\|U^2 \phi\|_{L^2_t W^{1/2, 6/5}_x} &\leq C(\|U\|_{L^4_t W^{1/2, 3}_x}^2 \|\phi\|_{L^{\infty}_t L^3_x} + \|U\|_{L^2_t W^{1/2, 6}_x} \|U\|_{L^{\infty}_t L^3_x} \|\phi\|_{L^{\infty}_t L^3_x})\\
&\leq C \|U\|_{L^{\infty}_t H^{1/2}_x \cap L^2_t W^{1/2, 6}_x}^2.
\end{aligned}\ee
As for the linear terms, a satisfactory estimate is
\be
\|V_{\pi^0} U\|_{W^{1/2, 6/5}_x} \leq C(\|V_{\pi^0}\|_{3/2} \|U\|_{W^{1/2, 6}} + \|V_{\pi^0}\|_{W^{1/2, 3/2}_x} \|U\|_6),
\ee
which implies that
\be
\|V_{\pi^0} U\|_{L^2_t W^{1/2, 6/5}_x} \leq C \|V_{\pi^0}\|_{L^{\infty}_t W^{1/2, 3/2}_x} \|U\|_{L^2_t W^{1/2, 6}_x}.
\ee
Therefore
\be\begin{aligned}
\|P_c U\|_{L^2_t W^{1/2, 6}_x \cap L^{\infty}_t H^{1/2}_x} \leq &C(\|U(0)\|_{H^{1/2}_x} + \\
&+ \delta\|U\|_{L^2_t W^{1/2, 6}_x \cap L^{\infty}_t H^{1/2}_x} +\delta^2 + \|U\|_{L^2_t L^6_x}).
\lb{cont2}
\end{aligned}\ee

Putting together estimates (\ref{root1}), (\ref{im}), and (\ref{cont1}), one has that
\be
\|U\|_{L^{\infty}_t L^2_x \cap L^2_t L^6_x} \leq C(\|U(0)\|_2 + \delta\|U\|_{L^{\infty}_t L^2_x \cap L^2_t L^6_x} + \delta^2),
\ee
which for sufficiently small $\delta$ and $\|U(0)\|_2$ implies that
\be
\|U\|_{L^{\infty}_t L^2_x \cap L^2_t L^6_x} \leq C \|U(0)\|_2 + C\delta^2 \leq \delta.
\lb{est_U}
\ee
This proves, after considering (\ref{root2}), (\ref{im}), and (\ref{cont2}), that
\be
\|U\|_{L^2_t W^{1/2, 6}_x \cap L^{\infty}_t H^{1/2}_x} \leq C(\|U(0)\|_{H^{1/2}_x} + \delta\|U\|_{L^2_t W^{1/2, 6}_x \cap L^{\infty}_t H^{1/2}_x} + \delta^2 + \|U\|_{L^2_t L^6_x})
\ee
and therefore
\be
\|U\|_{L^2_t W^{1/2, 6}_x \cap L^{\infty}_t H^{1/2}_x} \leq C\|U(0)\|_{H^{1/2}_x} + C\delta^2 \leq\delta.
\lb{est_grad_U}
\ee

Lastly, for $\pi$ the estimates (see Lemma \ref{lemma27})
\begin{multline}
\|\langle U, \dot \xi_F\rangle - i\langle U, E\xi_F\rangle - i\langle N(U^0, \pi^0), \xi_F\rangle\|_{L^1_t} \leq \\
\leq \nu^0 \|U\|_{L^{\infty}_t L^6_x} + \|U^0\|_{L^{\infty}_t L^6_x}  \|U^0\|_{L^2_t L^6_x}^2 + \|U^0\|_{L^2_t L^6_x}^2
\end{multline}
and
\begin{multline}
\|\langle U, \dot \xi_F\rangle - i\langle U, E\xi_F\rangle - i\langle N(U^0, \pi^0), \xi_F\rangle\|_{L^{\infty}_t} \leq \\
\leq \nu^0 \|U\|_{L^{\infty}_t L^6_x} + \|U^0\|_{L^{\infty}_t L^6_x}^3 + \|U^0\|_{L^{\infty}_t L^6_x}^2
\end{multline}
lead by the modulation equations (\ref{mod}) to the straightforward inequality (where we used $\nu^0 \leq \delta$)
\be\begin{aligned}
\|\dot \pi\|_{\infty} + \|\dot \pi\|_1 \leq &C \sum_F (\|\langle U, \dot \xi_F\rangle - i\langle U, E\xi_F\rangle - i\langle N(U^0, \pi^0), \xi_F\rangle\|_{\infty} +\\
&+ \|\langle U, \dot \xi_F\rangle - i\langle U, E\xi_F\rangle - i\langle N(U^0, \pi^0), \xi_F\rangle\|_1) \\
\leq& C\delta \|U\|_{L^{\infty}_t L^6_x} + C\delta^2 \leq C\delta (\|U(0)\|_{W^{1, 2}} + \delta^2) + C\delta^2 \leq C\delta^2,
\end{aligned}\ee
which, for small $\delta$, proves some of the desired bounds on $\pi$.

Next we obtain decay estimates for $U$ and $\pi$, which are necessary in order to bound the quantities $\rho_{\infty}(t)$ and $y(t)-y_{\infty}$. The important ingredient in the proof is the evaluation of $\|\langle t \rangle^{\beta} U\|_{L^2_t L^{6+\infty}_x}$ by means of Lemma \ref{biest}. Here it becomes necessary to assume that the initial data $U(0)$ is in $L^q$, $q<4/3$, in addition to being in $H^{1/2}$. We make the choice of $\beta = 2/q-1-\epsilon>1/2$.

The reason why we need decay is
\begin{multline}
\|\langle t \rangle\dot\pi(t)\|_{1} \leq \int_0^{\infty} \langle t \rangle |\langle U, \dot \xi_F\rangle - i\langle U, E\xi_F\rangle - i\langle N(U^0, \pi^0), \xi_F\rangle| \dd t \\
\begin{aligned}
\leq &C\int_0^{\infty} \langle t \rangle |\pi^0(t)-\pi^0(\infty)| \|U(t)\|_{1+\infty} + \langle t \rangle \|(U^0)^2(t)\|_{1+\infty} + \\
&+ \langle t \rangle\|(U^0)^2(t)\|_{3/2+\infty} \|U^0(t)\|_3 \dd t\\
\leq &C(\|\langle t \rangle^{1/2}(\pi^0(t)-\pi^0(\infty))\|_{L^2_t} \|\langle t \rangle^{1/2} U(t)\|_{L^2_t L^{1+\infty}_x} + \\
&+\|\langle t \rangle^{1/2} U^0(t)\|_{L^2_t L^{1+\infty}_x}^2(1+\delta))\\
\leq &\frac 1 2 \nu^0 + C \|U(t)\|_{\langle t \rangle^{-1/2} L^2_t L^{1+\infty}_x}^2 + C \delta^2.
\end{aligned}
\end{multline}
This inequality is the means to prove that $\nu(t)$ stays bounded.

In order to apply Lemma \ref{biest} and bound this quantity, we evaluate the right-hand side terms of (\ref{nlsu}), beginning with the worst:
\be
\|VU\|_{\langle t \rangle^{-\beta}L^2_t L^{1\cap 6/5}_x} \leq C \delta  \|U\|_{\langle t \rangle^{-\beta}L^2_t L^{6+\infty}_x}.
\ee
Here we used the fact that $\nu^0(t)$ is bounded. Exactly the same works for the other linear term in $U$,
\be
\|L_{\pi^0}U\|_{\langle t \rangle^{-\beta}L^2_t L^{1\cap 6/5}_x} \leq C \delta \|U\|_{\langle t \rangle^{-\beta}L^2_t L^{6+\infty}_x}.
\ee

The remaining terms, $N(U^0, \pi^0)$ and $N_{\pi^0}(U^0, \pi^0)$, work out in the same fashion, provided that $U$ is uniformly bounded in time, which it is. Special attention has to be paid to the nonlocalized term, for which note that
\be\begin{aligned}
\|(U^0)^3\|_{\langle t \rangle^{-\beta} L^2_t L^{1\cap 6/5}_x} \leq& C\|U^0\|_{\langle t \rangle^{-\beta} L^2_t L^{6+\infty}_x} \|U^0\|_{L^{\infty}_t L^{2\cap3}_x}^{2} \\
\leq& C \|U^0\|_{\langle t \rangle^{-\beta} L^2_t L^{6+\infty}_x} \delta^2.
\end{aligned}\ee
This is another place where the $\dot H^{1/2}$-critical nature of the problem comes into play, since $\dot H^{1/2}$ embeds in $L^3$.

In conclusion, the nonlinear right-hand side terms have the same behaviour as the linear ones:
\be
\|N(U^0, \pi^0) + N_{\pi^0}(U^0, \pi^0)\|_{\langle t \rangle^{-\beta} L^2_t L^{1\cap 6/5}_x} \leq C  \delta \|U^0\|_{\langle t \rangle^{-\beta}L^2_t L^{6+\infty}_x}.
\ee

After applying Lemma \ref{biest}, the result is
\be
\int_0^{\infty} \langle t \rangle^{4/q-2-2\epsilon} \|P_c U(t)\|_{6+\infty}^2 \dd t \leq C(\|U(0)\|_{H^{1/2}\cap L^q}^2 + \delta^2\|U\|_{\langle t \rangle^{1-2/q+\epsilon}L^2_t L^{6+\infty}_x}^2 + \delta^4).
\ee
The hyperbolic part $P_{im} U$ decays just as fast by (\ref{eq_imw}) and the projection on the generalized $0$ eigenspace $P_{root} U$ is dominated by the other two components, so we can add them both in for free. By making $\delta$ and the initial data small we obtain the desired estimates
\be
\|U\|_{\langle t \rangle^{1-2/q+\epsilon} L^2_t L^{6+\infty}_x} \leq \delta \text{ and } \|\langle t\rangle \dot \pi\|_1 \leq C\delta^2 \leq \delta.
\ee

This finishes proving that the mapping $\Phi$ takes $X_{\delta}$ to itself, provided that $\delta$ and the initial data are sufficiently small. Next, we need to show that $\Phi$ really is a contraction within this set.

\subsection{Contraction} Here we prove that $\Phi$ is a contraction on $X_{\delta}$ for small $\delta <<1$.

Consider two solutions $(Z_1, \pi_1)$ and $(Z_2, \pi_2)$ of the linearized equation (\ref{ec_liniara}) corresponding to two different pairs of auxiliary functions, $(Z_1^0, \pi_1^0)$ and respectively $(Z_2^0, \pi_2^0)$. In our previous notation, we have that $(Z_1, \pi_1) = \Phi(Z_1^0, \pi_1^0)$ and $(Z_2, \pi_2) = \Phi(Z_2^0, \pi_2^0)$.

Assume that $(Z_1^0, \pi_1^0)$, $(Z_2^0, \pi_2^0) \in X_{\delta}$. Then, it follows from the first part of the proof that the same holds for $(Z_1, \pi_1)$ and $(Z_2, \pi_2)$.

We perform the contraction in the following seminorm: 
\be
Y = \{(Z, \pi) \mid \|(Z, \pi)\|_Y = \|\dot \pi\|_{L^{1 \cap \infty} \cap \langle t \rangle^{-\beta} L^2_t} + \|Z\|_{t^{1-\beta} L^{2}_t L^{6+\infty}_x} < \infty\}.
\ee
Here $\beta=2/q-1-\epsilon$ is the same as before and $-\beta<-1/2$, $1-\beta<1/2$.

It is straightforward to note that this seminorm defines a metric on $X_{\delta}$ because the elements of $X_{\delta}$ have well-determined initial parameter values. Furthermore, since the seminorm is weaker than the one that defines $X_{\delta}$, it makes $X_{\delta}$ a complete metric space.

We prove that the map $\Phi$ is a contraction on $X_{\delta}$ in this metric, more precisely that $\ds\|(Z_1, \pi_1) - (Z_2, \pi_2)\|_Y \leq \frac 1 2 \|(Z^0_1, \pi^0_1) - (Z^0_2, \pi^0_2)\|_Y$, for some sufficiently small~$\delta$.

Since this result will be reused later in the proof, it is convenient to state it in the form of a lemma:
\begin{lemma}
Consider $(Z_1, \pi_1) = \Phi(Z_1^0, \pi_1^0)$, $(Z_2, \pi_2) = \Phi(Z_2^0, \pi_2^0)$ solving the linearized equation (\ref{Z}), such that $(Z_1, \pi_1)$, $(Z_1^0, \pi_1^0)$, $(Z_2, \pi_2)$,  $(Z_2^0, \pi_2^0) \in X_{\delta}$, $\pi_1^0(0) = \pi_2^0(0) = (0, 0, 0, \alpha_0)$, and with initial data
\be
Z_i(0) = \bpm R_0^i\\ \overline{R_0^i}\epm + h_i {f^+(0)}
\ee
such that $R_0^i \perp \tilde f^+(0)$, $R_0^i \perp \xi_F(0)$. Further assume the orthogonality conditions $\langle Z_i(t), \xi_F^{Z_i}(t) \rangle = 0$ for all $t \geq 0$ and $F \in \{v_i, D_i, \alpha, \Gamma\}$, with $\xi^{Z_i}$ pertaining to $\pi_i^0$. Then, if $\delta$ is sufficiently small,
\be
\|(Z_1, \pi_1) - (Z_2, \pi_2)\|_Y \leq C(\delta \|(Z^0_1, \pi^0_1) - (Z^0_2, \pi^0_2)\|_Y + \|R_0^1-R_0^2\|_2).
\ee
\lb{compar}
\end{lemma}

\begin{proof}
Subtract the two copies of equation (\ref{Z}) corresponding to $Z_1$ and to $Z_2$ from one another and introduce the new function $U = \frak g_{\pi_1^0} (Z_1-Z_2)$.  Then
\begin{multline}
i U_t + \mc H_{\pi_1^0} U = V_{\pi^0_1} U + \\
+ \frak g_{\pi_1^0} (-i\dot\pi_1 \partial_{\pi}  W(\pi_1^0) +i \dot\pi_2 \partial_{\pi} W(\pi_2^0) + (\mc H^Z_{\pi^1_0} - \mc H^Z_{\pi^2_0}) Z_1 + N^Z(Z^0_1, \pi^0_1) - N^Z(Z^0_2, \pi^0_2)),
\label{U}
\end{multline}
with initial data
\be
U(0) = \frak g_{\pi^0}(0) \Big(\bpm R_0^1\\ \overline{R_0^1}\epm - \bpm R_0^2\\ \overline{R_0^2}\epm + (h_1-h_2) {f^+(0)} \Big).
\ee
It is worth noting that both paths start at the same point, so the operators $\mc H^Z_{\pi^0_i}$ are the same and have the same eigenfunctions for $i=1, 2$.

Henceforth, we shall use the names $\xi_F$, etc., in relation to $\pi_1^0$.

The vectors ${f^+}(0)$ and $\xi_F(0)$ do not depend on the whole parameter paths $\pi_i^0$, but only on the starting point $W(0)$. We do not assume that $R_0^1 = R_0^2$, in order to prove the Lipschitz continuity of $\Phi$ at the same time.

The modulation equations for $Z$ are obtained from the orthogonality condition for $Z$ (\ref{ort_Z}) and have the form, similar to (\ref{modula}),
\be
-i\dot\pi \partial_{\pi} W(\pi^0) = L^Z_{\pi^0}Z + N^Z_{\pi^0}(Z^0, \pi^0),
\ee
where
\be
L^Z_{\pi^0} Z = \sum_{F \in \{\alpha, \Gamma, v_i, D_i\}} \|\phi\|_2^{-2} (\langle Z, \dot \xi_F^Z\rangle - i\langle Z, E^Z\xi^Z_F\rangle) \eta^Z_F,
\ee
\be
N^Z_{\pi^0}(Z^0, \pi^0) = - \sum_{F \in \{\alpha, \Gamma, v_i, D_i\}} \|\phi\|_2^{-2} i\langle N^Z(Z^0, \pi^0), \xi^Z_F\rangle \eta^Z_F,
\ee
and
\be
E^Z = \bpm -|v(t)|^2 - 2iv(t)\dl & 0 \\ 0 & |v(t)|^2 - 2iv(t)\dl\epm.
\ee

Split equation (\ref{U}) into three parts, corresponding to our decomposition of the spectrum of $\mc H_{\pi_1^0}$ into the absolutely continuous part, the generalized eigenspace at $0$, and the imaginary eigenvalues, respectively. Then, we estimate each component separately.

The ranges of projections on the generalized eigenspace at $0$ and on the imaginary eigenspace are spanned by finitely many Schwartz functions,
\be
P_{root} U(t) = \sum_F a_F(t) \eta_F(\infty),\ P_{im} U(t) = b^+(t) f^+ + b^-(t) f^-.
\ee

Firstly, we deal with $P_{root} U$. Both $Z_1$ and $Z_2$ satisfy orthogonality conditions of the form $\langle Z_i(t), \xi^{Z_i}_F(t) \rangle = 0$. Applying the transformation $\frak g$ and taking the difference, one has that
\be
\langle U, \xi_F \rangle = \langle Z_2, \xi^{Z_2}_F - \xi^{Z_1}_F \rangle
\ee
and therefore
\be
\langle Z_2(t), \xi^{Z_2}_F(t) - \xi^{Z_1}_F(t) \rangle = \sum_F a_F(t) \langle \eta_{F}(\infty), \xi_{G}(t)\rangle + \langle (P_{im} + P_c) U(t), \xi_G(t)\rangle.
\ee
The matrix with entries $\langle \eta_{F}(\infty), \xi_{G}(t)\rangle$ is invertible with bounded norm. Hence, the following holds:
\begin{multline}
\|P_{root} U(t)\|_{1\cap\infty} \leq C \Big[\|(P_c+P_{im}) U(t)\|_{1+\infty} + \|Z_2(t)\|_{1+\infty} \Big(|\pi_1^0(t)-\pi_2^0(t)| + \\
+ \int_0^t |\pi^0_1(s) - \pi^0_2(s)| \dd s\Big)\Big].
\end{multline}
An immediate consequence is that
\be
\|P_{root} U(0)\|_{1 \cap \infty} \leq C\Big[\Big\|\bpm R_0^1\\ \overline{R_0^1}\epm - \bpm R_0^2\\ \overline{R_0^2}\epm\Big\|_{1 +\infty} + |h_1-h_2| + \delta \|\pi_1-\pi_2\|_{\infty}\Big].
\lb{initial_1}
\ee
The last term is not strictly necessary here, since the two parameter paths start at the same point. We still include it, in order to keep the argument general.

Note that
\begin{multline}
\|(\|Z_2(t)\|_{1+\infty} |\int_0^t \pi_1^0(s) - \pi_2^0(s) \dd s|)\|_{\langle t \rangle^{1-\beta}L^2_t} \leq\\
\leq C\|Z_2\|_{\langle t \rangle^{1-2/q} L^2_t L^{2+\infty}_x} \|\pi^0_1-\pi^0_2\|_{\infty} \leq C\delta \|\dot\pi^0_1-\dot\pi^0_2\|_1.
\end{multline}
Therefore
\be
\|P_{root} U\|_{\langle t \rangle^{1-\beta} L^2_t L^{6+\infty}_x} \leq C(\|(P_c+P_{im}) U(t)\|_{\langle t \rangle^{1-\beta} L^2_t L^{6+\infty}_x} + \delta \|\dot\pi^0_1-\dot\pi^0_2\|_{1 \cap \infty}).
\lb{ProotZ}
\ee

Now we evaluate the terms on the right-hand side of equation (\ref{U}), in order to bound the two remaining components, $P_{im} U$ and $P_c U$. Since $P_{root} U$ grows like $t^{2-2/q-\epsilon}<1/2$ in $L^2$ in time, we prove that the other two components have the same growth rate in $L^2$ in time.

We evaluate the projection on the countinuous spectrum. Just as in the stability part of the proof, Lemma \ref{biest} leads to the required weighted estimate. From equation (\ref{U}) we have that
\be
\|P_c U\|_{\langle t \rangle^{1-\beta}L^{2}_t L^{6+\infty}_x} \leq C(\|U(0)\|_{2} + \|RHS\|_{\langle t \rangle^{1-\beta}L^{2}_t L^{1 \cap 6/5}_x}).
\ee

Thus, if we can establish the bound for the right-hand side of the equation, we then retrieve it for $P_c U$.

Observe that
\begin{align}
\ds\|V_{\pi^0_1} U\|_{\langle t \rangle^{1-\beta}L^{2}_t L^{1\cap 6/5}_x} \leq& C\delta \|U\|_{\langle t \rangle^{1-\beta}L^{2}_t L^{6+\infty}_x};\nonumber\\
\ds\|(L^Z_{\pi^0_1} - L^Z_{\pi^0_2}) Z_2\|_{\langle t \rangle^{1-\beta}L^{2}_t L^{1 \cap 6/5}_x} \leq& C \|d\pi^0\|_Y \|Z_2\|_{\langle t \rangle^{-\beta}L^{2}_t L^{6+\infty}_x} \nonumber\\
\leq& C \|d\pi^0\|_Y \delta;\\
\ds\|\frak g_{\pi^0_1} L^Z_{\pi^0_1} U\|_{\langle t \rangle^{1-\beta}L^{2}_t L^{1\cap 6/5}_x} \leq& \delta \|U\|_{\langle t \rangle^{1-\beta}L^{2}_t L^{6+\infty}_x}; \nonumber\\
\ds\|(\mc H^Z_{\pi^1_0} - \mc H^Z_{\pi^2_0}) Z_1\|_{\langle t \rangle^{1-\beta}L^{2}_t L^{1\cap 6/5}_x} \leq& C\|d\pi^0-d\pi^0(\infty)\|_{\infty} \|Z_1\|_{\langle t \rangle^{-\beta}L^{2}_t L^{6+\infty}_x} \nonumber\\
\leq& C\|d\dot\pi^0\|_{1\cap\infty} \delta,\nonumber
\end{align}
as well as
\be\begin{aligned}
&\|N^Z(Z^0_1, \pi^0_1) - N^Z(Z^0_1, \pi^0_2)\|_{\langle t \rangle^{1-\beta}L^{2}_t L^{1\cap 6/5}_x} \leq\\
&\qquad\begin{aligned}
&\leq C\|\pi^0_1-\pi^0_2\|_{\infty} \|Z^0_1\|_{\langle t \rangle^{-\beta}L^{2}_t L^{1\cap 6/5}_x} \|Z^0_1\|_{L^{\infty}_t L^2_x}\\
&\leq C\|\dot\pi^0_1 - \dot \pi^0_2\|_{1\cap\infty} \delta^2;
\end{aligned}
\\
&\|N^Z(Z^0_1, \pi^0_2) - N^Z(Z^0_2, \pi^0_2)\|_{\langle t \rangle^{1-\beta}L^{2}_t L^{1\cap 6/5}_x} \leq \\
&\qquad\begin{aligned}
\leq &C\|Z^0_1-Z^0_2\|_{\langle t \rangle^{1-\beta}L^{2}_t L^{6+\infty}_x} (\delta^2+\delta);
\end{aligned}
\\
&\|N^Z_{\pi^0_1}(Z^0_1, \pi^0_1) - N^Z_{\pi^0_2}(Z^0_2, \pi^0_2)\|_{\langle t \rangle^{1-\beta}L^{2}_t L^{1\cap 6/5}_x} \leq \\
&\qquad\begin{aligned}
\leq &\|N^Z_{\pi^0_1}(Z^0_1, \pi^0_1) - N^Z_{\pi^0_2}(Z^0_1, \pi^0_1)\|_{\langle t \rangle^{1-\beta}L^{2}_t L^{1\cap 6/5}_x} + \\
&+ \|N^Z_{\pi^0_2}(Z^0_1, \pi^0_1) - N^Z_{\pi^0_2}(Z^0_1, \pi^0_2)\|_{\langle t \rangle^{1-\beta}L^{2}_t L^{1\cap 6/5}_x} + \\
&+ \|N^Z_{\pi^0_2}(Z^0_1, \pi^0_2) - N^Z_{\pi^0_2}(Z^0_2, \pi^0_2)\|_{\langle t \rangle^{1-\beta}L^{2}_t L^{1\cap 6/5}_x} \\
\leq &C\delta(\|\dot\pi^0_1 - \dot \pi^0_2\|_1 + \|Z^0_1-Z^0_2\|_{\langle t \rangle^{1-\beta}L^{2}_t L^{6+\infty}_x}).
\end{aligned}
\end{aligned}\ee
We again used the $\dot H^{1/2}$-critical nature of the equation and the Keel-Tao endpoint Strichartz estimate, in the next-to-last inequality, concerning unlocalized terms, under the following guise:
\be
\|(Z^0_i)^2(t)\|_{L^{3/2\cap\infty}_x} \leq C\|Z^0_i(t)\|^2_{H^{1/2}_x} \leq C\delta^2.
\ee

After examining each term, the overall conclusion is that
\be\begin{aligned}
\|RHS\|_{\langle t \rangle^{1-\beta}L^{2}_t L^{1\cap 6/5}_x} \leq C\delta (&\|U\|_{\langle t \rangle^{1-\beta}L^{2}_t L^{16+\infty}_x} + \|\dot\pi_1^0-\dot\pi_2^0\|_{1\cap\infty} + \\
&+ \|Z^0_1-Z^0_2\|_{\langle t \rangle^{1-\beta}L^{2}_t L^{6+\infty}_x})
\end{aligned}\ee
and thus
\begin{multline}
\|P_c U\|_{\langle t \rangle^{1-\beta}L^{2}_t L^{6+\infty}_x} \leq\\
\leq C(\|U(0)\|_{2} + \delta \|U\|_{t^{1/10} L^{\infty}_t L^{3+\infty}_x} + \delta \|(dU^0, d\pi^0)\|_Y).
\lb{PcZ}
\end{multline}

Next, we bound $P_{im} U$. Note that $\|P_{im} U(t)\|_2 \leq \|Z_1\|_{L^{\infty}_t L^2_x} + \|Z_2\|_{L^{\infty}_t L^2_x} \leq C\delta$, so $P_{im} U$ is a bounded solution to a hyperbolic ODE system. Therefore, by yet another application of Lemma \ref{hyp},
\be
0 = \langle U(0), \tilde f^+ \rangle + \int_0^{\infty} e^{-t\sigma} N_+(t) \dd t.
\ee
where $\pm i\sigma$ are the imaginary eigenvalues of $\mc H_{\pi^0_1}$.
Thus
\be\begin{aligned}
\langle U(0), f^+_{\infty} \rangle \leq &C\|N_+\|_{\langle t \rangle^{1-\beta}L^{2}_t L^{1\cap 6/5}_x} \leq C\|RHS\|_{\langle t \rangle^{1-\beta}L^{2}_t L^{1\cap 6/5}_x}\\
\leq &C\delta (\|U\|_{\langle t \rangle^{1-\beta}L^{2}_t L^{6+\infty}_x} + \|(Z^0_1 - Z^0_2, \pi^0_1-\pi^0_2)\|_Y)
\lb{initial_2}
\end{aligned}\ee
and, by the same reasoning as in the stability proof (see \ref{eq_imw}),
\begin{multline}
\|P_{im} U\|_{\langle t \rangle^{1-\beta}L^{2}_t L^{6+\infty}_x} \leq \\
\leq C(\|U(0)\|_{1+\infty} + \delta \|U\|_{\langle t \rangle^{1-\beta}L^{2}_t L^{6+\infty}_x} + \delta\|(U^0_1-U^0_2, \pi^0_1-\pi^0_2)\|_Y).
\lb{PimZ}
\end{multline}
Note $\sigma$ belongs to a compact subset of $(0, \infty)$, because $\alpha$ belongs to a compact subset of $(0, \infty)$, so the constants are independent of $\sigma$.

Now we deal with the initial data $U(0)$:
\be
\|U(0)\|_p \leq C\Big(\Big\|\bpm R_0^1\\ \overline{R_0^1}\epm - \bpm R_0^2\\ \overline{R_0^2}\epm\Big\|_p + |h_1 - h_2|\Big).
\ee
Furthermore, subtracting the two copies of (\ref{bplus}) from one another, one has
\begin{multline}
b_+^1(0) - b_+^2(0) = \langle g_{\pi^0}(0) (Z_1(0) - Z_2(0)), \tilde f^+\rangle \\
\begin{aligned}
&= \Big\langle \frak g_{\pi^0} \Big(\bpm R_0^1 \\ \ov{R_0^1} \epm -  \bpm R_0^2 \\ \ov{R_0^2} \epm + (h_1-h_2) f^+(0)\Big), \tilde f^+\Big\rangle\\
&= \Big\langle \bpm R_0^1 \\ \ov{R_0^1} \epm -  \bpm R_0^2 \\ \ov{R_0^2} \epm, \tilde f^+ - \frak g_{\pi^0}^{-1} \tilde f^+(0) \Big\rangle + (h_1-h_2)(1 + \langle (g_{\pi^0} f^+(0)) - f^+, \tilde f^+ \rangle)
\end{aligned}\end{multline}
and thus, for small $\delta$, by  (\ref{initial_2})
\be\begin{aligned}
|h_1 - h_2| &\leq C \Big(|b_+^1(0) - b_+^2(0)| + \delta \Big\|\bpm R_0^1 \\ \ov{R_0^1} \epm -  \bpm R_0^2 \\ \ov{R_0^2} \epm\Big\|_{1\cap \infty} \Big)\\
&\leq C \delta (\|U\|_{\langle t \rangle^{1-\beta}L^{2}_t L^{6+\infty}_x} + \|(Z^0_1 - Z^0_2, \pi^0_1-\pi^0_2)\|_Y + \|R_0^1-R_0^2\|_{1\cap\infty}).
\lb{h1h2}
\end{aligned}\ee
By (\ref{initial_1}) it follows that, for small $\delta$,
\be
\|U(0)\|_p \leq C\big(\|R_0^1 - R_0^2\|_p + \delta \|U\|_{\langle t \rangle^{1-\beta}L^{2}_t L^{6+\infty}_x} + \delta \|(Z^0_1-Z^0_2, \pi^0_1-\pi^0_2)\|_Y\big).
\ee

Putting (\ref{ProotZ}), (\ref{PimZ}), and (\ref{PcZ}) together, one has that
\begin{multline}
\|U\|_{\langle t \rangle^{1-\beta}L^{2}_t L^{6+\infty}_x} \leq \\
\leq C\big(\delta \|U\|_{\langle t \rangle^{1-\beta}L^{2}_t L^{6+\infty}_x} + \delta \|(Z^0_1, \pi^0_1) - (Z^0_2, \pi^0_2)\|_Y + \|R_0^1 - R_0^2\|_2\big),
\end{multline}
whence, for small $\delta$,
\be
\|U\|_{\langle t \rangle^{1-\beta}L^{2}_t L^{6+\infty}_x} \leq C\delta \|(Z^0_1-Z^0_2, \pi^0_1-\pi^0_2)\|_Y + C\|R_0^1 - R_0^2\|_{2}.
\ee

In studying the difference $\pi_1-\pi_2$, we switch back to using $U_1$ and $U_2$ instead of $Z_1$ and $Z_2$. The reason why this is possible is that
\begin{multline}
\|U_1 - U_2\|_{\langle t \rangle^{1-\beta} L^{2}_t \langle x \rangle L^{6+\infty}_x} = \|\frak g_{\pi^0_1} Z_1 - \frak g_{\pi^0_2} Z_2\|_{\langle t \rangle^{1-\beta} L^{2}_t \langle x \rangle L^{6+\infty}_x} \\
\begin{aligned}
&\leq \|\frak g_{\pi^0_1} Z_1 - \frak g_{\pi^0_2} Z_1\|_{\langle t \rangle^{1-\beta} L^{2}_t \langle x \rangle L^{6+\infty}_x} + \|\frak g_{\pi^0_2} Z_1 - \frak g_{\pi^0_2} Z_2\|_{\langle t \rangle^{1-\beta} L^{2}_t \langle x \rangle L^{6+\infty}_x}\\
&\leq C\|\pi^0_1-\pi^0_2\|_{\infty} \|\langle t \rangle Z_1\|_{\langle t \rangle^{1-\beta} L^{2}_t L^{6+\infty}_x} + \|Z_1 - Z_2\|_{t^{1-\beta} L^{2}_t L^{6+\infty}_x}\\
&\leq \max(C, C\delta) \|(Z_1, \pi_1) - (Z_2, \pi_2)\|_Y.
\end{aligned}
\end{multline}
To put it otherwise, this norm is sufficiently weak not to see such small symmetry transformations. This helps us insofar as all the terms that appear in the modulation equations are localized.

The modulation equations for $\pi_1-\pi_2$ are of the form, derived from (\ref{mod}),
\be
-i\dot\pi_1 \bpm W(\pi_1^0) \\ \ov{W(\pi_1^0)} \epm + i\dot\pi_2 \bpm W(\pi_2^0) \\ \ov{W(\pi_2^0)} \epm = L_{\pi_1^0}U_1 - L_{\pi_2^0}U_2 + N_{\pi_1^0}(U^1_0, \pi^1_0) - N_{\pi_2^0}(U^2_0, \pi^2_0).
\ee
Denoting for convenience $d\pi = \pi_1-\pi_2$, one has that
\be\begin{aligned}
\|d\dot \pi\|_{L^1_t} \leq \|d\dot \pi\|_{\langle t \rangle^{-\beta} L^2_t} \leq &\|L_{\pi^0_1} U_1 - L_{\pi^0_2} U_2\|_{\langle t \rangle^{-\beta} L^2_t L^{1+\infty}_x} + \\
& + \|N_{\pi^0_1}(U^0_1, \pi^0_1) - N_{\pi^0_2}(U^0_2, \pi^0_2)\|_{\langle t \rangle^{-\beta} L^2_t L^{1+\infty}_x}\\
&+ C\Big\|\dot \pi_2 \bpm W(\pi_1^0) \\ \ov{W(\pi_1^0)} \epm - \bpm W(\pi_2^0) \\ \ov{W(\pi_2^0)} \epm\Big\|_{\langle t \rangle^{-\beta} L^2_t L^{1+\infty}_x}.
\end{aligned}\ee
An important fact is that
\begin{multline}
\|L_{\pi^0_1}(T) - L_{\pi^0_2}(T)\|_{L^{1+\infty}_x \to L^{1\cap \infty}_x} \leq \\
\begin{aligned}
\leq &\sum_{F \in \{v_i, D_i\}} \big\|\|\phi^1\|_2^{-2} (\dot \xi_F^{1} - iE\xi^{1}_F) \otimes \eta^{1}_F - \|\phi^2\|_2^{-2} (\dot \xi_F^{2} - iE^{Z}\xi^{2}_F) \otimes \eta^{2}_F\big\|_{L^{1+\infty}_x \to L^{1\cap \infty}_x} +\\
&\sum_{F \in \{\alpha, \Gamma\}} \big\| 2\alpha_1\|\phi^1\|_2^{-2} (\dot \xi_F^{1} - iE\xi^{1}_F) \otimes \eta^{1}_F - 2\alpha_2\|\phi^2\|_2^{-2} (\dot \xi_F^{2} - iE^{Z}\xi^{2}_F) \otimes \eta^{2}_F\big\|_{L^{1+\infty}_x \to L^{1\cap \infty}_x}\\
\leq &C\Big((|\pi^0_1(T)-\pi^0_1(\infty)| + |\pi^0_2(T)-\pi^0_2(\infty)|)\Big(\int_0^T |d\pi^0(t)| \dd t + |d\pi^0(T)|\Big) + \\
&+ |d\pi^0(T) - d\pi^0(\infty)| + |d\dot\pi^0(T)-d\dot\pi^0(\infty)|\Big)\\
\leq &C \delta \|d\pi^0\|_{\infty}.
\end{aligned}
\end{multline}
Indeed, despite the fact that the parameter paths converge to different final values, the difference of these two operators decays in time.

Using the estimates
\begin{multline}
\|(L_{\pi^0_1} - L_{\pi^0_2}) U_1\|_{\langle t \rangle^{-\beta} L^2_t L^{1\cap \infty}_x} \leq \\
\begin{aligned}
\leq &C\ds\| \Big((|\pi^0_1(t)-\pi^0_1(\infty)| + |\pi^0_1(t)-\pi^0_1(\infty)|) \big(\int_0^t |d\pi^0(s)| \dd s + |d\pi^0(t)|\big) + \\
&+ |d\pi^0(t) - d\pi^0(\infty)| + |d\dot \pi^0(t) - d\dot \pi^0(\infty)| \Big) U_1(t)\|_{\langle t \rangle^{-\beta} L^2_t L^{1+ \infty}_x} \\
\ds\leq &C\|\dot\pi^0_1-\dot \pi^0_2\|_{1\cap \infty} \big(\|t (\pi^0_1(t) - \pi^0_1(\infty))\|_{L^{1\cap\infty}_t} + 1\big)\|U_1\|_{\langle t \rangle^{-\beta} L^2_t L^{2+\infty}_x}\\
\leq &C\delta \|d\dot \pi^0\|_1,
\end{aligned}
\end{multline}
\begin{multline}
\|L_{\pi^0_2} (U_1-U_2)\|_{\langle t \rangle^{-\beta} L^2_t L^{1\cap\infty}_x} \leq \\
\begin{aligned}
\ds \leq &C\Big\|(|\pi^2_0(t) - \pi^2_0(\infty)| + |\dot \pi^2_0(t)|) \|U_1(t)-U_2\|_{L^{1+\infty}_x}\Big\|_{\langle t \rangle^{-\beta} L^2_t} \\
\ds \leq &C\|U_1-U_2\|_{\langle t \rangle^{1-\beta} L^{2}_t \langle x \rangle L^{6+\infty}_x} \|t(\pi^2_0(t) - \pi^2_0(\infty))\|_{L^{1\cap\infty}_t}\\
\leq &C\delta \|U_1-U_2\|_{\langle t \rangle^{1-\beta} L^{2}_t \langle x \rangle L^{6+\infty}_x},
\end{aligned}\end{multline}
\begin{multline}
\|N_{\pi^0_1}(U^0_1, \pi^0_1) - N_{\pi^0_2}(U^0_2, \pi^0_2)\|_{\langle t \rangle^{-\beta} L^{2}_t L^{1\cap\infty}_x} \leq \\
\begin{aligned}
\leq &\|N_{\pi^0_1}(U^0_1, \pi^0_1) - N_{\pi^0_2}(U^0_1, \pi^0_1)\|_{\langle t \rangle^{-\beta} L^{2}_t L^{1\cap\infty}_x} + \\
&+ \|N_{\pi^0_2}(U^0_1, \pi^0_1) - N_{\pi^0_2}(U^0_1, \pi^0_2)\|_{\langle t \rangle^{-\beta} L^{2}_t L^{1\cap\infty}_x} + \\
&+ \|N_{\pi^0_2}(U^0_1, \pi^0_2) - N_{\pi^0_2}(U^0_2, \pi^0_2)\|_{\langle t \rangle^{-\beta} L^{2}_t L^{1\cap\infty}_x} \\
\leq &C\delta(\|d\dot\pi\|_1 + \|U^0_1-U^0_2\|_{\langle t \rangle^{1-\beta} L^{2}_t \langle x \rangle L^{6+\infty}_x}),
\end{aligned}
\end{multline}
and
\be
\Big\|\dot \pi_2 \bpm W(\pi_1^0) \\ \ov{W(\pi_1^0)} \epm - \bpm W(\pi_2^0) \\ \ov{W(\pi_2^0)} \epm\Big\|_{\langle t \rangle^{-\beta} L^2_t L^{1+\infty}_x} \leq C \delta \|d\pi^0\|_{\infty} \leq C \delta \|d\dot\pi^0\|_1,
\ee
one gets the desired conclusion that
\begin{multline}
\|d\dot \pi\|_{1\cap\infty} \leq C\delta (\|(Z^0_1-Z^0_2, \pi^0_1-\pi^0_2)\|_Y + \|U_1-U_2\|_{\langle t \rangle^{1-\beta} L^{2}_t L^{6+\infty}_x}).
\end{multline}
Therefore,
\be
\|(Z_1-Z_2, \pi_1-\pi_2)\|_Y \leq C\delta (\|(Z^0_1-Z^0_2, \pi^0_1-\pi^0_2)\|_Y + \|R_0^1 - R_0^2\|_{2})
\lb{lip}
\ee
for sufficiently small $\delta$.
\end{proof}

Since $\Phi$ is a contraction in a complete metric space, it has a fixed point $\Big(Z =\bpm R\\ \overline{R}\epm, \pi\Big)~\in~X_{\delta}$, such that $R + W(\pi)$ is a global solution for (\ref{NLS}).

\subsection{Remaining bounds} The correction term $\mc F(R_0) = h(R_0) {f^+}^Z $ satisfies the appropriate bounds since
\be
|h(R_0)| \leq C\delta^2
\ee
by (\ref{h}) and (\ref{est_U}), for $\|U(0)\|_{L^{4/3-\epsilon} \cap H^{1/2}} \leq C\delta$, and
\be
|h(R^1_0) - h(R^2_0)| \leq C \delta\|R^1_0-R^2_0\|_{2}
\ee
by (\ref{h1h2}) and (\ref{lip}). Same goes for the $a_F$ coefficients.

If we did not assume that $R_0^1 = R_0^2$ in (\ref{lip}), then it implies that $\Phi(R_0)$ is a Lipschitz function of $R_0$. Indeed
\be
\|(Z_1, \pi_1) - (Z_2, \pi_2)\|_Y \leq C\delta \|R_0^1 - R_0^2\|_{2}
\ee
implies
\be
\|R_1 - R_2\|_{\langle t \rangle^{1-\beta} L^2_t L^{6+\infty}_x} \leq C\delta \|R_0^1 - R_0^2\|_{2}
\ee
and
\be
\|W(\pi_1) - W(\pi_2)\|_{{\langle t \rangle^{1-\beta} L^2_t \langle x \rangle L^{6+\infty}_x}} \leq C\delta \|R_0^1 - R_0^2\|_{3/2 \cap 2}.
\ee

Only the scattering is left to prove and it follows in a standard manner using the Strichartz inequalities. Indeed, observe that
\be
U_1 = \int_0^{\infty} e^{-is\mc H} RHS(s) \dd s
\ee
is in $L^2$ and the integral converges in the $L^2$ norm, because $RHS$ is in $L^2_t L^6_x$. Therefore,
\be
P_c U(t) - e^{it\mc H} P_c (U(0) + U_1) = -e^{it \mc H} P_c \int_t^{\infty} e^{-is\mc H} P_c RHS(s) \dd s = o_{L^2}(1).
\ee
The other two components of $U$, $P_{im} U$ and $P_{root} U$ converge to zero in the $L^2$ norm. Indeed, they easily converge to zero in other norms and, being given by Schwartz functions in the space variable, all of their Lebesgue norms are comparable. Thus $U$ behaves like $e^{it\mc H} P_c (U(0) + U_1) + o_{L^2}(1)$.

Let $\mc H_0(\alpha_{\infty}) = (\Delta - \alpha_{\infty})\sigma_3$, where $\sigma_3 = \bpm 1 & 0 \\ 0 & -1 \epm$. We want to establish that
\be
L = \lim_{t \to \infty} e^{it\mc H_0(\alpha_{\infty})} e^{-it\mc H} P_c(U(0) + U_1)
\ee
exists as a strong $L^2$ limit. But, letting $M = \bpm 2\phi_{\infty}^2 & \phi_{\infty}^2 \\ -\phi_{\infty}^2 & -2\phi_{\infty}^2 \epm = \mc H - \mc H_0$,
\be
\frac d {\dd t} e^{it\mc H_0(\alpha_{\infty})} e^{-it\mc H} P_c (U(0) + U_1) = e^{it\mc H_0(\alpha_{\infty})} M e^{-it\mc H} (U(0) + U_1).
\ee
In other words,
\be
L = P_c (U(0) + U_1) + \lim_{t \to \infty} \int_0^t e^{it\mc H_0(\alpha_{\infty})} M e^{-it\mc H} P_c (U(0) + U_1) \dd t.
\ee
However, we note that
\be\begin{aligned}
\Big\|\int_{T_1}^{T_2} e^{it\mc H_0(\alpha_{\infty})} M e^{-it\mc H} P_c (U(0) + U_1) \dd t\Big\|_2 \leq& C\int_{T_1}^{T_2} \|M e^{-it\mc H} P_c (U(0) + U_1)\|_{6/5}^2 \dd t\\
\leq &C\int_{T_1}^{T_2} \|e^{-it\mc H} P_c (U(0) + U_1)\|_{6}^2 \dd t\\
\leq &C\|U(0) + U_1\|_2^2.
\end{aligned}\ee
Since this last integral is absolutely convergent, the initial one also converges. Therefore $L$ exists as a strong $L^2$ limit and
\be
U(t) = e^{it \mc H_0(\alpha_{\infty})} L + o_{L^2}(1).
\ee
Switching back to $Z(t) = \frak g_{\pi}^{-1}(t) U(t)$, one has
\be
Z(t) = \frak g_{\pi}^{-1}(t) (e^{it \mc H_0(\alpha_{\infty})} L + o_{L^2}(1)) = e^{it \mc H_0(0)} \frak g_{\pi}^{-1}(0) L + o_{L^2}(1),
\ee
which finishes the proof of scattering.
\end{proof}

\begin{proof}[Proof of Corollary \ref{cor17}]
Firstly, apply a symmetry transformation to the whole equation in order to make $W(0) = \phi(\cdot, \alpha_0)$ a positive ground state. The result that holds for the transformed problem is also valid for the original one.

The previous Theorem \ref{main} proves the existence of a codimension-nine Lipschitz submanifold of $L^q \cap H^{1/2}$, $1 \leq q < 4/3$, tangent to $S_{\delta}$ at $W(0)$, made of initial data $\Psi(R_0)(0)$, $R_0 \in S_{\delta}$, for global solutions of (\ref{NLS}). Since the equation is invariant under symmetry transformations, we let symmetry transformations act on this submanifold and retrieve a larger set of good initial data. More precisely, let
\be
\mc N_{L^q \cap H^{1/2}} = \bigcup_{\mc G} \{\mc G (\Psi(R_0)(0)) \mid R_0 \in S_{\delta}\}.
\ee
The action $\mc G (\Psi(R_0)(0))$ maps the product between an 8-dimensional manifold of symmetry transformations and the codimension-nine submanifold into $L^q \cap H^{1/2}$. The matrix
\be
\Big\langle \frac{\partial \frak g}{\partial G} \bpm W(0) \\ \ov{W(0)} \epm, \xi_F(W(0)) \Big\rangle = \langle \eta_G(W(0)), \xi_F(W(0)) \rangle,
\ee
where $F$ and $G$ stand for parameters $v_k$, $D_k$, $\alpha$, and $\Gamma$ is invertible by Lemma \ref{orto}, which implies that the action of the symmetry transformations is transverse to the codimension-nine manifold. Therefore, the range of the map is locally (in a neighborhood of $W(0)$) a codimension-one submanifold of $L^q \cap H^{1/2}$.
\end{proof}

\begin{proof}[Proof of Theorem \ref{tichie}]
The first new claim is that if a solution obtained by the previous corollary has initial data $\Psi(R_0)(0) = W(0) + R_0 + \mc F_1(R_0) \in H^1$, then $\Psi(R_0)(T)$ is still in $H^1$ at any time $T>0$.

Using the change of coordinates $U = \frak g_{\pi}(t) \Big(\bpm \Psi(R_0)(t) \\ \ov{\Psi(R_0)(t)} \epm - \bpm W(\pi(t)) \\ \ov{W(\pi(t))} \epm\Big)$, we transform the equation in the same manner as in the proof of Theorem \ref{main} to the form (\ref{nlsu}). Taking derivatives, one has that
\be
i\partial_t \dl U + \mc H_{\pi} \dl U = [\mc H_{\pi}, \dl] U  + \dl (V_{\pi} U) + \dl (L_{\pi} U) + \dl N(U, \pi) + \dl (N_{\pi}(\pi, U)). 
\ee
This equation is linear in $\dl U$. By means of Strichartz estimates, for $U \in X_{\delta}$ we obtain
\be\begin{aligned}
\|\dl U\|_{L^2_t L^6_x \cap L^{\infty}_t L^2_x} \leq& C(\|\dl U(0)\|_2 + \|RHS\|_{L^{2}_t L^{6/5}_x + L^1_t L^2_x})\\
\leq& C(\|\dl U(0)\|_2 + \delta\|\dl U\|_{L^2_t L^6_x \cap L^{\infty}_t L^2_x} + \delta) \leq C(\|\dl U(0)\|_2 + \delta)
\end{aligned}\ee
for sufficiently small $\delta$. By approximating $\dl U(0)$ with functions of better regularity, the bound is seen to hold everywhere in time, instead of almost everywhere. This proves the statement.

Secondly, we prove that if a $H^1$ solution has initial data $\Psi(R_0)(0) = W(0) + R_0 + \mc F_1(R_0) \in |x|^{-1} L^2$, then the property of being in $|x|^{-1} L^2$ is kept at any later time $T>0$.
Using the same machinery as in the previous proof, we can reduce this to showing that if $|x| U(0) \in L^2$ in the equation
\be
i U_t + \mc H_{\pi} U = -i\dot \pi \partial_{\pi}\frak g_{\pi} \bpm W(\pi) \\ \ov{W(\pi)} \epm + N(U, \pi) + V_{\pi}U,
\ee
then the property is preserved at time $T>0$. However, $xU$ satisfies an equation of its own, namely
\be
i (xU)_t + \mc H_{\pi^0} (xU) = 2\dl\sigma_3 U - ix\dot \pi \partial_{\pi}\frak g_{\pi} \bpm W(\pi) \\ \ov{W(\pi)} \epm + xN(U, \pi) + xV_{\pi}U.
\ee
The local terms on the right-hand side are already bounded in Strichartz norms by our knowledge about $U$. This leaves
\be
\|\dl\sigma_3 U\|_{L^1(0, T; L^2_x)} \leq C T \|\dl U\|_{L^{\infty}_t L^2_x}
\ee
and
\be
\|xU^3\|_{L^1_t L^{6/5}_x} \leq C \|xU\|_{L^{\infty}_t L^2_x} \|U\|_{L^2_t L^6_x}^2 \leq C \delta \|xU\|_{L^{\infty}_t L^2_x}.
\ee
As expected, the gradient term grows linearly in time and we cannot do any better. By a standard argument, for sufficiently small $\delta$ it follows that
\be
\|xU(T)\|_{2} \leq C(\|xU(0)\|_{2} + T \delta). 
\ee

Note that $\|\dl U(0)\|_2 \leq C \|\dl R_0\|_2$ and likewise $\|x U(0)\|_2 \leq C \|xR_0\|_2$. Thus the solution stays for all times in the $\Sigma$ space, but the norm may grow linearly in time.

The second claim was that the manifold of global solutions $\mc N$ is locally invariant under the flow. Define the manifold as
\be\begin{aligned}
\mc N  = \{\Psi(0) \mid & \Psi(0) \in \mc N_{L^q \cap H^{1/2}},\ \Psi(0) = \mc G(W(0) + R_0 + \mc F_1(R_0)),\\
&\text{for $R_0$ with }\|R_0 + \mc F_1(R_0)\|_{\Sigma} < \delta_0\}
\end{aligned}\lb{defN}\ee
with $\mc N_{L^q \cap H^{1/2}}$ being the codimension-one submanifold from Corollary \ref{cor17}. The only new condition pertains to the size of the initial data in $\Sigma$. Clearly $\mc N$ is still a codimension-one submanifold of $\Sigma$.

Note that every globally bounded solution sufficiently close to the manifold of standing waves must actually start on $\mc N_{L^q \cap H^{1/2}}$. We phrase this as the following lemma:
\begin{lemma}
Consider a solution $\Psi$ to equation (\ref{NLS}) such that $\bpm \Psi \\ \ov{\Psi} \epm = Z + \bpm W(\pi) \\ \ov{W(\pi)} \epm$, $Z(t) \perp \xi_F^Z(\pi(t))$ for all $t\geq 0$, and $(Z, \pi) \in X_{\delta}$ for some fixed, sufficiently small $\delta$. Further assume that $\|Z(0)\|_{L^q \cap H^{1/2}} \leq C \delta$ for the same $\delta$.

Then $\Psi(0) \in \mc N_{L^q \cap H^{1/2}}$.
\lb{comparatie}
\end{lemma}
\begin{proof}
One can perform small symmetry transformations in this situation, since $\mc N_{L^q \cap H^{1/2}}$ is invariant under them. Therefore assume, without loss of generality, that $W(\pi(0))$ is a positive ground state, meaning $\pi(0) = (0, 0, 0, \alpha_0)$. Let $\bpm R_0 \\ \ov{R_0} \epm = (P_c(0)+P_-(0)) Z(0)$; then $\|R_0\|_{L^q \cap H^{1/2}} \leq C\delta$.

Consider the global solution $\Psi(R_0)$, of initial data
\be
\Psi(R_0)(0) = W(\pi(0)) + R_0 + \mc F_1(R_0).
\ee
It exists by Theorem \ref{main} if $\delta$ is sufficiently small and it has the property that, for some path $\pi(R_0)$ with $\pi(R_0)(0) = (0, 0, 0, \alpha_0)$,
\be
\bpm \Psi(R_0) \\ \ov{\Psi(R_0)} \epm = Z(R_0) + \bpm W(\pi(R_0))\\ \ov{W(\pi(R_0))} \epm,\ (Z(R_0), \pi(R_0)) \in X_{\delta},
\ee
and $Z(R_0)(t) \perp \xi_F^Z(\pi(R_0)(t))$.

Thus the conditions of Lemma \ref{compar} are met and one can compare the global solutions $\Psi$ and $\Psi(R_0)$. The immediate result is that $\Psi = \Psi(R_0)$. However, $\Psi(R_0)$ belongs to $\mc N_{L^q \cap H^{1/2}}$ by construction, which finishes the proof.
\end{proof}

With the help of Lemma \ref{comparatie} it is straightforward to prove that $\mc N$ is locally in time invariant. Indeed, consider the truncated solution $\Psi_\tau$ obtained by restricting some global solution $\Psi$, $\Psi(0) \in \mc N$, to the time interval $[\tau, \infty)$. Keeping the same notations, one has for every $t$
\be
\bpm \Psi \\ \ov{\Psi} \epm = Z + \bpm W(\pi) \\ \ov{W(\pi)} \epm,\ Z(t) \perp \xi_F^Z(\pi(t)),\ (Z, \pi) \in X_{\delta}.
\ee
Note that $\|Z(t)\|_{\Sigma}$ grows at most linearly in time, so if the condition $\|Z(0)\|_{\Sigma} <\delta_0$ from the definition (\ref{defN}) of $\mc N$ is met then it still holds for $\tau$ close to $0$. Assume that $\delta_0$ is sufficiently small that by H\"older's inequality
\be
\|Z(t)\|_{L^{4/3} \cap H^{1/2}} \leq C \|Z(t)\|_{\Sigma} < C\delta_0.
\ee
Therefore we can apply Lemma \ref{comparatie} and obtain that $\Psi(\tau) = \Psi_{\tau}(0) \in \mc N_{L^{q} \cap H^{1/2}}$. Since $\|Z(\tau)\|_{\Sigma} < \delta_0$, it follows that $\Psi(\tau) \in \mc N$.

The same proof cannot yield global in time invariance, because after a strictly positive time the solution may escape the $\delta_0$-neighborhood of the standing wave manifold in which it has to be for the argument to hold. Existence results for initial data in a weaker invariant norm, one that does not grow in time, are needed in order to approach the global result.

Similar considerations apply to the proof of local in time invariance for negative time. Consider a global solution $\Psi$ with $\Psi(0) \in \mc N$. Instead of truncating, now we use local existence theory to get a continuation of $\Psi$ to a small interval $[-t, 0]$. Let $\Psi_{-t}$ be the solution obtained by pasting this onto the original solution. By means of the modulation equations
\be
\Big\langle\bpm \Psi(s) \\ \ov{\Psi(s)} \epm - \bpm W(\pi(s)) \\ \ov{W(\pi(s))} \epm, \xi_F(\pi(s)) \Big\rangle = 0
\ee
we also extend the parameter path $\pi$ to a small time interval $[-t, 0]$. At this point we apply Lemma \ref{comparatie} exactly as above, eventually concluding that $\Psi(-t) \in \mc N$.

There is still the issue of checking that $\mc N$ is a centre-stable manifold as in \cite{bates}. We prove it in the neighborhood of each standing wave $W$, since being a centre-stable manifold is a local property.

To begin with, we rewrite equation (\ref{NLS}) to make it fit the framework of the theory of Bates, Jones \cite{bates}. Consider a fixed ground state $\phi(\cdot, \alpha_0)$ and the constant path $\pi_0 = (0, 0, 0, \alpha_0)$. Linearizing the equation around this constant path and applying a symmetry transformation, as in Lemma \ref{transformare}, yields, for $U = \frak g_{\pi_0} \Big(\bpm \Psi \\ \ov{\Psi} \epm - \bpm W(\pi_0) \\ \ov{W(\pi_0)} \epm\Big)$, that
\be
i \partial_t U + \mc H U = N(U, \pi_0),
\lb{eq_2188}
\ee
where
\be
\mc H = \bpm \Delta+2\phi^2(\cdot, \alpha_0) -\alpha_0 & \phi^2(\cdot, \alpha_0) \\ -\phi^2(\cdot, \alpha_0) & -\Delta-2\phi^2(\cdot, \alpha_0)+\alpha_0 \epm
\ee
and $N(U, \pi_0)$ is as in (\ref{defNUpi}). Note that here all the right-hand side terms are at least quadratic in $U$, due to linearizing around a constant path.

The spectrum of $\mc H$ is known, see Section \ref{spectru}, namely $\sigma(\mc H) = (-\infty, -\alpha] \cup [\alpha, \infty) \cup \{0, \pm i \sigma\}$. The stable spectrum is $-i\sigma$, the unstable spectrum is $i\sigma$, and everything else belongs to the centre. It is easy to check that all the conditions of \cite{bates} are met, leading to the existence of a centre-stable manifold.

In the sequel we prove that $\mc N$ is a centre-stable manifold, namely that it fulfills the three three properties enumerated in Definition \ref{centr}: $\mc N$ is $t$-invariant with respect to a neighborhood of $\phi(\cdot, \alpha_0)$, $\pi^{cs}(\mc N)$ contains a neighborhood of $0$ in $X^c \oplus X^s$, and $\mc N \cap W^u = \{0\}$. All of this is relative to a specific neighborhood of $0$, $\mc V = \{U \mid \|U\|_{\Sigma} < \delta_0\}$ for some small $\delta_0$.

The $t$-invariance of $\mc N$ relative to $\mc V$ follows from definition and Lemma \ref{comparatie}, in the same manner in which we proved the local in time invariance of $\mc N$ under the Hamiltonian flow.

Then,
\be
\pi^{cs}(\mc N) = \{(P_- + P_c + P_{root})(U(0))\} = \Big\{\frak g \Big(\bpm \Psi_0 \\ \ov{\Psi_0} \epm\Big) - \phi(\cdot, \alpha_0)\Big\}.
\ee
Since $\bpm R_0 \\ \ov{R_0} \epm$ covers a whole neighborhood of $0$ in ${\rm Ker}(P_+) \cap {\rm Ker}(P_{root})$ and the action of $\frak g$ is transverse, with the range of its differential at $0$ spanned by $\eta_F(\pi_0)$, it follows that $\pi^{cs}(\mc N) $ is a neighbothood of $0$ in ${\rm Ker}(P_+)$, the second property that the centre-stable manifold must have.

Finally, consider a solution $U \in W^u$ of (\ref{eq_2188}), meaning that $\|U(t)\|_{\Sigma} \leq \delta$ for all negative $t$ and that it decays exponentially as $t \to -\infty$, $\|U(t)\|_{\Sigma} \leq C e^{Ct}$ (even though polynomial decay is sufficient). Decompose $U$ into its projections on the continuous, imaginary, and zero spectrum of $\mc H$ and let
\be
\delta(T) = \|U\|_{\langle t \rangle^{-1} L^1_t(-\infty, T] L^{6/5 \cap 6}_x} + \|U\|_{L^2_t(-\infty, T] L^{6}_x \cap L^{\infty}_t(-\infty, T] L^2_x}.
\ee
Observe that $\delta(t) \to 0$ as $t \to -\infty$, so we can assume it to be arbitrarily small, $\delta(t) < 1$ to begin with.

By means of Strichartz estimates one obtains that
\be
\begin{aligned}
\|P_c U\|_{L^2_t(-\infty, T] L^{6}_x \cap L^{\infty}_t(-\infty, T] L^2_x} &\leq C \|N(U, \pi_0)\|_{L^2_t(-\infty, T] L^{6/5}_x + L^{1}_t(-\infty, T] L^2_x} \\
& \leq C \delta(T) \|U\|_{L^{\infty}_t (-\infty, T] L^2_x},
\end{aligned}
\ee
because now the right-hand side $N$ contains only quadratic or higher terms.

The same estimate holds for $P_{im} U$, because it is bounded at $-\infty$, so we can use Lemma \ref{hyp}. We write it in the form
\be\begin{aligned}
P_- U(t) &= -\int_{-\infty}^t e^{-\sigma(t-s)} P_- N(U(s), \pi_0) \dd s\\
P_+ U(t) &= e^{(t-T)\sigma} P_+U(T) - \int_{t}^T e^{(t-s)\sigma} P_+ N(U(s), \pi_0) \dd s.
\end{aligned}\lb{PplusU}\ee
Therefore
\be
\|P_{im} U\|_{L^{\infty}_t(-\infty, T] L^2_x} \leq C (\|P_+ U(T)\|_2 + \delta(T) \|U\|_{L^{\infty}_t(-\infty, T] L^2_x}).
\ee

We assumed no orthogonality condition. The modulation equations (\ref{mod}) now give $P_{root} U$ and contain only quadratic or higher terms. Note that $\ds\lim_{t \to -\infty} P_{root} U(t) = 0$. Therefore
\be
\|P_{root} U\|_{L^{\infty}_t(-\infty, T] L^2_x} \leq C \|N(U, \pi_0)\|_{\langle t \rangle^{-1} L^1_t(-\infty, T] L^{1+\infty}_x} \leq C \delta(t) \|U\|_{L^{\infty}_t(-\infty, T] L^2_x}.
\ee

Putting these estimates together, one has that
\be
\|U\|_{L^{\infty}_t(-\infty, T] L^2_x} \leq C (\delta(T) \|U\|_{L^{\infty}_t(-\infty, T] L^2_x} + \|P_+ U(T)\|_{2}).
\ee
For sufficiently negative $T_0$, it follows that $\|U(t)\|_2 \leq C \|P_+ U(t)\|_2$, for any $t \leq T_0$. The converse is obviously true, so the two norms are comparable.

Furthermore, by reiterating this argument one has that $\|(1-P_+)U(t)\|_2 \leq C \delta(t) \|P_+ U(t)\|_2$.

Next, assume that $U$ is on the stable manifold, meaning that, for some $\Psi(0) \in \mc N$, $U = \frak g_{\pi_0} \Big(\bpm \Psi \\ \ov{\Psi} \epm - \bpm W(\pi_0) \\ \ov{W(\pi_0)} \epm\Big)$. $\Psi$ has some moving soliton path $\pi$ associated to it such that $\Big(\bpm \Psi \\ \ov{\Psi} \epm - \bpm W(\pi) \\ \ov{W(\pi)} \epm=Z, \pi\Big) \in X_{\delta}$. Then,
\begin{align}
|\pi(t)-\pi_0(t)| \leq &C\Big(\Big\|\bpm \Psi(t) \\ \ov{\Psi(t)} \epm - \bpm W(\pi(t)) \\ \ov{W(\pi(t))} \epm\Big\|_{\Sigma} +  \Big\|\bpm \Psi(t) \\ \ov{\Psi(t)} \epm - \bpm W(\pi_0(t)) \\ \ov{W(\pi_0(t))} \epm\Big\|_{\Sigma}\Big)\nonumber\\
\leq &C (\delta+|\pi(\infty) - (0, 0, 0, \alpha_0)|)
\end{align}
and therefore $\|\pi -\pi_0\|_{\infty} \leq C(\delta+|\pi(\infty) - (0, 0, 0, \alpha_0)|)$. This implies that
\be
\|U(t)\|_{L^{\infty}_t L^{2\cap 3}_x} \leq C (\delta+|\pi(\infty) - (0, 0, 0, \alpha_0)|)
\ee
for positive $t$ and therefore for all $t$, due to the fact that $U \in W^u$.

If $(0, 0, 0, \alpha_0)$ is not the terminal value of the path $\pi$, then, since $U$ is a bounded solution to (\ref{eq_2188}) at $+\infty$, it follows by Lemma \ref{hyp} that
\be
\|P_+ U(t)\|_2 \leq \int_t^{\infty} e^{(t-s)\sigma} \|P_+ N(U, \pi)\|_2 ds \leq C (\delta^2 + |\pi(\infty) - (0, 0, 0, \alpha_0)|^2).
\ee
Also, because $\|P_{root} U(t)\|_2$ is bounded from below,
\be
\|(1-P_+) U(t)\|_2 \geq \|P_{root} U(t)\|_2 \geq C |\pi(\infty) - (0, 0, 0, \alpha_0)|.
\ee
Therefore
\be
\|P_+ U(t)\|/\|(1-P_+)U(t)\| \leq C(\delta + |\pi(\infty) - (0, 0, 0, \alpha_0)|) \leq C(\delta+\delta_0).
\ee
If, on the other hand, $\pi$ approaches $(0, 0, 0, \alpha_0)$ in the limit, note that in the modulated equation (\ref{nlsu}) valid for $(U_Z = \frak g_{\pi} Z, \pi)$, if $\delta$ is sufficiently small, it follows that
\be
C\|P_c U_Z(0)\|_{\Sigma} \geq \|P_c U_Z(t)\|_2 \geq C(\|P_c U_Z(0)\|_{\Sigma} - \|P_c U_Z(0)\|_{\Sigma}^2) \geq C \|P_c U_Z(0)\|_{\Sigma} \geq C \delta
\ee
by means of the Strichartz estimates. Such a lower bound is also implied by scattering. Furthermore, since $\pi$ and $\pi_0$ grow near in the $+\infty$ limit, it follows that $\|P_c U_Z - P_c U\|_2 \to 0$.

Therefore the norm of $(1-P_+)U$ is bounded from below at $+\infty$ in either case, unless $P_c Z(0) = 0$. However, if it were so, it would imply that $Z=0$ and $\Psi$ is constant, equal to a soliton. Since $\Psi - W(\pi_0)$ decays exponentially at $-\infty$, this would imply that they are equal and then $U=0$.

For $P_+U$, now we have that $\|P^+_{\pi} U_Z(t) - P^+_{\pi_0} U(t)\|_{1 \cap \infty} \to 0$ as $t \to \infty$ and therefore $\|P_+ U\|_{1 \cap \infty} \to 0$.

Excluding the trivial case $U=0$, we obtain that either $\|P_+ U(t)\|/\|(1-P_+)U(t)\| \leq C(\delta + \delta_0)$ or $\|P_+ U(t)\|/\|(1-P_+)U(t)\| \to 0$ as $t \to \infty$. This implies that one can make this ratio as small as necessary for some large $T_0$.

Assume that $\delta<1$. Then Lemma 2.4 from \cite{bates} states, under even more general conditions, that if the ratio $\|P_+ U(T_0)\|/\|(1-P_+)U(T_0)\|$ is small enough, it will stay bounded for all $t \leq T_0$. The proof of this result is based on Gronwall's inequality.

However, this contradicts our previous conclusion that $\|(1-P_c)U(t)\|_2/\|P_c U(t)\|_2$ goes to $0$ as $t$ goes to $-\infty$. Therefore, $U$ can only be $0$.

This proves that $\mc N \cap W^u = \{0\}$. In other words, there are no exponentially unstable solutions in $\mc N$ in the sense of \cite{bates}. The final requirement for $\mc N$ to be a centre-stable manifold is thus met.
\end{proof}

\section{Linear Estimates}\lb{lin_est}
\subsection{The endpoint Strichartz estimate}
Consider operators in $\set R^3$ of the form $\mc H = \mc H_0 + V$, where
\be \mc H_0 = \bpm -\Delta + \mu & 0\\0&\Delta-\mu\epm,\ V = \bpm -U & -W \\ W & U\epm.
\lb{presupuneri}
\ee
We assume that $-\sigma_3 V$ is a positive matrix, that $L_- = -\Delta + \mu + U + W \geq 0$, that $|V|  \leq C \langle x \rangle^{-7/2-}$, that the spectral Assumption \ref{assum} holds, and that the edges of the spectrum $\pm \mu$ are neither eigenvalues nor resonances.

The operator $\mc H$ has $\sigma(\mc H) \subset \set R \cup i\set R$ and $\sigma_{ess}(\mc H) = (-\infty, -\mu] \cup [\mu, \infty)$. We make the spectral assumption that $\mc H$ has no eigenvalues in the set $(-\infty, -\mu) \cup (\mu, \infty)$ and that the thresholds $\pm \mu$ are also regular, meaning that $I+(\mc H_0-\mu\pm i0)^{-1}V:\langle x \rangle^{1+\epsilon}L^2 \to \langle x \rangle^{1+\epsilon}L^2$ is invertible.

Firstly, we need the main result of \cite{tao}:
\begin{theorem}
Let $(X, d\mu)$ be a measure space, $L^p_x = L^p(X, d\mu)$. Suppose that for each $t$ one has  an operator $U(t)$ such that
\be
\|U(t)f\|_2 \leq C \|f\|_2,\ \|U(s) U^*(t) f\|_{\infty} \leq C |t-s|^{-\sigma} \|f\|_1.
\ee
Let $\sigma>1$. Call $(q, r)$ sharp $\sigma$-admissible if $q, r \geq 2$, $\ds\frac 1 q + \frac {\sigma} r = \frac \sigma 2$, and let $q'$ be the exponent such that $\frac 1 {q'} + \frac 1 q = 1$. Then
\begin{align}
\|U(t) f\|_{L^q_t L^r_x} \leq& C\|f\|_2,\lb{r1}\\
\ds\Big\|\int U^*(s) F(s) \dd s\Big\|_2 \leq& C\|F\|_{L^{q'}_t L^{r'}_x},\lb{r2}\\
\ds\Big\|\int_{s<t} U(t)U^*(s) F \dd s\Big\|_{L^q_t L^r_x} \leq& C\|F\|_{L^{\tilde q'}_t L^{\tilde r'}_x}\lb{r3},
\end{align}
for any sharp $\sigma$-admissible $(q, r)$, $(\tilde q, \tilde r)$.
\lb{tao}
\end{theorem}

The following lemma is a straightforward generalization of Theorem \ref{tao}. Since it is important in the sequel, however, a short proof will be given.
\begin{lemma}
Let $(X, d\mu)$ be a measure space, $L^p_{\mu} = L^p(X, d\mu)$. Suppose that for each $t$ one has operators $U(t)$ and $V(t)$ that satisfy
\be
\|U(t)f\|_2 \leq C \|f\|_H,\ \|V(t)f\|_2 \leq C \|f\|_H
\ee
and
\be\ba{c}
\|U(s) V(t) f\|_{\infty} \leq C |t-s|^{-\sigma} \|f\|_1, \\
\|U(s) U^*(t) f\|_{\infty} \leq C |t-s|^{-\sigma} \|f\|_1,\\
\|V^*(s) V(t) f\|_{\infty} \leq C |t-s|^{-\sigma} \|f\|_1.
\ea\ee
Let $\sigma>1$. Call $(q, r)$ sharp $\sigma$-admissible if $q, r \geq 2$, $\ds\frac 1 q + \frac {\sigma} r = \frac \sigma 2$, and let $q'$ be the exponent such that $\frac 1 {q'} + \frac 1 q = 1$. Then, in addition to the Strichartz estimates (\ref{r1}-\ref{r3}) for $U$ and $V$, one has that
\be
\ds\Big\|\int_{s<t} U(t)V(s) F \dd s\Big\|_{L^q_t L^r_x} \leq C\|F\|_{L^{\tilde q'}_t L^{\tilde r'}_x} \lb{retard2}
\ee
for any sharp $\sigma$-admissible $(q, r)$, $(\tilde q, \tilde r)$.
\lb{St}
\end{lemma}
\begin{proof}
The proof is a rephrasing of the one given in \cite{tao}. Inequalities (\ref{r1}-\ref{r3}) are already provided, so only (\ref{retard2}) is left.
Consider the bilinear form
\be
T(F, G) = \iint_{s<t} \langle U^*(s) F(s), V(t) G(t)\rangle \dd s \dd t.
\ee
By interpolation between
\be
|\langle U^*(s) F(s), V(t) G(t)\rangle| \leq C \|F(s)\|_2 \|G(t)\|_2
\lb{end1}
\ee
and
\be
|\langle U^*(s) F(s), V(t) G(t)\rangle| \leq C |t-s|^{-\sigma} \|F(s)\|_1 \|G(t)\|_1
\ee
we obtain
\be
|\langle U^*(s) F(s), V(t) G(t)\rangle| \leq C |t-s|^{-1-\beta(r, r)} \|F(s)\|_{r'} \|G(t)\|_{r'}
\ee
where $\ds\beta(r, \tilde r) = \sigma - 1 - \frac \sigma r - \frac \sigma {\tilde r}$.

Let $\ds T_j(F, G) = \iint_{t-2^{j+1}<s\leq t-2^j} \langle U^*(s) F(s), V(t) G(t)\rangle \dd s \dd t$. Then the estimate
\be
|T_j(F, G)| \leq C 2^{-j\beta(a, b)} \|F\|_{L^2_t L^{a'}_x} \|G\|_{L^2_t L^{b'}_x}.\lb{bilin}
\ee
holds for all $j \in\set Z$ and all $(\frac 1 a, \frac 1 b)$ in a neighborhood of $(\frac 1 r, \frac 1 r)$. The proof goes through showing (\ref{bilin}) for the exponents $a=b=\infty$, $(a, 2)$ with $2\leq a < r$, and $(2, b)$ with $2\leq b <r$.

One can now infer that
\be
|T(F, G)| \leq \|F\|_{L^{q'}_t L^{r'}_x} \|G\|_{L^{q'}_t L^{r'}_x}
\lb{end2}
\ee
for the endpoint $\ds(q, r) = (2, \frac {2\sigma}{\sigma-1})$. Since (\ref{end2}) was already true for the other endpoint by (\ref{end1}), it is true for all admissible $(q, r)$. The general retarded estimate (\ref{retard2}) follows immediately as in \cite{tao} by interpolation between (\ref{end1}), (\ref{end2}), and
\be
|T(F, G)| \leq \sup_t \|\int_{s<t} U^*(s) F(s) \dd s\|_2 \|G\|_{L^1_t L^2_x} \leq \|F\|_{L^{q'}_t L^{r'}_x} \|G\|_{L^1_t L^2_x},
\ee
and the symmetric inequality, which are both consequences of (\ref{r2}).
\end{proof}

Applying Lemma \ref{St} to the families of evolution operators $U(t)=e^{it\mc H} P_c$, $V(s) = e^{-is \mc H} P_c$, we have obtained
\begin{corollary}
Assume that
\be
\|e^{it \mc H} P_c\|_{2 \to 2} \leq C
\ee
and
\be\ba{c}
\|e^{it \mc H^*} P_c^* e^{-is\mc H} P_c \|_{1 \to \infty} \leq C|t-s|^{-3/2},\\
\|e^{it \mc H} P_c e^{-is\mc H^*} P_c^*\|_{1 \to \infty} \leq C|t-s|^{-3/2}.
\ea\ee
Then
\be\begin{aligned}
\|e^{it \mc H} P_c f\|_{L^q_t L^r_x} \leq &C\|f\|_2,\\
\ds\Big\|\int e^{-is \mc H} P_c F(s) \dd s\Big\|_2 \leq &C\|F\|_{L^{q'}_t L^{r'}_x},\\
\ds\Big\|\int_{s<t} e^{it \mc H} P_c e^{-isH^*} P_c^* F \dd s\Big\|_{L^q_t L^r_x} \leq &C\|F\|_{L^{\tilde q'}_t L^{\tilde r'}_x},\\
\ds\Big\|\int_{s<t} e^{i(t-s) \mc H} P_c F \dd s\Big\|_{L^q_t L^r_x} \leq &C\|F\|_{L^{q'}_t L^{r'}_x}.
\lb{strichartz}
\end{aligned}\ee
for any sharply admissible $(q, r)$ (that is, such that $q$, $r \geq 2$, $\ds\frac 1 q + \frac 3 {2r} = \frac 3 4$) and $(\tilde q, \tilde r)$. The same estimates hold after swapping $\mc H$ and $\mc H^*$.
\end{corollary}

Coming back to the particular operator $\mc H$ given in (\ref{presupuneri}), let $P_p$ be the Riesz projection on the point spectrum and $1-P_p = P_c$ be the projection on the continuous spectrum of $\mc H$. The evolution applied to eigenfunctions may lead to exponential growth in any norm; leaving that aside, one can achieve the above bound for $\mc H P_c$. Indeed, the bounds
\be
\|e^{it\mc H} P_c \|_{2 \to 2} \leq C
\lb{l2}
\ee
and
\be
\|e^{it\mc H} P_c \|_{1 \to \infty} \leq C|t|^{-3/2}
\lb{I1}
\ee
were proved by Schlag in \cite{schlag} and Schlag, Erdogan in \cite{erdogan2}. What is left to prove is
\be
\|e^{it\mc H} P_c e^{-is\mc H^*} P_c^*\|_{1 \to \infty} \leq C|t-s|^{-3/2},
\lb{I2}
\ee
as well as the symmetric estimate.

\subsection{Proof of the strengthened dispersive estimate}

\begin{proof}
We begin with the following explicit representation derived from \cite{schlag}: for $f$, $g \in L^{2, 1+}$,
\be
\langle \phi, e^{it\mc H} P_c \psi \rangle= \frac 1 {2\pi i} \int_{\Gamma^-_{\epsilon} \cup \Gamma^+_{\epsilon}} e^{it\lambda} \langle \phi, R_V(\lambda) \psi \rangle \dd \lambda,
\ee
where $R_V(\lambda) = (\mc H - \lambda)^{-1}$ and $\Gamma^{\pm}_{\epsilon}$ are the counterclockwise contours given by $\Gamma^{\pm}_{\epsilon} = \{z \mid d(z, [\pm\mu, \pm\infty)) = \epsilon\}$. The integral can be taken improper, but it is more helpful to consider instead the mollified version
\be
\int_0^{\infty} f(\lambda) \dd x = \lim_{R \to \infty} \int_0^{\infty} f(\lambda) \chi(\lambda/R) \dd \lambda,
\ee
where $\chi$ is a smooth cutoff function with $\chi(x) = 1$ for $|x|<1$ and $\chi(x) = 0$ for $|x|>2$.
A similar formula holds for $e^{-is\mc H^*}P_c^*$, with $V$ replaced by $V^*$. The expression we need to estimate becomes
\be
\langle e^{-it\mc H^*} P_c^* \phi,  e^{-is\mc H^*} P_c^* \psi\rangle = -\frac {1}{4\pi^2} \int_{\Gamma^-_{\epsilon_1} \cup \Gamma^+_{\epsilon_1}} \int_{\Gamma^-_{\epsilon_2} \cup \Gamma^+_{\epsilon_2}} e^{i(t\lambda - s\eta)} \langle \phi, R_V(\lambda) R_{V^*}(\eta) \psi \rangle \dd \eta \dd \lambda.
\ee
We make the arbitrary choice $\epsilon_1 > \epsilon_2$. After splitting each contour into $\Gamma^+$ and $\Gamma^-$, we obtain $4$ terms to be treated separately. We begin with the $\Gamma^+ \Gamma^+$ term. Expand both $R_V$ and $R_{V^*}$ into finite Born sums consisting of $2m$ terms and a remainder. Let $R_0(\lambda) = (\mc H_0 - \lambda)^{-1}$. The expression becomes
\be
\langle\phi, e^{it\mc H} P_+ e^{-is\mc H^*} P_+^* \psi\rangle=
\ee
\be
=\sum_{\ell=0}^{2m-1} \sum_{k=0}^{2m-1} (-1)^{k+\ell} \int_{\Gamma^+_{\epsilon_1}} \int_{\Gamma^+_{\epsilon_2}} e^{i(t\lambda - s\eta)} \big\langle \phi, R_0(\lambda)(VR_0(\lambda))^{\ell} R_0(\eta)(V^*R_0(\eta))^k \psi\big\rangle \dd \eta \dd \lambda \lb{unu}
\ee
\begin{multline}
+ \sum_{\ell=0}^{2m-1} (-1)^{\ell} \int_{\Gamma^+_{\epsilon_1}} \int_{\Gamma^+_{\epsilon_2}} e^{i(t\lambda-s\eta)} \big\langle\phi, \\
R_0(\lambda)(VR_0(\lambda))^{\ell} (R_0(\eta)V^*)^m R_{V^*}(\eta) (V^* R_0(\eta))^m \psi\big\rangle \dd \eta \dd \lambda \lb{doi}
\end{multline}
\begin{multline}
+ \sum_{k=0}^{2m-1} (-1)^k \int_{\Gamma^+_{\epsilon_1}} \int_{\Gamma^+_{\epsilon_2}} e^{i(t\lambda-s\eta)} \big\langle\phi, \\
(R_0(\lambda)V)^m R_{V}(\lambda) (V R_0(\lambda))^m R_0(\eta)(V^*R_0(\eta))^k \psi\big\rangle \dd \eta \dd \lambda \lb{trei}
\end{multline}
\begin{multline}
+ \int_{\Gamma^+_{\epsilon_1}} \int_{\Gamma^+_{\epsilon_2}} e^{i(t\lambda-s\eta)} \big\langle\phi, \\
(R_0(\lambda)V)^m R_{V}(\lambda) (V R_0(\lambda))^m (R_0(\eta)V^*)^m R_{V^*}(\eta) (V^* R_0(\eta))^m \psi\big\rangle \dd \eta \dd \lambda.\lb{patru}
\end{multline}

In each term, the localizing potentials $V$ or $V^*$ alternate with the resolvents $R_0$, with the exception of exactly two resolvent operators following one another. Since this is a potentially dangerous situation, we apply the resolvent identity. For the very simplest term, this means
\begin{multline}
\int_{\Gamma^+_{\epsilon_1}} \int_{\Gamma^+_{\epsilon_2}} e^{i(t\lambda-s\eta)} \big\langle\phi, R_0(\lambda) R_0(\eta) \psi \big\rangle \dd \eta \dd \lambda = \\
=\int_{\Gamma^+_{\epsilon_1}} \int_{\Gamma^+_{\epsilon_2}} \frac{e^{i(t\lambda-s\eta)}}{\lambda-\eta} \big\langle\phi, R_0(\eta) \psi \big\rangle \dd \eta \dd \lambda - \int_{\Gamma^+_{\epsilon_1}} \int_{\Gamma^+_{\epsilon_2}} \frac{e^{i(t\lambda-s\eta)}}{\lambda-\eta} \big\langle\phi, R_0(\lambda) \psi \big\rangle \dd \eta \dd \lambda.
\end{multline}
Every resulting term has a kernel that can be written, by means of the resolvent identity, as a sum of two parts of the form
\be
T=\int_{\Gamma^+_{\epsilon_1}} \int_{\Gamma^+_{\epsilon_2}} \frac{e^{i(t\lambda-s\eta)}}{\lambda-\eta} f(\lambda) g(\eta) \dd \eta \dd \lambda.
\lb{expr_gen_T}
\ee

We make $f$ and $g$ explicit later, but for now we continue in this general setting. 
After making $\epsilon_2=0$, we get
\be
T=\int_{\Gamma^+_{\epsilon_1}} \int_0^{\infty} \frac{e^{i(t\lambda-s(\eta+\mu))}}{\lambda-(\eta+\mu)} f(\lambda) (g(\eta+\mu+i0)-g(\eta+\mu-i0)) \dd \eta \dd \lambda.
\ee
Let $f_s(\lambda) = \frac 1 2(f(\lambda+i0)+f(\lambda-i0))$ and $f_a(\lambda) = \frac 1 2(f(\lambda+i0)-f(\lambda-i0))$ and same for $g$. We assume that $g(\lambda \pm i \epsilon) = g(\lambda \pm i0) + O(\epsilon)$ uniformly on compact intervals. Letting $\epsilon_1$ go to $0$ we have
\be\ba{ll}
T&=\ds\int_0^{\infty} e^{it(\lambda+\mu)} f_s(\lambda+\mu) \lim_{\epsilon \to 0} \Big(\int_0^{\infty} e^{-is(\eta + \mu)} g_a(\eta+\mu) \frac \epsilon {(\lambda-\eta)^2+\epsilon^2} \dd \eta\Big) \dd \lambda+ \\
&\ds+\int_0^{\infty} e^{it(\lambda+\mu)} f_a(\lambda+\mu) \lim_{\epsilon \to 0} \Big(\int_0^{\infty} e^{-is(\eta + \mu)} g_a(\eta+\mu) \frac {\lambda-\eta} {(\lambda-\eta)^2+\epsilon^2} \dd \eta\Big) \dd \lambda\\
&\ds= e^{i(t-s)\mu} \int_0^{\infty} \Big(e^{i(t-s)\lambda} f_s(\lambda+\mu) g_a(\lambda+\mu) + \\
&\ds+ e^{it\lambda} f_a(\lambda+\mu) H\big(\chi_{[0, \infty)} e^{-is\eta} g_a(\eta+\mu)\big)(\lambda)\Big) \dd \lambda,
\ea\lb{HP}\ee
where $H$ is the Hilbert transform. The limit exists because we are applying singular kernels to integrable functions of compact support. Any further error terms are in the order of $\epsilon$ and vanish.

Therefore, we need to examine oscillatory integrals of the following form:
\begin{lemma} Assume that $F_s$ is even and $F_a$ and $G_a$ are odd functions and that $\widehat {F_s}, \widehat {F'_s} \in \mc M$ (the space of finite measures) and likewise for $F_a$ and $G_a$. Then
\be
\int_0^{\infty} e^{i(t-s)\lambda} F_s(\sqrt \lambda) G_a(\sqrt\lambda) \dd \lambda \leq C |t-s|^{-3/2} (\|\widehat {F_s}\|_1 \|\widehat {G'_a}\|_1 + \|\widehat {F'_s}\|_1 \|\widehat {G_a}\|_1), \lb{poisson}
\ee
as well as
\begin{multline}
\int_0^{\infty} e^{it\lambda} F_a(\sqrt \lambda) H(e^{-is\eta} \chi_{[0, \infty)} G_a(\sqrt\eta))(\lambda) \dd \lambda \leq \\
\leq C|t-s|^{-3/2} (\|\widehat {F_a}\|_1 \|\widehat {G'_a}\|_1 + \|\widehat {F'_a}\|_1 \|\widehat {G_a}\|_1). \lb{hilbert}
\end{multline}
\lb{baza}
\end{lemma}

The integrals on the left-hand side are improper and computed with the help of a smooth cutoff. Under the assumptions, $F_a \in C^1(\set R)$, so $\chi_{[0, \infty)} F_a(\sqrt\lambda)\in C^{1/2}$, and likewise for $G_a$.

\begin{proof} For (\ref{poisson}) we have
\be\ba{ll}
(\ref{poisson}) &\ds= \int_{-\infty}^{\infty} {e^{i(t-s)\lambda^2} \lambda F_s(\lambda) G_a(\lambda) \dd \lambda} \\
&\ds=\frac 1 {t-s} \int_{-\infty}^{\infty} {e^{i(t-s)\lambda^2} (F_s(\lambda) G_a(\lambda))' \dd \lambda}\\
&\ds\leq C |t-s|^{-3/2} (\|\widehat {F_s}\|_1 \|\widehat {G'_a}\|_1 + \|\widehat {F'_s}\|_1 \|\widehat {G_a}\|_1).
\ea
\ee

For terms of the form (\ref{hilbert}) a slight refinement is needed. First note that
\begin{multline}
(\chi_{[0, \infty)} F_a(\sqrt \cdot))^{\wedge}(\tau) = \int_0^{\infty} e^{-i\tau\lambda} F_a(\sqrt \lambda) \dd \lambda = \frac 1 2 \int_{-\infty}^{\infty} e^{-i\tau\lambda^2} \lambda F_a(\lambda) \dd \lambda \\
= \frac 1 {2\tau} \int_{-\infty}^{\infty} e^{-i\tau\lambda^2} F'_a(\lambda) \dd \lambda = \frac 1 {4\sqrt 2 \pi^{3/2}} \tau^{-3/2} \int_{-\infty}^{\infty} e^{i\xi^2/\tau} \widehat {F'_a}(\xi) \dd \xi.
\end{multline}
The same goes for $G_a$. Letting $s-t=k$ we have
\begin{multline}
(\ref{hilbert}) = \int_{-\infty}^{\infty} (e^{it\cdot} \chi_{[0, \infty)} F_a(\sqrt \cdot))^{\wedge}(\tau) (H(e^{-is\cdot} \chi_{[0, \infty)} G_a(\sqrt \cdot)))^{\wedge}(-\tau) \dd \tau\\
= \int_{-\infty}^{\infty} (\chi_{[0, \infty)} F_a(\sqrt \cdot))^{\wedge}(\tau-t) (\chi_{[0, \infty)} G_a(\sqrt \cdot))^{\wedge}(-\tau+s) \sgn(\tau) \dd \tau\\
\leq 2 \sup_{a, b} \bigg|\int_a^b (\chi_{[0, \infty)} F_a(\sqrt \cdot))^{\wedge}(\tau) (\chi_{[0, \infty)} G_a(\sqrt \cdot))^{\wedge}(k-\tau) \dd \tau\bigg| \\
= \frac 1 {16\pi^3} \sup_{a, b} \bigg|\int_a^b \int_{-\infty}^{\infty} \int_{-\infty}^{\infty} |\tau|^{-3/2} |k-\tau|^{-3/2} e^{i(\frac {\lambda^2}{\tau} + \frac {\eta^2}{k-\tau})} \widehat {F_a'}(\xi) \widehat {G_a'}(\nu) \dd \nu \dd \xi \dd \tau\bigg|.
\end{multline}
We apply the stationary phase method to this integral (see, for example, \cite{stein}, p. 332). We write the proof explicitly because neither the phase, nor the integrand is absolutely integrable.

Consider $\ds\int_a^b e^{i\phi} \psi \dd \tau$, where $\ds\phi = \frac {\xi^2}{\tau} + \frac {\nu^2}{k-\tau}$ and $\psi = |\tau|^{-3/2} |k-\tau|^{-3/2}$. The aim is to prove that
\be
\bigg|\int_a^b e^{i\phi} \psi \dd \tau\bigg| \leq C |k|^{-3/2}\Big(\frac 1 {|\xi|}+\frac 1 {|\nu|}\Big).
\lb{model}
\ee

Without loss of generality, let $\xi$, $\nu$, $k>0$. First assume $[a, b] \subset [0, k]$ and note that
\be
\left|\int_{\tau_1}^{\tau_2} e^{i\phi} \psi \dd \tau\right| \leq \int_{\tau_1}^{\tau_2} \frac{\dd \tau}{\tau^{3/2} (k-\tau)^{3/2}} = \Psi_1(\tau) \Big|_{\tau_1}^{\tau_2},
\lb{abso}
\ee
where the antiderivative is $\ds\Psi_1(\tau) = \frac 1 {2k^2} (\tau^{1/2} (k-\tau)^{-1/2} - \tau^{-1/2} (k-\tau)^{1/2})$. Indeed,
\be
\begin{aligned}
(\tau^{1/2} (k-\tau)^{-1/2} - \tau^{-1/2} (k-\tau)^{1/2})' = &\frac 1 2 (\tau^{-1/2}(k-\tau)^{-1/2} + \tau^{1/2}(k-\tau)^{-3/2} + \\
&+\tau^{-3/2} (k-\tau)^{1/2} + \tau^{-1/2} (k-\tau)^{-1/2})\\
=& \frac{(k-\tau+\tau)^2}{2\tau^{3/2} (k-\tau)^{3/2}}.
\end{aligned}
\ee
On any interval not containing a stationary point, moreover, one has
\be
\bigg|\int_{\tau_1}^{\tau_2} e^{i\phi} \psi \dd \tau\bigg| \leq C \sup_{\tau \in [\tau_1, \tau_2]} \Big| \frac {\psi(\tau)} {\phi'(\tau)} \Big|,
\lb{part}
\ee
where, by the convexity of $1/x$,
\be
\Big|\frac {\psi(\tau)}{\phi'(\tau)} \Big|= \frac 1 {\nu\tau+\xi(k-\tau)} \frac{\tau^{1/2}(k-\tau)^{1/2}}{|\nu \tau-\xi (k-\tau)|} \leq \frac 1 k \Big(\frac 1 \nu + \frac 1 \xi \Big) \frac {\tau^{1/2}(k-\tau)^{1/2}}{|\nu \tau-\xi (k-\tau)|}.
\lb{344}
\ee
The reason is that, after integration by parts, $\ds \Big(\frac {\psi(\tau)}{\phi'(\tau)}\Big)'$ changes sign at most a constant number of times, so one can integrate and only lose some constant.

The last expression in $(\ref{344})$ is $0$ at the endpoints $0$ and $k$. Note that the phase derivative $\phi'$ vanishes at exactly one point in the interval $[0, k]$, namely $\tau_0 = \ds\frac{\xi k}{\nu+\xi}$. Surround $\tau_0$ with a small interval $[\tau_1, \tau_2]$ on which we estimate the integral in absolute value and otherwise integrate by parts as above. We obtain that
\be
\bigg|\int_{a}^{b} e^{i\phi} \psi \dd \tau\bigg| \leq C \Big(\max_{\tau\in[0, \tau_1]\cup [\tau_2, k]} \Big|\frac {\psi(\tau)}{\phi'(\tau)}\Big| + |\Psi_1(\tau_2)-\Psi_1(\tau_1)|\Big).
\ee
Choose $\tau_1$ and $\tau_2$ such that
\be
\Psi_1(\tau_2) - \Psi_1(\tau_0) = k^{-3/2} \Big(\frac 1 \nu + \frac 1 \xi\Big)
\ee
and likewise for $\tau_1$, which takes care of the last term. This choice is possible and is unique due to the fact that $\Psi_1(\tau) \to \pm \infty$ as $\tau \to k-$ and $0+$, respectively, and $\Psi_1$ is strictly increasing.

Therefore, for $\tau \in [0, \tau_1] \cup [\tau_2, k]$,
\be
|\Psi_1(\tau) - \Psi_1(\tau_0)| \geq k^{-3/2} \Big(\frac 1 \nu + \frac 1 \xi\Big).
\lb{condititie}
\ee
It is left to prove that, under this condition,
\be
\Big|\frac {\psi(\tau)}{\phi'(\tau)}\Big| \leq C k^{-3/2} \Big(\frac 1 \nu + \frac 1 \xi\Big).
\ee
Rewrite condition (\ref{condititie}) as
\be
\frac 1 {k^2} |\nu^{1/2}\tau^{1/2} - \xi^{1/2}(k-\tau)^{1/2}| \cdot |\nu^{-1/2}(k-\tau)^{-1/2} - \xi^{-1/2}\tau^{-1/2}| \geq C k^{-3/2} \Big(\frac 1 \nu + \frac 1 \xi\Big)
\ee
and further as
\be
\frac{|\nu \tau - \xi(k-\tau)|}{(k-\tau)^{1/2}\tau^{1/2}} \geq C k^{1/2}\Big (\frac 1 {\nu} + \frac 1 {\xi}\Big) \frac {(\nu^{1/2}\tau^{1/2} + \xi^{1/2}(k-\tau)^{1/2})\nu^{1/2}\xi^{1/2}}{|\nu^{1/2}(k-\tau)^{1/2} - \xi^{1/2}\tau^{1/2}|}.
\ee
Then it suffices to prove that
\be
\Big(\frac 1 {\nu} + \frac 1 {\xi}\Big) \frac {(\nu^{1/2}\tau^{1/2} + \xi^{1/2}(k-\tau)^{1/2})\nu^{1/2}\xi^{1/2}}{|\nu^{1/2}(k-\tau)^{1/2} - \xi^{1/2}\tau^{1/2}|} \geq C
\ee
or equivalently
\be
(\nu+\xi)(\nu \xi^{1/2} \tau^{1/2} + \xi \nu^{1/2} (k-\tau)^{1/2}) \geq C \nu\xi |\nu^{1/2}(k-\tau)^{1/2} - \xi^{1/2}\tau^{1/2}|
\ee
which is obvious. The same goes for $\tau_2$.

Next, assume that $[a, b] \subset [k, \infty)$. The proof goes along the same lines, based on (\ref{abso}) and on (\ref{part}), with the difference that the antiderivative is now written
\be
\Psi_2 = \frac 1 {2k^2}(\tau^{1/2}(\tau-k)^{-1/2} + \tau^{-1/2}(\tau-k)^{-1/2})
\ee
and that
\be
\frac {\psi(\tau)}{\phi'(\tau)} = \frac 1 {\nu\tau+\xi(\tau-k)} \frac{\tau^{1/2}(\tau-k)^{1/2}}{\nu \tau-\xi (\tau-k)}.
\ee

If $\nu>\xi$ the phase has no stationary points in this interval. We divide $[k, \infty)$ into the subintervals on which $\ds \Psi_2(\tau) \leq \xi^{-1} k^{-3/2}$ and the rest. On the former the integral is trivially bounded and the number of such intervals is bounded from above. On the latter (also the union of a bounded number of intervals) one has, following integration by parts,
\be
\bigg|\int_{\tau_1}^{\tau_2} e^{i\phi} \psi \dd \tau\bigg| \leq C \max_{\substack{\tau \in [k, \infty) \\ \Psi_2(\tau) \geq \xi^{-1} k^{-3/2}}} \Big|\frac {\psi(\tau)}{\phi'(\tau)}\Big|.
\ee
However, the condition $\Psi_2(\tau) \geq \xi^{-1} k^{-3/2}$ implies that
\be
\frac{\xi(2\tau-k)}{\tau^{1/2}(\tau-k)^{1/2}} \geq C k^{1/2}
\ee
and therefore
\be
\frac {\tau^{1/2} (\tau-k)^{1/2}} {\nu\tau+\xi(\tau-k)} \frac 1 {(\nu -\xi) \tau + \xi k} \leq C k^{-1/2} \frac 1 {\xi k}.
\ee

A stationary point occurs if $\nu<\xi$, namely $\ds \tau_0 = \frac{\xi k} {\xi-\nu}$, which becomes infinite if $\xi=\nu$. Surround it with an interval $[\tau_1, \tau_2]$ on which we integrate the absolute value, otherwise integrate by parts. Overall, the integral is bounded by
\be
\bigg|\int_{a}^{b} e^{i\phi} \psi \dd \tau\bigg| \leq C (|\Psi_2(\tau_2) - \Psi_2(\tau_1)| + \max_{\tau \in [k, \tau_1] \cup [\tau_2, \infty)} \Big|\frac {\psi(\tau)}{\phi'(\tau)}\Big|).
\ee
Choose $\tau_1$ such that
\be
|\Psi_2(\tau_1) - \Psi_2(\tau_0)| = k^{-3/2} \frac 1 \nu
\ee
and likewise for $\tau_2$ (or $\tau_2 = \infty$ if there is no such value). This takes care of the integral on $[\tau_1, \tau_2]$. As for the remaining portion, note that for any $\tau \in [k, \tau_1] \cup [\tau_2, \infty)$ one has
\be
|\Psi_2(\tau_1) - \Psi_2(\tau_0)| \geq k^{-3/2} \frac 1 \nu
\ee
and therefore
\be
\frac 1 {k^2} |\nu^{1/2}\tau^{1/2} - \xi^{1/2}(\tau-k)^{1/2}| \cdot |(\tau-k)^{-1/2} \nu^{-1/2} - \tau^{-1/2} \xi^{-1/2}| \geq C k^{-3/2} \frac 1 \nu.
\ee
Equivalently, under the assumption $\nu <\xi$,
\be
\frac{|\nu\tau - \xi(\tau-k)|} {(\tau-k)^{1/2} \tau^{1/2}} \geq C\frac {k^{1/2}}{\nu} \frac{(\tau^{1/2}\nu^{1/2} + (\tau-k)^{1/2}\xi^{1/2})\xi^{1/2}\nu^{1/2}}{\tau^{1/2} \xi^{1/2} - (\tau-k)^{1/2} \nu^{1/2}}.
\ee
Note that
\be
\Big|\frac {\psi(\tau)}{\phi'(\tau)}\Big| \leq \frac 1 {\nu k} \frac {(\tau-k)^{1/2} \tau^{1/2}}{|\nu\tau - \xi(\tau-k)|}.
\ee
It is left to prove that
\be
\frac {\nu(\tau^{1/2} \xi^{1/2} - (\tau-k)^{1/2} \nu^{1/2})}{(\tau^{1/2}\nu^{1/2} + (\tau-k)^{1/2}\xi^{1/2})\xi^{1/2}\nu^{1/2}} \leq C.
\ee
However, the last statement is equivalent to
\be
\tau^{1/2} \xi^{1/2} \nu - (\tau-k)^{1/2} \nu^{3/2} \leq C (\tau^{1/2} \xi^{1/2} \nu + (\tau-k)^{1/2}\xi \nu^{1/2})
\ee
which is again obvious.

The third case $[a, b] \subset (-\infty, 0]$ is identical to the second case $[a, b] \subset [k, \infty)$.

Cutting the interval $[a, b]$ into at most three pieces according to this partition, one obtains the conclusion (\ref{hilbert}).
\end{proof}

The lemma relates to (\ref{HP}) in the following manner: when applying it to (\ref{HP}), we take $F_{s,a}(\lambda) =  \frac 1 2 f((\lambda+i0)^2+\mu)\pm f((\lambda-i0)^2+\mu)$ and likewise for $G_a$. This leads to a bound for the conditionally converging integral $T$ (\ref{expr_gen_T}).

After this general discussion of (\ref{expr_gen_T}), we return to the concrete Born sum expansion (\ref{unu}-\ref{patru}), expanded again, as previously stated, by means of the resolvent identity. Now we identify $F_s$, $F_a$, and $G_a$ for this case. From the Born sum expansion (\ref{unu}-\ref{patru}) and from the known expression for the kernel of the free resolvent
\be
R_0(\lambda^2+\mu \pm i0)(x, y) = \frac 1 {4\pi|x-y|} \left(\ba{cc} \exp(\pm |x-y| \lambda)& 0 \\ 0 & \exp(-|x-y|\sqrt{\lambda^2+2\mu})\ea\right)
\ee
(where $\sqrt{}$ is given by the main branch of the logarithm) we get the following types of factors:
\be
\ba{ll}
F_s(\lambda) = & \cos(A\lambda) \exp(-B\sqrt{\lambda^2+2\mu})\ (a) \\
& \text {or } ((R_0(\lambda^2+\mu)V)^m R_{V}(\lambda^2) (V R_0(\lambda^2+\mu))^m)_s\ (b) \\
& \text{or } ((R_0(\lambda^2+\mu)V)^m R_{V}(\lambda^2+\mu) (V R_0(\lambda^2+\mu))^{m-1})_s\ (c),\\
F_a(\lambda) = & \sin(A\lambda) \exp(-B\sqrt{\lambda^2+2\mu})\ (d) \\
& \text {or } ((R_0(\lambda^2+\mu)V)^m R_{V}(\lambda^2+\mu) (V R_0(\lambda^2+\mu))^m)_a\ (e) \\
& \text{or } ((R_0(\lambda^2+\mu)V)^m R_{V}(\lambda^2+\mu) (V R_0(\lambda^2+\mu))^{m-1})_a\ (f),\\
G_a(\lambda) = & \sin(C\lambda) \exp(-D\sqrt{\lambda^2+2\mu})\ (g) \\
& \text {or } ((R_0(\lambda^2+\mu)V^*)^m R_{V^*}(\lambda^2+\mu) (V^* R_0(\lambda^2+\mu))^m)_a\ (h) \\
& \text{or } ((R_0(\lambda^2+\mu)V^*)^{m-1} R_{V^*}(\lambda^2+\mu) (V R_0(\lambda^2+\mu))^m)_a\ (i),
\ea\lb{enum}
\ee
where we again used the notation $f_s(\lambda) = \frac 1 2(f(\lambda+i0) + f(\lambda-i0))$ and $f_a(\lambda) = \frac 1 2(f(\lambda+i0) - f(\lambda-i0))$.

Factors of the form (a), (d), and (g) stem from the general terms in the Born sum expansion, as in (\ref{unu}), while the others represent the contribution of the remainders and mixed terms (\ref{doi}-\ref{patru}).

We evaluate the first type of factors, coming from (\ref{unu}), in view of applying Lemma \ref{baza}. Note that, uniformly in $A$ and $B$,
\begin{multline}
\|\big(\sin(A\lambda) \exp(-B\sqrt{\lambda^2+2\mu})\big)^{\wedge}\|_1 \leq \\
\leq \|(\sin(A\lambda))^{\wedge}\|_{\mc M} \|(\exp(-B\sqrt{\lambda^2+2\mu}))^{\wedge}\|_1 \leq C,
\lb{est_sin1}
\end{multline}
where $\mc M$ is the space of finite measures. This is true because $(\sin(A\lambda))^{\wedge}$ is the sum of two point measures of mass $1/2i$, while
\begin{multline}
\|(\exp(-B\sqrt{\lambda^2+2\mu}))^{\wedge}\|_1 = \Big\|\int_{-\infty}^{\infty} e^{-B\sqrt{\lambda^2+2\mu}} e^{-i\tau\lambda} \dd \lambda\Big\|_{L^1_\tau}\\
= \Big\|\int_{-\infty}^{\infty} e^{-\sqrt{\eta^2+2\mu B^2}} e^{-i\tau\eta} \dd \eta Big\|_{L^1_\tau} \leq C \|e^{-\sqrt{\eta^2+2\mu B^2}}\|_{H^1_{\eta}} \leq C \|e^{-|\eta|}\|_{H^1} < \infty.
\end{multline}
Similarly,
\begin{multline}
\Big\|\Big(\frac d {\dd \lambda}\big(\sin(A\lambda) \exp(-B\sqrt{\lambda^2+2\mu})\big)\Big)^{\wedge}\Big\|_1 \leq \\
\begin{aligned}
\leq &\|\big(A\cos(A\lambda) \exp(-B\sqrt{\lambda^2+2\mu})\big)^{\wedge}\|_1 + \\
& + \Big\|\Big(\frac{\sin(A\lambda)}{\lambda} \frac {B\lambda^2}{\sqrt{\lambda^2+2\mu}} \exp(-B\sqrt{\lambda^2+2\mu})\Big)^{\wedge}\Big\|_1 \leq C A,
\end{aligned}
\lb{est_sin2}
\end{multline}
because
\be
\Big\|\Big(\frac{\sin(A\lambda)}{\lambda}\Big)^{\wedge}\Big\|_1 = \frac 1 2 \|\chi_{[-A, A]}\|_1 \leq CA
\ee
and
\be
\begin{aligned}
\Big\|\big(\frac {B\lambda^2}{\sqrt{\lambda^2+2\mu}} \exp(-B\sqrt{\lambda^2+2\mu})\big)^{\wedge}\Big\|_1 = &\Big\|\int_{-\infty}^{\infty} \frac {B\lambda^2}{\sqrt{\lambda^2+2\mu}} e^{-B\sqrt{\lambda^2+2\mu}} e^{-i\tau\lambda} \dd \lambda\Big\|_{L^1_\tau}\\
= &\Big\|\int_{-\infty}^{\infty} \frac {\eta^2}{\sqrt{\eta^2+2\mu B^2}} e^{-\sqrt{\eta^2+2\mu B^2}} e^{-i\tau\eta} \dd \eta\Big\|_{L^1_\tau} \\
\leq &C \Big\|\frac {\eta^2}{\sqrt{\eta^2+2\mu B^2}} e^{-\sqrt{\eta^2+2\mu B^2}}\Big\|_{W^{1, 2}_{\eta}}\\
\leq &\||\eta| e^{-|\eta|}\|_{W^{1, 2}} < \infty.
\end{aligned}\ee
With sine replaced by cosine, the analogous estimate is
\be
\Big\|\Big(\cos(A\lambda) \lambda \frac d {\dd \lambda} \exp(-B\sqrt{\lambda^2+2\mu})\Big)^{\wedge}\Big\|_1 \leq C.
\lb{est_cos}
\ee
This takes care of factors of the form (a), (d), or (g).

We proceed to obtain a bound for the remaining factors. Intuitively, since they represent the remainder in the Born sum expansion, the bound should be less sharp. Let
\be
\ba{ll}
G_x(\lambda^2)(x_1) &= \bpm e^{-i\lambda|x|} & 0 \\ 0 & 1 \epm R_0(\lambda^2+\mu)(x, x_1) \\
&= \ds\frac 1 {4\pi|x-x_1|} \bpm \exp(i(|x-x_1|-|x|)\lambda)& 0 \\ 0 & \exp(-|x-x_1|\sqrt{\lambda^2+2\mu})\epm
\ea
\ee
and
\be
a_{x, y}(\lambda^2) = \left\langle V R_V(\lambda^2+\mu) V (R_0(\lambda^2+\mu) V)^m G_y(\lambda^2), (R_0^*(\lambda^2+\mu) V^*)^m G^*_x(\lambda^2)\right\rangle.
\ee
The following estimates hold for $G$:
\be
\sup_x \|\partial^j_{\lambda} G_x(\lambda^2)\|_{\langle x_1 \rangle^{3/2+j+\epsilon} L^2_{x_1}} \leq C \langle x \rangle^{-1},\ \sup_x \|\partial^j_{\lambda} G_x(\lambda^2)\|_{\langle x_1 \rangle^{1/2+j+\epsilon} L^2_{x_1}} \leq C,
\ee
implying that
\be
\begin{aligned}
\left|\frac{d^j}{\dd \lambda^j}a_{x, y}(\lambda^2)\right| \leq& C \langle \lambda \rangle^{-2-\epsilon}\langle x\rangle^{-1}\langle y \rangle^{-1} \text{ for } j=0, 1 \\
\left|\frac{d^j}{\dd \lambda^j}a_{x, y}(\lambda^2)\right| \leq& C \langle\lambda\rangle^{-2-\epsilon} \text{ for } j=2,
\end{aligned}
\ee
provided that $m$ is sufficiently large and $V$ has sufficient decay, $|V| \leq C \langle x \rangle^{-7/2-\epsilon}$.

We used the limiting absorbtion principle Lemma \ref{limabs} to bound this quantity. Incrementing $m$ by $1$ increases the decay by $\langle \lambda \rangle^{-1/2}$. The decay condition on $V$ arises as follows: if, for example, two derivatives fall on the $R_V$ factor, then $V$ has to compensate for $3+\epsilon$ powers of x from its output and for another $1/2+\epsilon$ power coming from the pairing with $R_0$.

Then the factors of the form (b) and (e) can be written as
\be
K_{Vs}(\lambda) = \frac 1 2 (K_V((\lambda+i0)^2) + K_V((\lambda-i0)^2),
\ee
respectively
\be
K_{Va}(\lambda) = \frac 1 2 (K_V((\lambda+i0)^2) - K_V((\lambda-i0)^2),
\ee
where
\be\begin{aligned}
K_V(\lambda^2)(x, y) &= ((R_0(\lambda^2+\mu)V)^m R_{V}(\lambda^2+\mu) (V R_0(\lambda^2+\mu))^m)(x, y) \\
&= \bpm e^{i\lambda|x|} & 0 \\ 0 & 1 \epm a_{x, y}(\lambda^2) \bpm e^{i\lambda|y|} & 0 \\ 0 & 1 \epm.
\end{aligned}\ee

Clearly
\be
\|\widehat{K_{Va}}\|_1 = \|(a_{x, y}(\lambda^2))^{\wedge}\|_1 \leq \|a_{x, y}(\lambda^2)\|_{W^{1, 2}(\set R)} \leq C \langle x\rangle^{-1}\langle y \rangle^{-1}
\lb{est_k1}
\ee
and, taking into account the fact that $K_{Va}(\lambda)$ is of the form $F_a(\sqrt{\lambda})$ for antisymmetric $F_a$,
\be
\|\widehat{K_{Va}'}\|_1 \leq \|\frac d {\dd \lambda} a_{x, y}(\lambda^2)\|_{W^{1, 2}(\set R)} + (|x|+|y|)\|(a_{x, y}(\lambda^2))^{\wedge}\|_1 \leq C.
\lb{est_k2}
\ee
Moreover,
\be
\Big\|\Big(\frac {K_{Va}}{\lambda}\Big)^{\wedge}\Big\|_1 \leq C \|\xi \widehat{K_{Va}}(\xi)\|_1 \leq C \|\widehat{K_{Va}'}\|_1 \leq C.
\lb{est_k3}
\ee
These decay estimates carry on to factors of the form (b) or (e) and, with minimal modifications, to factors of the form (c), (f), (h), and (i), provided that $m$ is sufficiently large.

After showing that all types of terms arising from (\ref{enum}) satisfy the prerequisites for applying Lemma \ref{baza}, we apply it to each of them, in turn.

We first deal with the $4m^2$ terms in (\ref{unu}). As mentioned previously, each splits into two parts after using the resolvent identity. Let us introduce the following notations:
\be
K(x_0, \ldots, x_{k+\ell-1}) = \frac {\prod_{j=1}^{k+\ell-2} {V_j(x_j)}} {\prod_{j=0}^{k+\ell-2}|x_{j+1}-x_j|},
\ee
where $V_j$ are entries of $V$ (that is, $\pm U$ or $\pm W$);
\be
A = \sum_{j \in J} |x_{j+1}-x_j|,\ B = \sum_{j \in J^c}|x_{j+1}-x_j|,\ \mc C = \sum_{j \in I} |x_{j+1}-x_j|,\ D = \sum_{j \in I^c}|x_{j+1}-x_j|,
\ee
where $J \subset \{0, \ldots, \ell-1\}$ and $J^c = \{0, \ldots, \ell-1\} \setminus J$; $I \subset \{\ell, \ldots, \ell+k-2\}$ and $I^c = \{\ell, \ldots, \ell+k-2\} \setminus I$. Also, let
\be\begin{aligned}
F_s(\lambda) = &\cos(A\lambda) \exp(-B\sqrt{\lambda^2+2\mu}),\\
F_a(\lambda) = &\sin(A\lambda) \exp(-B\sqrt{\lambda^2+2\mu}),\\
G_a(\lambda) = &\sin(\mc C\lambda) \exp(-D\sqrt{\lambda^2+2\mu}).
\lb{fafs}
\end{aligned}\ee
After applying the resolvent identity and performing the matrix multiplication, we find that each term of (\ref{unu}) is a sum, for all possible choices of $J$ and $I$, of terms with kernel of the form
\begin{multline}
K_{k-1\ \ell}(x_0, x_{k+\ell-1}) = \int_{\set R^{3(k+\ell-2)}} K(x_0, \ldots, x_{k+\ell-1}) \int_0^{\infty} \bigg(e^{i(t-s)\lambda} F_s(\sqrt \lambda) G_a(\sqrt\lambda) + \\
+ e^{it\lambda} F_a(\sqrt\lambda) H_{\eta}(\chi_{[0, \infty]} e^{-is \eta} G_a(\sqrt\eta))(\lambda) \bigg) \dd \lambda \dd x_1 \ldots \dd x_{k+\ell-2}.
\end{multline}
Then, in view of Lemma \ref{baza} and the bounds (\ref{est_sin1}), (\ref{est_sin2}), and (\ref{est_cos}), we have that
\begin{multline}
\Big| \int_0^{\infty} e^{i(t-s)\lambda} F_s(\sqrt \lambda) G_a(\sqrt\lambda) + e^{it\lambda} F_a(\sqrt\lambda) H(\chi_{[0, \infty]} e^{-is \eta} G_a(\sqrt\eta))(\lambda) \dd \lambda \Big| \leq \\ \leq C |t-s|^{-3/2} (A+\mc C) \leq C |t-s|^{-3/2} \sum_{j=0}^{k+\ell-2} |x_{j+1}-x_j|.
\end{multline}
Combining this with Lemma 2.5, p. 12 of \cite{rod}, which states that
\be
\Big|\int_{\set R^{3(k+\ell-2)}} K(x_0, \ldots, x_{k+\ell-1}) \sum_{j=0}^{k+\ell-2} |x_{j+1}-x_j| \dd x_1 \ldots \dd x_{k+\ell-2} \Big| \leq (k+\ell)\|V\|_{\mathcal K}^{k+\ell-1},
\ee
where $\ds \|V\|_{\mathcal K} = \sup_{x \in \set R^3} \int \frac {|V(y)|}{|x-y|} dy$ is the Kato norm, we obtain that
\be
\|K_{k-1\ \ell}\|_{1 \to \infty} \leq \sup_{x_0, x_{k+\ell-1}} K_{k-1\ \ell}(x_0, x_{k+\ell-1}) \leq C |t-s|^{-3/2}.
\ee

The same method can be applied to the remaining terms in (\ref{doi}), (\ref{trei}), and (\ref{patru}). Since (\ref{doi}) and (\ref{trei}) are similar, we look at a typical term of one of these two sums, consisting of the product between some term of the Born sum expansion, on one hand, and the remainder, on the other. The kernel of such a term is of the form
\begin{multline}
K_{k-1\ V^*}(x_0, x_{k-1+2m}) = \int_{\set R^{3(k-2+2m)}} K(x_0, \ldots, x_{k-1}) V(x_{k-1}) \cdot \\
\begin{aligned}&\cdot \int_0^{\infty} e^{i(t-s)\lambda} F_s(\sqrt\lambda) e^{it\lambda} K_{V^*a}(x_{k-1}, \ldots, x_{k-1+2m}, \sqrt\lambda) + \\
&+ e^{it\lambda} F_a(\sqrt\lambda) H_{\eta} (\chi_{[0, \infty]} e^{-is \eta} K_{V^*a}(x_{k-1}, \ldots, x_{k-1+2m}, \sqrt \eta))(\lambda) \dd \lambda \dd x_1 \ldots x_{k-2+2m},
\end{aligned}
\lb{termenmixt}
\end{multline}
where $F_a$, $F_s$ are as in (\ref{fafs}).

The kernel $K_{V^*a}$ involves $R_{V^*}$ and is not given by an explicit formula as $F_a$ is. However, we still have estimates (\ref{est_sin1}), (\ref{est_sin2}), (\ref{est_cos}), as well as (\ref{est_k1}), (\ref{est_k2}), (\ref{est_k3}), based on the limiting absorbtion principle.

We could treat treat the pairing $F_s K_{V^*a}$ using only inequalities of the form (\ref{est_k1}), (\ref{est_k2}), because $K(x_0, \ldots, x_{k-1}) V(x_{k-1}) F_s(\sqrt \lambda) K_{V^*a}(x_{k-1}, \ldots, x_{k-1+2m}, \sqrt(\lambda))$ is the sum of two antisymmetric kernels. Otherwise, keeping within the previous framework, note that we can take a factor of $\lambda$ from the antisymmetric part over to the symmetric part whenever that is needed.

By Lemma \ref{baza},
\begin{multline}
\Big|\int_0^{\infty} e^{i(t-s)\lambda} F_s(\sqrt\lambda) e^{it\lambda} K_{V^*a}(x_{k-1}, \ldots, x_{k-1+2m}, \sqrt\lambda) + \\
\begin{aligned}
&+ e^{it\lambda} F_a(\sqrt\lambda) H_{\eta} (\chi_{[0, \infty]} e^{-is \eta} K_{V^*a}(x_{k-1}, x_{k-1+2m}, \sqrt \eta))(\lambda)\\
&\dd \lambda \dd x_k \ldots \dd x_{k-2+2m}\Big| \leq \\
&\leq C|t-s|^{-3/2}\big (\sum_{j=0}^{k-2} |x_{j+1}-x_j| + 1\big).
\end{aligned}
\lb{cuomisiuni}
\end{multline}
This estimate, to which we now add back those factors of (\ref{termenmixt}) that we omitted for convenience in (\ref{cuomisiuni}), results in
\begin{multline}
|K_{k-1\ V^*}(x_0, x_{k-1+2m})| \leq \\
\begin{aligned}
\leq &C|t-s|^{-3/2} \int K(x_0, \ldots, x_{k-1}) |V(x_{k-1})| (\sum_{j=0}^{k-2} |x_{j+1}-x_j| + 1) \dd x_1 \ldots \dd x_{k-1}\\
\leq &C|t-s|^{-3/2} \|V\|_1 \sup_{x_{k-1}} \bigg|\int K(x_0, \ldots, x_{k}) \sum_{j=0}^{k-2} |x_{j+1}-x_j| \dd x_1 \ldots \dd x_{k-2}\bigg|  + \\
& + C|t-s|^{-3/2}\bigg|\int K(x_0, \ldots, x_{k-1}) V(x_{k-1}) \dd x_1 \ldots \dd x_{k-1}\bigg| \leq C|t-s|^{-3/2}.
\end{aligned}
\lb{imp}
\end{multline}
Thus we have proved that $\|K_{k-1\ V^*}\|_{1 \to \infty} \leq C|t-s|^{-3/2}$.

Finally, the last term appearing in (\ref{patru}), consisting of the product of the remainders in the Born sum expansion, yields to the same approach. The final step of the computation, instead of (\ref{imp}), is
\be\ba{l}
|K_{V\ V^*} (x_0, x_{4m-1})| \leq\\
\ds \leq C|t-s|^{-3/2} \int \langle x_0\rangle^{-1}\langle x_{2m}\rangle^{-1} V(x_{2m})+ V(x_{2m}) \langle x_{2m}\rangle^{-1}\langle x_{4m-1} \rangle^{-1} \dd x_{2m}\\
\leq C|t-s|^{-3/2} \|V\|_1 \leq C |t-s|^{-3/2}.
\ea\ee

This completes the proof of the fact that $\|e^{it\mc H} P_+ e^{-is\mc H^*} P_+^*\|_{1 \to \infty} \leq C|t-s|^{-3/2}$.

The other $3$ combinations are entirely analogous. Indeed, the $P_- P_-$ term can be treated by the same means. As for the mixed terms, a very similar approach works and we present the proof in brief. Again, we expand both factors into finite Born series and we obtain a sum with $(2m+1)^2$ terms.

Each term has a kernel that can be written, by means of the resolvent identity, as a sum of $2$ parts of the form
\be
T=\int_{\Gamma^+_{\epsilon_1}} \int_{\Gamma^-_{\epsilon_2}} \frac{e^{i(t\lambda-s\eta)}}{\lambda-\eta} f(\lambda) g(\eta) \dd \eta \dd \lambda.
\ee
After making $\epsilon_1=\epsilon_2=0$, we get, by analogy to (\ref{HP}),
\be\ba{ll}
T&\ds=\int_0^{\infty} e^{it\lambda} f_a(\lambda+\mu) H(\chi_{(-\infty, 0]} e^{-is \eta} g_a(\eta-\mu))(\lambda) \dd \lambda,
\ea\ee
where $H$ is the Hilbert transform. The other term, involving the Dirichlet kernel, cancels, because now the contours surround disjoint regoins.

However, now we note upon inspection that the terms stemming from $\chi_{(-\infty, 0]} g_a(\eta-\mu)$ have the same form as those we have already enumerated in (\ref{enum}). From here the proof proceeds in the same manner as in the previous case. The bound that we eventually obtain for these terms is even better, $C|t+s|^{-3/2}$ instead of $C|t-s|^{-3/2}$, because the phases add instead of cancelling.
\end{proof}

\subsection{Other linear estimates}

\begin{proof}[Proof of Corollary \ref{Strichartz}]
The dispersive estimate (\ref{I2}) holds for $\mc H^*$ as well as $\mc H$, since they are conjugated by $\sigma_3$. This, together with the $L^2$ bound (\ref{l2}), implies the Keel-Tao endpoint Strichartz estimate.
\end{proof}

Next, by interpolating between the endpoint Strichartz estimate for $L^2$ initial data and the $L^1 \to L^{\infty}$ decay estimates, we achieve an improved decay of the solution, in norm, for $L^p$ initial data, $1\leq p\leq 2$.

In the sequel, denote $\bpm a \\ b \epm^{\dg} = \bpm a & b \epm$, so that $\bpm a \\ b \epm \bpm c \\ d \epm^{\dg} = \bpm ac & ad \\ bc & bd \epm$.

\begin{lemma} For $1\leq q \leq 2$,
\be
\int_{T}^{\infty} \|e^{it\mc H} P_c U_1 (e^{it\mc H} P_c U_2)^{\dg}\|_{\frac{3q}{2(q-1)}} \dd t \leq CT^{2- \frac 4 q} \|U_1\|_{q} \|U_2\|_{q}.
\lb{biest1}
\ee
Likewise, for $1 \leq q <4/3$, $\beta <2/q-1$,
\be
\int_{T}^{\infty} \langle t \rangle^{2\beta} \|e^{it\mc H} P_c U_1 (e^{it\mc H} P_c U_2)^{\dag}\|_{3+\frac{3q}{2(q-1)}} \dd t \leq C\langle T \rangle^{2-\frac 4 q + \beta} \|U_1\|_{q\cap 2} \|U_2\|_{q\cap 2}.
\lb{biest2}
\ee
\end{lemma}
\begin{proof}
We obtain the first result by complex bilinear interpolation (see \cite{bergh}, p. 96, Theorem 4.4.1). We use it in the following form:
\begin{theorem}
For $i=0, 1$, let $\mc B: A_i \oplus B_i \to C_i$ be a bilinear mapping such that
\be
\|\mc B(a, b)\|_{C_i} \leq M_i \|a\|_{A_i} \|b\|_{B_i}.
\ee
Then, for each $0 \leq \theta \leq 1$, $\mc B$ can be extended uniquely to a bilinear mapping from $[A_0, A_1]_{[\theta]} \oplus [B_0, B_1]_{[\theta]}$ to $[C_0, C_1]_{[\theta]}$ with norm at most $M_0^{\theta} M_1^{1-\theta}$.
\end{theorem}
The first estimate (\ref{biest1}) then follows directly if we take
\be
\mc B(U_1, U_2) = e^{it\mc H} P_c U_1 (e^{it\mc H} P_c U_2)^{\dag}
\ee
and interpolate between
\be
\int_{T}^{\infty} \|e^{it\mc H} P_c U_1 (e^{it\mc H} P_c U_2)^{\dag}\|_{\infty} \dd t \leq CT^{-2} \|U_1\|_{1} \|U_2\|_{1}
\ee
and
\be
\int_{T}^{\infty} \|e^{it\mc H} P_c U_1 (e^{it\mc H} P_c U_2)^{\dag}\|_3 \dd t \leq C\|U_1\|_2 \|U_2\|_2.
\ee
The second statement (\ref{biest2}) follows from the first, once we write
\begin{multline}
\int_{T}^{\infty} \langle t \rangle^{2\beta} \|e^{it\mc H} P_c U_1 (e^{it\mc H} P_c U_2)^{\dg}\|_{3+\frac{3q}{2(q-1)}} \dd t \leq \\
\begin{aligned}
\leq &\langle T \rangle^{2\beta} \int_{T}^{\infty} \|e^{it\mc H} P_c U_1 (e^{it\mc H} P_c U_2)^{\dg}\|_{3+\frac{3q}{2(q-1)}} \dd t + \\
& + C \int_{T}^{\infty} \langle t \rangle^{2\beta-1} \int_t^{\infty} \|e^{is\mc H} P_c U_1 (e^{is\mc H} P_c U_2)^{\dg}\|_{3+\frac{3q}{2(q-1)}} \dd s \dd t.
\end{aligned}
\end{multline}
\end{proof}

From this we derive an inequality concerning the solution of the inhomogenous problem.
\begin{lemma}
Consider the equation
\be
i\partial_t U + \mc H P_c U = RHS(t),\ U(0) \text{ given}.
\ee
Then, for $q<4/3$, $\beta < 2/q-1$,
\be
\int_{0}^{\infty} \langle t \rangle^{2\beta} \|P_c U\|_{6+\infty}^2 \dd t \leq C (\|U(0)\|_{q \cap 2} + \|RHS\|_{\langle t \rangle^{-{\beta}} L^{2}_t L^{1\cap 6/5}_x}) ^2.
\ee
\lb{biest}
\end{lemma}
Note that for the exponent of interest, $\beta = 1/2$, it is possible to replace $\infty$ by $12+\epsilon$ and $1$ by $6/5-\epsilon$, if $q$ is sufficiently close to $4/3$.
\begin{proof}
Firstly, we examine the source terms. Setting $U_1=U_2$ in (\ref{biest2}), one has
\be
\int_{T}^{\infty} \langle t \rangle^{2\beta} \|e^{it\mc H} P_c U(0)\|^2_{6+\frac{3q}{q-1}} \dd t \leq C\langle T \rangle^{2-4/q + 2\beta} \|U(0)\|_{q \cap 2}^2
\ee
and note that $L^{6+\frac{3q}{q-1}} \subset L^{6+\infty}$.

We then evaluate the inhomogenous terms. They have $L^2$ in time decay, as will follow from (we henceforth denote $\min(a, b) = a \wedge b$)
\begin{multline}
\int_0^{\infty} \langle t \rangle^{2\beta} \Big\|\int_0^t e^{i(t-s) \mc H} P_c RHS (s) \dd s\Big\|_{6+\infty}^2 \dd t \leq \\
\begin{aligned}
\leq &C\Big(\int_0^{\infty} \langle t \rangle^{2\beta} \Big\|\int_{t/2 \wedge (t-1)}^t e^{i(t-s) \mc H} P_c RHS (s) \dd s\Big\|_{6}^2 \dd t + \\
&+ \int_0^{\infty} \langle t \rangle^{2\beta} \Big\|\int_0^{t/2 \wedge (t-1)} e^{i(t-s) \mc H} P_c RHS (s) \dd s\Big\|_{\infty}^2 \dd t\Big).
\end{aligned}
\end{multline}
Now we examine the two expressions separately. Concerning the first, note that what one needs to prove is equivalent to
\begin{multline}
\Big|\iint_{t/2 \wedge (t-1) \leq s \leq t} \langle t \rangle^{-2\beta} \langle s \rangle^{2\beta} \langle e^{-is \mc H} P_c  F_1 (s), e^{-it\mc H^*} F_2(t) \rangle \dd s \dd t \Big| \leq \\
\leq C \|F_1\|_{L^2_s L^{6/5}_x} \|F_2\|_{L^2_t L^{6/5}_x}.
\end{multline}
However, observe that this follows by the same means as the usual Keel-Tao endpoint Strichartz estimate, by a dyadic partition, because the extra $\langle t \rangle^{-2\beta} \langle s \rangle^{2\beta}$ factor is bounded.

The second term can be handled as follows:
\begin{multline}
\int_0^{\infty} \langle t \rangle^{2\beta} \|\int_0^{t/2 \wedge(t-1)} e^{i(t-s) \mc H} P_c RHS (s) \dd s\|_{\infty}^2 \dd t \leq \\
\begin{aligned}
\leq &\int_0^{\infty} t^{2\beta-3} \Big(\int_0^{t/2 \wedge(t-1)} \langle s \rangle^{2\beta} \|RHS(s)\|_1^2 \dd s\Big) \Big(\int_0^{t/2 \wedge(t-1)} \langle s \rangle^{-2\beta} \dd s\Big)\\
\leq &\int_0^{\infty} \langle s \rangle^{2\beta} \Big(\int_{s+1}^{\infty} t^{2\beta-3} \dd t\Big) \|RHS(s)\|_1^2 \dd s \leq \int_0^{\infty} \langle s \rangle^{2\beta} \|RHS(s)\|_1^2 \dd s,
\end{aligned}
\end{multline}
provided that $2\beta < 2$.

It follows that for $\beta < 1$
\be
\|\langle t \rangle^{2\beta} \int_0^t e^{i(t-s) \mc H} P_c RHS (s) \dd s\|_{L^2_t L^{6+\infty}_x} \leq \|\langle s \rangle^{2\beta} RHS (s)\|_{L^2_s L^{1 \cap 6/5}_x}.
\ee
This suffices to bound the product of the two inhomogenous terms.

The product of a source term and an inhomogenous term can be handled in the same manner,
\begin{multline}
\int_T^{\infty} \langle t \rangle^{2\beta} \|e^{it\mc H} P_c U(0)\|_{6+\infty} \Big\|\int_0^t e^{i(t-s)\mc H} P_c {RHS}(s) \dd s\Big\|_{6+\infty} \dd t \leq \\
\begin{aligned}
\leq &\Big(\int_T^{\infty} \langle t \rangle^{2\beta} \|e^{it\mc H} P_c U_1(0)\|_{\frac {3p}{2(p-1)}}^2 \dd t \Big)^{1/2}\\
&\Big(\int_T^{\infty} \langle t \rangle^{2\beta} \Big\|\int_0^t e^{i(t-s)\mc H} P_c {RHS}_2(s) \dd s\Big\|_{6+\infty}^2 \dd t\Big)^{1/2} \\
\leq &\|U(0)\|_{q\cap 2} \|RHS\|_{\langle t \rangle^{-{\beta}} L^{2}_t L^{1\cap 6/5}_x}.
\end{aligned}
\end{multline}
\end{proof}


\begin{thebibliography}{25}
\bibitem[Agm]{agmon} S. Agmon, \emph{Spectral properties of Schršdinger operators and scattering theory}, Ann.\ Scuola Norm.\ Sup.\ Pisa Cl.\ Sci.\ (4)  2 (1975), no.\ 2, pp. 151--218.

\bibitem[BatJon]{bates} P. W. Bates, C.\ K.\ R.\ T.\ Jones, \emph{Invariant manifolds for semilinear partial differential equations}, Dynamics Reported 2 (1989), pp. 1-38.

\bibitem[BerCaz]{bercaz} H.\ Berestycki, T.\ Cazenave, \emph{Instabilit\'e des \'etats stationnaires dans les \'equations de Schr\"odinger et de Klein-Gordon non lin\'eaires}, C.\ R.\ Acad.\ Sci.\ Paris S\'er.\ I Math. 293 (1981), no.\ 9, pp. 489-492.

\bibitem[BerLio]{bere} H.\ Berestycki, P.\ L.\ Lions, \emph{Nonlinear scalar field equations. I. Existence of a ground state}, Arch.\ Rational Mech.\ Anal.\ 82 (1983), no.\ 4, pp. 313-345.

\bibitem[BerL\"of]{bergh} J.\ Bergh, J.\ L\"ofstr\"om, \emph{Interpolation Spaces. An Introduction}, Springer-Verlag, 1976.

\bibitem[BusPer1]{buslaev1} V.\ S.\ Buslaev, G.\ S.\ Perelman, \emph{Scattering for the nonlinear Schr\"odinger equation: states that are close to a soliton} (Russian), Algebra i Analiz 4 (1992), no. 6, 63--102; translation in St.\ Petersburg Math.\ J.\ 4 (1993), no. 6, pp. 1111-1142.

\bibitem[BusPer2]{buslaev2} V.\ S.\ Buslaev, G.\ S.\ Perelman, \emph{On the stability of solitary waves for nonlinear Schr\"odinger equations}, Nonlinear evolution equations, 75--98, Amer.\ Math.\ Soc.\ Transl.\ Ser.\ 2, 164, Amer. Math. Soc., Providence, RI, 1995.

\bibitem[BusPer3]{buslaev3} V.\ S.\ Buslaev, G.\ S.\ Perelman, \emph{Nonlinear scattering: states that are close to a soliton} (Russian), Zap. Nauchn.\ Sem.\ S.-Peterburg.\ Otdel.\ Mat.\ Inst.\ Steklov, (POMI) 200 (1992), Kraev.\ Zadachi Mat.\ Fiz.\ Smezh.\ Voprosy Teor.\ Funktsii, 24, 38--50, 70, 187; translation in J.\ Math.\ Sci.\ 77 (1995), no. 3, 3161--3169.

\bibitem[Caz]{caz2} T.\ Cazenave, \emph{Semilinear Schr\"odinger equations}, Courant Lecture Notes in Mathematics, 10, New York University, Courant Institute of Mathematical Sciences, New York; AMS, Providence, RI, 2003.

\bibitem[CazLio]{cazenave} T.\ Cazenave and P.\ L.\ Lions, \emph{Orbital stability of standing waves for some nonlinear Schr\"odinger equations}, Commun.\ Math.\ Phys.\ 85, 549-561 (1982).

\bibitem[Cof]{coffman} C.\ V.\ Coffman, \emph{Uniqueness of positive solutions of $\Delta u - u + u^3 = 0$ and a variational characterization of other solutions}, Arch.\ Rat.\ Mech.\ Anal., 46 (1972), pp.\ 81-95.

\bibitem[Cuc]{cuc} S.\ Cuccagna, \emph{Stabilization of solutions to nonlinear Schr\"odinger equations}, Comm.\ Pure Appl.\ Math.\ 54 (2001), no.\ 9, pp.\ 1110-1145.

\bibitem[CucPelVou]{cuc2} S.\ Cuccagna, D.\ Pelinovsky, V.\ Vougalter, \emph{Spectra of positive and negative energies 
in the linearized NLS problem}, Comm.\ Pure Appl.\ Math.\ 58 (2005), no. 1, 1--29.

\bibitem[DemSch]{demanet} L.\ Demanet, W.\ Schlag, \emph{Numerical verification of a gap condition for a linearized NLS equation}, Nonlinearity 19 (2006), pp.\ 829-852.

\bibitem[Erdsch]{erdogan2} B.\ Erdogan, W.\ Schlag, \emph{Dispersive estimates for Schr\"odinger operators in the presence of a resonance and/or an eigenvalue at zero energy in dimension three: II}, to appear in Journal d'Analyse Mathematique.

\bibitem[GeJoLaSt]{ges} F.\ Gesztesy, C.\ K.\ R.\ T.\ Jones, Y.\ Latushkin \& M.\ Stanislavova, \emph{A spectral mapping theorem and invariant manifolds for nonlinear Schr\"odinger equations}, Indiana University Mathematics Journal Vol.\ 49, No.\ 1 (2000), pp.\ 221-243.

\bibitem[Gla]{glassey} R.\ T.\ Glassey, \emph{On the blowing-up of solutions to the Cauchy problem for the nonlinear Schr\"odinger equation}, J.\ Math.\ Phys.\ 18, (1977), pp.\ 1794-1797.

\bibitem[GrShSt1]{gril1} M.\ Grillakis, J.\ Shatah, W.\ Strauss, \emph{Stability theory of solitary waves in the presence of symmetry. I}, J. Funct. Anal. 74 (1987), no.\ 1, pp.\ 160-197.

\bibitem[GrShSt2]{gril2} M.\ Grillakis, J.\ Shatah, W.\ Strauss, \emph{Stability theory of solitary waves in the presence of symmetry. II}, J.\ Funct.\ Anal.\ 94 (1990), no.\ 1, pp.\ 308-348.

\bibitem[HunLee]{hund} D.\ Hundertmark, Y.-R.\ Lee, \emph{Exponential decay of eigenfunctions and generalized eigenfunctions of a non self-adjoint matrix Schr\"odinger operator related to NLS}, preprint 2006.

\bibitem[KeeTao]{tao} M.\ Keel, T.\ Tao, \emph{Endpoint Strichartz estimates}, Amer.\ Math.\ J.\ 120 (1998), pp.\ 955-980.

\bibitem[KenMer]{kenig} C.\ Kenig, F.\ Merle, \emph{Global well-posedness, scattering and blow-up for the energy-critical, focusing, non-linear Schr\"odinger equation in the radial case}, to appear, Inv.\ Math.

\bibitem[Kel]{kelley} P.\ L.\ Kelley, \emph{Self-focusing of optical beams}, Phys.\ Rev.\ Lett.\ 15, pp.\ 1005-1008 (1965).

\bibitem[KriSch1]{krieger} J.\ Krieger, W.\ Schlag, \emph{Stable manifolds for all monic supercritical NLS in one dimension}, Journal of the AMS, 
Volume 19, Number 4, October 2006, pp. 815-920.

\bibitem[KriSch2]{kri3} J.\ Krieger, W.\ Schlag, \emph{Non-generic blow-up solutions for the critical focusing NLS in 1-d}, preprint 2005.

\bibitem[KriSch3]{kri4} J.\ Krieger, W.\ Schlag, \emph{On the focusing critical semi-linear wave equation}, preprint 2005, to appear in Amer.\ J.\ Math.

\bibitem[Kwo]{kwong} M.\ K.\ Kwong, \emph{Uniqueness of positive solutions of $\Delta u - u + u^p = 0$ in $\set R^n$}, Arch.\ Rat.\ Mech.\ Anal.\ 65 (1989), pp.\ 243-266.

\bibitem[McLSer]{mcl} K.\ McLeod, J.\ Serrin, \emph{Nonlinear Schr\"odinger equation. Uniqueness of positive solutions of $\Delta u + f (u) = 0$ in $\set R^n$}, Arch.\ Rat.\ Mech.\ Anal.\ 99 (1987), pp.\ 115-145.

\bibitem[Mer]{merle} F.\ Merle, \emph{Determination of blow-up solutions with minimal mass for nonlinear Schrödinger equations with critical power}, Duke Math.\ J.\ 69 (1993), no.\ 2, pp.\ 427-454.

\bibitem[MerRap]{mer2} F.\ Merle, P.\ Raphael, \emph{On a sharp lower bound on the blow-up rate for the $L\sp 2$ critical nonlinear Schr\"odinger equation}, J.\ Amer.\ Math.\ Soc.\ 19 (2006), no. 1, pp. 37-90.

\bibitem[Per]{perel} G.\ Perelman, \emph{On the formation of singularities in solutions of the critical nonlinear Schr\"odinger equation}, Ann.\ Henri Poincar\'e 2 (2001), no.\ 4, pp.\ 605-673.

\bibitem[PilWay]{pillet} C.\ A.\ Pillet, C.\ E.\ Wayne, \emph{Invariant manifolds for a class of dispersive, Hamiltonian, partial differential equations}, J.\ Diff.\ Eq.\ 141 (1997), no.\ 2, pp.\ 310-326.

\bibitem[Rodsch]{rod} I.\ Rodnianski, W.\ Schlag, \emph{Time decay for solutions of Schr\"odinger equations with rough and time-dependent potentials}, Invent.\ Math.\ 155 (2004), no.\ 3, pp.\ 451-513.

\bibitem[RoScSo1]{rod2} I.\ Rodnianski, W.\ Schlag, A. Soffer, \emph{Dispersive analysis of charge transfer models}, Communications on Pure and Applied Mathematics, Volume 58, Issue 2, pp.\ 149-216.

\bibitem[RoScSo2]{rod3} I.\ Rodnianski, W.\ Schlag, A. Soffer, \emph{Asymptotic stability of N-soliton states of NLS}, preprint 2003.

\bibitem[Sch]{schlag} W.\ Schlag, \emph{Stable Manifolds for an orbitally unstable NLS}, preprint 2004, to appear in Annals of Mathematics.

\bibitem[Sch2]{schlag2} W.\ Schlag, \emph{Spectral theory and nonlinear partial differential equations: a survey}, Discrete Contin.\ Dyn.\ Syst.\ 15 (2006), no.\ 3, 703--723.

\bibitem[SofWei1]{soffer1} A.\ Soffer, M.\ I.\ Weinstein, \emph{Multichannel nonlinear scattering for nonintegrable equations}, Comm.\ Math.\ Phys.\ 133 (1990), pp.\ 119 - 146.

\bibitem[SofWei2]{soffer2} A.\ Soffer, M.\ I.\ Weinstein, \emph{Multichannel nonlinear scattering, II. The case of anisotropic potentials and data}, J.\ Diff.\ Eq.\ 98 (1992), pp.\ 376 - 390.

\bibitem[Ste]{stein} E.\ Stein, \emph{Harmonic Analysis}, Princeton University Press, Princeton, 1994.

\bibitem[SulSul]{sulem} C.\ Sulem, P.-L.\ Sulem, \emph{The Nonlinear Schr\"odinger Equation. Self-focusing and Wave Collapse}, Applied Mathematical Sciences, 139, Springer-Verlag, New York, 1999.

\bibitem[Tay]{taylor} M.\ E.\ Taylor, \emph{Tools for PDE.\ Pseudodifferential operators, paradifferential operators, and layer potentials}, Mathematical Surveys and Monographs, 81, American Mathematical Society, Providence, RI, 2000.

\bibitem[Wei1]{wein1} M.\ I.\ Weinstein, \emph{Lyapunov stability of ground states of nonlinear dispersive equations}, Comm.\ Pure Appl.\ Math.\ 39 (1986), no.\ 1, pp.\ 51-67.

\bibitem[Wei2]{wein2} M.\ I.\ Weinstein, \emph{Modulational stability of ground states of nonlinear Schr\"odinger equations}, SIAM J.\ Math.\ Anal. 16 (1985), no.\ 3, pp.\ 472-491.

\bibitem[Yaj]{yajima} K.\ Yajima, \emph{Dispersive estimate for Schr\"odinger equations with threshold resonance and eigenvalue}, preprint 2004, to appear in Comm.\ Math.\ Physics.

\end{thebibliography}
\end{document}